%%%%%%AMS-TeX file for%%%%%%%%%%%%%%%%%%%%%%%%%%%%%%%%%%%%%%%%%%%%%%
%
%   8 LECTURES ON QUANTUM GROUPS AND q-SPECIAL FUNCTIONS
%
%for the Fourth Colombian Summer School, Bogota,
%July 22-August 2, 1996,
%by Erik Koelink
%%%%%%%%%%%%%%%%%AMS TeXfile%%%%%%%%%%%%%%%%%%%%%%%%%%%%%%%%%%%%%%%%
\input amstex.tex
\documentstyle{amsppt}
\magnification=1200
\baselineskip=13pt
\hsize=6.5truein
\vsize=8.9truein
\parindent=20pt
%%%%%%%%%%%%%%%%% M a c r o s %%%%%%%%%%%%%%%%%%%%%%%%%%%%%%%%%%%%%%
%Section, Theorem and Formula Numbering%%%%%%%%%%%%%%%%%%%%%%%%%%%%%
%%%%%%%%%%%%%%%%%%%%%%%%%%%%%%%%%%%%%%%%%%%%%%%%%%%%%%%%%%%%%%%%%%%%
\countdef\sectionno=1
\countdef\eqnumber=10
\countdef\theoremno=11
\countdef\countrefno=12
\countdef\cntsubsecno=13
\sectionno=0

\def\newsection{\global\advance\sectionno by 1
                \global\cntsubsecno=0
                \the\sectionno.\ }

\def\newsubsection{\global\advance\cntsubsecno by 1
                   \global\eqnumber=1
                   \global\theoremno=1 
                   \S\the\sectionno.\the\cntsubsecno.}

\def\theoremname#1{\the\sectionno.\the\cntsubsecno.\the\theoremno
         \xdef#1{{\the\sectionno.\the\cntsubsecno.\the\theoremno}}
         \global\advance\theoremno by 1}

\def\eqname#1{\the\sectionno.\the\cntsubsecno.\the\eqnumber
         \xdef#1{{\the\sectionno.\the\cntsubsecno.\the\eqnumber}}
         \global\advance\eqnumber by 1}

\global\countrefno=1

\def\refno#1{\xdef#1{{\the\countrefno}}
             \global\advance\countrefno by 1}

\def\thmref#1{#1}

%%%%%%%%%%%%%%%%%%%%%%%%%%%%Abbreviations%%%%%%%%%%%%%%%%%%%%%%%%%%%
\def\R{{\Bbb R}}
\def\N{{\Bbb N}}
\def\C{{\Bbb C}}
\def\Z{{\Bbb Z}}
\def\T{{\Bbb T}}
\def\Zp{{\Bbb Z}_+}
\def\hZp{{1\over 2}\Zp}
\def\A{A_q(SL(2,\C))}
\def\Asu{A_q(SU(2))}
\def\U{U_q({\frak{sl}}(2,\C))}
\def\Usu{U_q({\frak{su}}(2))}
\def\Unc{U_q({\frak{su}}(1,1))}
\def\PG{{\text{Pol}}(G)}
\def\al{\alpha}
\def\be{\beta}
\def\ga{\gamma}
\def\de{\delta}
\def\De{\Delta}
\def\et{\eta}
\def\ep{\varepsilon}
\def\si{\sigma}
\def\ta{\tau}
\def\th{\theta}
\def\la{\lambda}
\def\rts{\rho_{\ta ,\si}}
\def\rti{\rho_{\ta ,\infty}}
\def\ris{\rho_{\infty ,\si}}
\def\rii{\rho_{\infty ,\infty}}
\def\vp{\varphi}
\def\Hi{\ell^2(\Zp)}

%%%%%%%%%%%%%%%%%%%%%%%%%Reference Numbering%%%%%%%%%%%%%%%%%%%%%%%%
\refno{\Abe}
\refno{\AlSa}
\refno{\AlSaC}
\refno{\Andr}
\refno{\AskeI}
\refno{\AskeIMAMS}
\refno{\AskeRS}
\refno{\AskeW}
\refno{\Bere}
\refno{\Berg}
\refno{\BergI}
\refno{\BiedL}
\refno{\Bres}
\refno{\BurbK}
\refno{\CharP}
\refno{\Chih}
\refno{\DijkK}
\refno{\DijkKGD}
\refno{\DijkN}
\refno{\Domb}
\refno{\Drin}
\refno{\DunfS}
\refno{\FlorVPLA}
\refno{\FlorVLMP}
\refno{\Flor}
\refno{\FlorK}
\refno{\GaspR}
\refno{\GranZ}
\refno{\Isma}
\refno{\IsmaW}
\refno{\Jant}
\refno{\Jose}
\refno{\Kake}
\refno{\KakeMU}
\refno{\KalnMMJMP}
\refno{\KalnMRMJM}
\refno{\KalnMJMP}
\refno{\KalnMMuJMP}
\refno{\KalnMMSIAM}
\refno{\Kass}
\refno{\KlimK}
\refno{\KoekS}
\refno{\KoelComp}
\refno{\KoelITSF}
\refno{\KoelSIAM}
\refno{\KoelDMJ}
\refno{\Koelunpub}
\refno{\KoelCJM}
\refno{\KoelIM}
\refno{\KoelAAM}
\refno{\KoelFIC}
\refno{\KoelHAlg}
\refno{\KoelVdJ}
\refno{\KoelV}
\refno{\Koorold}
\refno{\KoorIM}
\refno{\KoorOPTA}
\refno{\KoorAF}
\refno{\KoorBC}
\refno{\KoorZSE}
\refno{\KoorTrento}
\refno{\KoorMLI}
\refno{\KoorS}
\refno{\Lusz}
\refno{\Macdpp}
\refno{\Macdbook}
\refno{\MacdSB}
\refno{\Maji}
\refno{\MasuMNNSU}
\refno{\MasuMNNU}
\refno{\Mill}
\refno{\Noum}
\refno{\NoumAM}
\refno{\NoumDS}
\refno{\NoumMPJA}
\refno{\NoumMunpub}
\refno{\NoumMCMP}
\refno{\NoumMDMJ}
\refno{\NoumMLNM}
\refno{\NoumMCompM}
\refno{\NoumYM}
\refno{\Pal}
\refno{\Podl}
\refno{\RahmV}
\refno{\SchmCMP}
\refno{\Schmrep}
\refno{\StokSIAM}
\refno{\Stokrep}
\refno{\StokK}
\refno{\Sugi}
\refno{\Swee}
\refno{\Szeg}
\refno{\Temm}
\refno{\VaksKSMD}
\refno{\VaksKFAA}
\refno{\VaksS}
\refno{\VAsscK}
\refno{\VDaelBLMS}
\refno{\VDaelPAMS}
\refno{\VdJeug}
\refno{\Vile}
\refno{\VileK}
\refno{\Woro}
\refno{\Worotwee}
\refno{\Worodrie}

%%%%%%%%%%%%%%%%%%%%%%%%%%%%%%%%%%%%%%%%%%%%%%%%%%%%%%%%%%%%%%%%%%%%
%%%%%%%%%%%%%%%%%%%%%%%Beginning of Text%%%%%%%%%%%%%%%%%%%%%%%%%%%%
%%%%%%%%%%%%%%%%%%%%%%%%%%%%%%%%%%%%%%%%%%%%%%%%%%%%%%%%%%%%%%%%%%%%
\topmatter
\title 8 Lectures on\\Quantum Groups and $q$-Special
Functions\endtitle
\author Erik Koelink\endauthor
\rightheadtext{Quantum Groups and $q$-Special Functions}
\leftheadtext{Quantum Groups and $q$-Special Functions}
\address Vakgroep Wiskunde, Universiteit van Amsterdam,
Plantage Muidergracht 24,
1018 TV Amsterdam, the Netherlands\endaddress
\email koelink\@fwi.uva.nl\ \ {\sl WWW}{\rm:}\  
http://turing.fwi.uva.nl/\~{}koelink/
\endemail
\affil Report 96-10,
Universiteit van Amsterdam\endaffil
\date Updated version of August 22, 1996\enddate
\thanks Supported by the Netherlands
Organization for Scientific Research (NWO)
under project number 610.06.100. \endthanks
\abstract The lecture notes contains an introduction
to quantum groups, $q$-special functions and their
interplay. After generalities on Hopf algebras,
orthogonal polynomials and basic hypergeometric series
we work out the relation between the quantum $SU(2)$
group and the Askey-Wilson polynomials 
out in detail as the main example. As an application
we derive an addition formula for a two-parameter
subfamily of Askey-Wilson polynomials. A relation between
the Al-Salam and Chihara polynomials and the quantised
universal enveloping algebra for $su(1,1)$ is given.
Finally, more examples and other approaches
as well as some open problems are given.
\endabstract
%Contents
\toc
\widestnumber\head{8}
\head 1. Hopf algebras\endhead
\head 2. The quantum $SL(2,\C)$ group\endhead
\head 3. Orthogonal polynomials and basic hypergeometric
series\endhead
\head 4. Quantum subgroups and the Haar functional\endhead
\head 5. Askey-Wilson polynomials and generalised matrix
elements\endhead
\head 6. Addition formulas for Askey-Wilson polynomials\endhead
\head 7. Convolution theorem for Al-Salam and Chihara
polynomials\endhead
\head 8. More examples\endhead
\endtoc
\endtopmatter
\document

%%%%%%%%%%%%%%%%%%%%%%%%%%%%%%%%%%%%%%%%%%%%%%%%%%%%%%%%%%%%%%%%%%%%
%%N E W   S E C T I O N%%%%%%%%%%%%%%%%%%%%%%%%%%%%%%%%%%%%%%%%%%%%%
%%%%%%%%%%%%%%%%%%%%%%%%%%%%%%%%%%%%%%%%%%%%%%%%%%%%%%%%%%%%%%%%%%%%
\newpage
\head Introduction\endhead

In the latter half of this century it has become clear that there 
is an intimate relation between special functions of hypergeometric
type, such as Jacobi, Hahn and Krawtchouk polynomials (or other
polynomials within the Askey scheme of hypergeometric orthogonal
polynomials), Bessel functions and representation theory of
groups, notably Lie groups. See the books by Vilenkin \cite{\Vile}
and its succesor by Vilenkin and Klimyk \cite{\VileK} for
a large number of relations between special functions and
group representations, and how group properties can imply
interesting formulae for the special functions involved.

Basic, or $q$-hypergeometric series and $q$-special functions
are almost as ancient as their $q=1$ counterparts. But the
development for the $q$-special functions has been much slower
than the development for special functions. In the 1970's
the works of Andrews and Askey initiated a great number
of papers on $q$-hypergeometric series and $q$-special
functions. One of the highlights in this development
is the introduction by Askey and Wilson of a very general
four parameter set of orthogonal polynomials, nowadays
called the Askey-Wilson polynomials. However, before the
introduction of quantum groups there was not much known
about where these $q$-special functions `live' in an
analogous way as special functions live on Lie groups.
Of course, there were some isolated results such as
the interpretation of certain $q$-analogues of special
functions on finite groups of Lie type.

Since the introduction of quantum groups by Drinfeld, Jimbo
and Woronowicz in the mid-eighties, a lot of research has been
going on developing the relations between $q$-special functions
and quantum groups. A large number of such relations
are known by now. Since then a huge amount of papers
on quantum groups have appeared, and also a number of books
on quantum groups have been published. Since we are not
dealing with all kinds of aspects of quantum groups, we refer
to the books by Chari and Pressley \cite{\CharP},
Jantzen \cite{\Jant}, Joseph
\cite{\Jose}, Kassel \cite{\Kass}, Lusztig \cite{\Lusz}, Majid
\cite{\Maji} for more information. Chari and Pressley \cite{\CharP}
and also Vilenkin and Klimyk \cite{\VileK} have a chapter
on the relation between
quantum groups and $q$-special functions. For a
more physical point of view to quantum groups and related
special functions Biedenharn and Lohe \cite{\BiedL}
can be consulted.

It is the purpose of these lectures to get some feeling
for quantum groups and its relation with $q$-special functions.
I have chosen to do so by studying various aspects of one
particular simple well-known example, namely the quantum
group analogue of the compact Lie group of $2\times 2$-unitary
matrices, $SU(2)$. At the end we briefly discuss some other
examples which are related to the treatment of the
quantum $SU(2)$ group of these lectures. After reading these
course notes it must be clear that there is a very nice and
important interplay between quantum groups
and $q$-special functions.

What are quantum groups? Firstly, they are {\sl not} groups,
but they are somehow related to groups in the sense that
they are deformations of certain structures reflecting the
group properties. This is rather vague, but more precise
information can be found in the books on quantum groups
mentioned. We can think of a quantum group as a deformation
of the algebra of functions on a group. For a Lie group
we have the duality between functions on the group
and the universal enveloping algebra of the corresponding
Lie algebra. It turns out that for a large class of
quantum groups, there exist deformations of universal
enveloping algebras such that the duality survives.
These deformations are known as quantised universal enveloping
algebras, or quantum algebras, or $q$-algebras, and for
each simple Lie algebra there is a `canonical'
deformation, due to Jimbo.

In these lectures we first discuss the fundamental concept
in quantum group theory, namely Hopf algebras. Two
important examples, namely the algebra of functions
on a (finite) group and the universal enveloping algebra, are
discussed. Duality is an important concept.
In the second lecture we investigate in detail deformations
of the universal enveloping algebra $U({\frak{sl}}(2,\C))$
and the algebra of polynomials on $SL(2,\C)$.
Orthogonal polynomials and $q$-hypergeometric series
are discussed in lecture~3, and we combine the two
in a discussion of the Askey-Wilson polynomials.
In lecture~4 we describe how we can characterise
quantum subgroups of the quantum $SU(2)$ group,
and we derive an explicit expression
for the analogue of the Haar measure for
left and right invariant (with respect to such a subgroup)
functions. In lecture~5 we then show how the full four-parameter
family of Askey-Wilson polynomials can be interpreted
on the quantum $SU(2)$ group as generalised matrix elements.
Using this interpretation an addition formula for a
two-parameter family of Askey-Wilson polynomials is derived
in lecture~6. In lecture~7 we discuss how a general
convolution formula for a subclass of Askey-Wilson
polynomials can be derived from the quantised
universal enveloping algebra for ${\frak{su}}(1,1)$.
An overview of related results is finally given in
the last lecture.

The main line in the relation between quantum groups and
$q$-special functions as presented in these lectures
uses the duality between deformed
function algebras and quantised universal enveloping algebras.
It is shown that results in the theory of quantum groups can be
proved using $q$-special functions, and that identitites for
$q$-special functions can be derived using their interpretation
on quantum groups. Another approach, due to Kalnins, Miller and
coworkers, and Floreanini and Vinet, only uses the quantum algebra,
and we discuss this alternative shortly in
lecture~8. Lecture~7 contains a result which is a kind
of mixture of these two approaches.

Each lecture ends with a number of references to the literature,
which, due to the enormous amount of papers in this area,
cannot be complete. For each lecture a number of exercises
is given.
There are more published lecture notes on this subject,
notably Koornwinder \cite{\KoorOPTA}, \cite{\KoorTrento}
and Noumi \cite{\Noum}. For the general lecture~1
I have used \cite{\KoorTrento} a lot. Lectures~2, 4 and
5 follow my survey paper \cite{\KoelAAM} with a
different proof of the expression for the Haar functional
in lecture~4, which is taken from joint work with
Verding \cite{\KoelV}. Lecture~6 is based on
unpublished work \cite{\Koelunpub}. It is contained
a limiting case of the far more
complex and computational result of \cite{\KoelFIC}.
Lecture~7 is a special case of
joint work \cite{\KoelVdJ} with Van der Jeugt.
Lecture~8 gives a biased look to other cases, and discusses
some open ends.

\demo{Acknowledgement} These lecture notes were used for a course
at the IV Escuala de Verano, Universidad Nacional de Colombia 
and Universidad de los Andes, Bogot\'a, 
Colombia, July 22 -- August 2, 1996. I thank the organisors,
and in particular Jairo Charris and Ernesto Acosta, for the
invitation and their kind hospitality. I thank Tom Koornwinder and
Jasper Stokman for pointing out typos in a previous version.
\enddemo

%%%%%%%%%%%%%%%%%%%%%%%%%%%%%%%%%%%%%%%%%%%%%%%%%%%%%%%%%%%%%%%%%%%%
%%N E W   S E C T I O N%%%%%%%%%%%%%%%%%%%%%%%%%%%%%%%%%%%%%%%%%%%%%
%%%%%%%%%%%%%%%%%%%%%%%%%%%%%%%%%%%%%%%%%%%%%%%%%%%%%%%%%%%%%%%%%%%%
\newpage

\head\newsection Hopf algebras\endhead

The concept of a {\sl Hopf algebra} is fundamental for
the theory of quantum groups.
In the first lecture we study this concept
in some detail and we treat some important examples.

The ground field is $\C$, although everything goes through
when working over a commutative ring with unit. The tensor product
$V\otimes W$ of two linear spaces is the {\sl algebraic tensor 
product}, which means that elements of $V\otimes W$ consist of 
finite linear combinations of the form $v\otimes w$, 
$v\in V$, $w\in W$.

%%%%%%%%%%%%%%%%%%%%%%%%%%%%%%%%%%%%%%%%%%%%%%%%%%%%%%%%%%%%%%%%%%%%
%%N E W   S U B S E C T I O N%%%%%%%%%%%%%%%%%%%%%%%%%%%%%%%%%%%%%%%
%%%%%%%%%%%%%%%%%%%%%%%%%%%%%%%%%%%%%%%%%%%%%%%%%%%%%%%%%%%%%%%%%%%%
\subhead\newsubsection 
Algebras, bi-algebras and Hopf algebras
\endsubhead
Recall that an {\sl algebra}, or better, an associative algebra with
unit, is a linear space $A$ with a bilinear mapping
$A\times A\to A$, $(a,b)\mapsto ab$,
called {\sl multiplication},
and a distinguished non-zero element $1\in A$, called the 
{\sl unit}, such that $a(bc)=(ab)c$ and $1a=a=a1$ for all 
$a,b,c\in A$. This leads to two mappings, 
$m\colon A\otimes A\to A$, also called
multiplication, and $\et\colon \C\to A$, also called unit, defined
by $m(a\otimes b) = ab$ and $\et(z)=z1$. Then we can rephrase the
associatitvity and unit in terms of the following commuting 
diagrams;
$$
\CD
A\otimes A\otimes A @>m\otimes id>> A\otimes A
@.\qquad @.\C\otimes A @>\et\otimes id>> A\otimes A
@<id\otimes\et << A\otimes \C\\
@Vid\otimes mVV @VmVV
@.  @V\cong VV @VmVV @V\cong VV\\
A\otimes A @>m>> A
@.\qquad @. A @>id>> A @<id<< A
\endCD
\tag\eqname{\vgldiagramalgebra}
$$
where we use $\C\otimes A\cong A \cong A\otimes \C$ by
identifying $z\otimes a$ and $a\otimes z$ with $za$ for $z\in\C$ 
and $a\in A$. An algebra homomorphism means a unital algebra 
homomorphism, i.e. mapping unit onto unit.

For an algebra $A$ the tensor product $A\otimes A$ is again an 
algebra with multiplication 
$(a\otimes b)(c\otimes d)=ac\otimes bd$ and
unit $1\otimes 1$. Note that
$$
m_{A\otimes A} \colon A^{\otimes 4} \to A^{\otimes 2}, \qquad
m_{A\otimes A} = (m_A \otimes m_A)\circ (id\otimes \si \otimes id)
\tag\eqname{\vgldefmAtensorA}
$$
where $\si\colon A\otimes A\to A\otimes A$ is the flip automorphism,
$\si(a\otimes b)=b\otimes a$. Finally, note that commutativity
of $A$ is equivalent to the condition $\si\circ m = m$.

\proclaim{Definition \theoremname{\defcoalgebra}} A {\bf
coalgebra}, or better, a coassociative coalgebra with counit,
is a linear space $A$ with a linear mapping 
$\De\colon A\to A\otimes A$,
called the comultiplication, and a non-zero linear mapping
$\ep\colon A\to\C$, called the counit, such that
the following diagram is commutative;
$$
\CD
A\otimes A\otimes A @<\De\otimes id<< A\otimes A
@.\qquad @.\C\otimes A @<\ep\otimes id<< A\otimes A
@>id\otimes\ep >> A\otimes \C\\
@AAid\otimes \De A @AA\De A
@.  @V\cong VV @AA\De A @V\cong VV\\
A\otimes A @<\De << A
@.\qquad @. A @<id<< A @>id>> A
\endCD
\tag\eqname{\vgldiagramcoalgebra}
$$
\endproclaim

The commutative diagram of \thetag{\vgldiagramcoalgebra}
is obtained from \thetag{\vgldiagramalgebra} by reversing
arrows.

We say that the coalgebra is cocommutative if $\si\circ\De=\De$.

\proclaim{Definition \theoremname{\defbialgebra}} A {\bf bialgebra}
is an algebra $A$, such that $A$ is also a coalgebra and the
comultiplication $\De$ and counit $\ep$ are algebra homomorphisms.
\endproclaim

\demo{Remark \theoremname{\remdefbialgebra}} An equivalent 
definition of a bialgebra can be obtained replacing the condition 
that $\De$ and $\ep$ are algebra homomorphism by the condition that
$m$ and $\et$ are coalgebra homomorphisms.
\enddemo

\proclaim{Definition \theoremname{\defHopfalgebra}} A {\bf Hopf 
algebra} is a bialgebra $A$ with a linear mapping $S\colon A\to A$,
the {\bf antipode}, such that the following diagram is commutative;
$$
\CD
A\otimes A @>S\otimes id>> A\otimes A @<id \otimes S << 
A\otimes A \\
@AA\De A @VmVV @AA\De A \\
A @>\et\circ\ep >> A @< \et\circ\ep << A
\endCD
\tag\eqname{\vgldiagramS}
$$
A {\bf Hopf algebra morphism} $\phi\colon A\to B$ of two Hopf 
algebras $A$ and $B$ is an algebra morhism $\phi\colon A \to B$ 
such that $\ep_B\circ\phi=\ep_A$, 
$\De_B\circ\phi=\phi\otimes\phi\circ\De_A$
and $S_B\circ\phi=\phi\circ S_A$.
\endproclaim

\proclaim{Proposition \theoremname{\proppropsofHopfalgebra}}
{\rm (i)} If $A$ is a bialgebra and $S$ an antipode making $A$ into
a Hopf algebra, then $S$ is unique.

{\rm (ii)} Let $A$ be a Hopf algebra, then the antipode $S$ is
unital, counital, antimultiplicative and anticomultiplicative. Or,
with $\si$ the flip automorphism,
$$
S(1)=1,\quad \ep\circ S = \ep, \quad
S\circ m = m\circ \si\circ (S\otimes S),\quad
\De\circ S = \si\circ(S\otimes S)\circ \De.
$$
\endproclaim

\demo{Proof} (i) Let $F$ and $G$ be linear mappings of $A$ into
itself, then we define the convolution product $F\ast G$ by
$F\ast G= m\circ(F\otimes G)\circ \De$. This convolution product
is associative, which follows from the associativity of $m$ and
the coassociativity of $\De$. Moreover, $\et\circ\ep\colon A\to
A$ is the unit for the convolution product.
So the endomorphism algebra of
$A$, $End(A)$, becomes algebra.
Now assume that $A$ is a Hopf algebra with antipode $S$, then
$S\ast id = \et\circ\ep=id\ast S$. Or, $S$ is a two-sided
inverse of the identity mapping in $End(A)$ with respect to
the convolution product and thus unique.

(ii) Use \thetag{\vgldiagramS} to $1\in A$ to find $S(1)=1$.
To ease notation we introduce
$$
\De(a) = \sum_{(a)} a_{(1)}\otimes a_{(2)}, \qquad
(id\otimes \De)\De(a) = \sum_{(a)} a_{(1)}\otimes a_{(2)} 
\otimes a_{(3)},
\tag\eqname{\vglSweedlernot}
$$
which is well-defined by \thetag{\vgldiagramcoalgebra}.
Then by \thetag {\vgldiagramcoalgebra} we have
$\sum_{(a)}\ep(a_{(1)})a_{(2)}=a$. Apply $S$ and next $\ep$
to see that
$$
\ep\bigl(S(a)\bigr) = \ep\otimes \ep\bigl( (id\otimes S)
\De(a)\bigr) =
\ep\circ m\circ(id\otimes S)\circ \De (a) = 
\ep\circ\et\circ\ep (a) = \ep(a).
$$
Next
$$
\align
S(b)S(a)&= \sum_{(a),(b)} S(b_{(1)})S(a_{(1)}) 
\ep(a_{(2)}b_{(2)}) \\
&= \sum_{(a),(b)}  S(b_{(1)})S(a_{(1)}) a_{(2)}b_{(2)} 
S(a_{(3)}b_{(3)}) \\
&= \sum_{(a),(b)}  \ep(a_{(1)})\ep(b_{(1)}) S(a_{(2)}b_{(2)})
=S(ab).
\endalign
$$
The last statement is left as an exercise.\qed
\enddemo

\demo{Example \theoremname{\exaHopfalgebraone}} Let $G$ be a 
finite group and $A=C(G)$, the space of (continuous) 
complex-valued functions
on $G$. Then $A$ is a commutative algebra under pointwise 
multiplication and its unit is the constant function equal to $1$. 
Since $G$ is finite we have $A\otimes A\cong C(G\times G)$, and 
then $m(F)(g)=F(g,g)$, $\et(z)(g)=z$. Moreover, $A$ is a bialgebra 
if we define the comultiplication and counit by
$$
\align
&\De(f)(g,h) = f(gh), \qquad f\in A, \ g,h\in G,\\
&\ep(f) = f(e), \qquad f\in A, \text{$e\in G$ unit of the group}.
\endalign
$$
We define $S(f)(g)=f(g^{-1})$ and then $A$ is a Hopf algebra
by a straightforward, but instructive, check. Note that $A$
is a commutative Hopf algebra with $S^2=1$.

The important observation is that the group structure 
--multiplication, unit and inverse-- is stored in the Hopf 
algebra structure --comultiplication, counit and antipode-- 
and so instead of studying $G$ we can study $C(G)$.

This works nicely for a finite group. Now suppose that $G$ is an
algebraic subgroup of $SL(n,\C)$, the group of $n\times n$-matrices
with complex entries and determinant one, so $G$ is some
complex matrix group. Then we can take $A=\PG$
consisting of complex-valued functions in $g$, such that considered
as functions of its matrix entries $g_{ij}$ they are polynomials.
Then ${\text{Pol}}(G\times G)\cong \PG\otimes \PG$.
Let $t_{ij}\in A$ be defined by $t_{ij}(g)=g_{ij}$, then the
$t_{ij}$ generate $A$ and
$\De(t_{ij})=\sum_{k=1}^nt_{ik}\otimes t_{kj}$,
$\ep(t_{ij})=\de_{ij}$, $S(t_{ij})=T_{ji}$, where $T_{ji}$ is the
cofactor of the $(ji)$-th entry in the matrix
$(t_{ij})_{1\leq i,j\leq n}$.
\enddemo

\demo{Example \theoremname{\exaHopfalgebratwo}} Let $\frak g$ be
a complex Lie algebra, i.e. a vector space over $\C$ equipped
with a Lie bracket $[\cdot,\cdot]\colon {\frak g}\times{\frak g}
\to {\frak g}$, which is a bilinear mapping
such that $[X,Y]=-[Y,X]$ and
$[X,[Y,Z]]+[Y,[Z,X]]+[Z,[X,Y]]=0$ (Jacobi identity). An example
of a Lie algebra is ${\frak g}={\frak{sl}}(n,\C)$, the space of
$n\times n$-matrices with complex coefficients with zero trace.
The Lie bracket is $[X,Y]=XY-YX$, where the product on the right
hand side is matrix multiplication.

Let $A$ be the universal enveloping algebra of $\frak g$,
$A=U({\frak g})$. This is a unital algebra generated by
$X$, $X\in{\frak g}$, with relations $XY-YX=[X,Y]$ for
$X,Y,\in{\frak g}$. Then we define
$$
\De(X) = 1\otimes X+X\otimes 1, \ \ep(X)=0, \ S(X)=-X,
\qquad X\in{\frak g}
$$
and extend $\De$ and $\ep$ to $U({\frak g})$ as algebra 
homomorphisms and $S$ as an antihomomorphism.
Then $U({\frak g})$ is a Hopf algebra. Note that it is
a cocommutative Hopf algebra, but that it is not commutative
unless ${\frak g}$ is abelian, i.e. $[X,Y]=0$ for all
$X,Y\in{\frak g}$.
\enddemo

%%%%%%%%%%%%%%%%%%%%%%%%%%%%%%%%%%%%%%%%%%%%%%%%%%%%%%%%%%%%%%%%%%%%
%%N E W   S U B S E C T I O N%%%%%%%%%%%%%%%%%%%%%%%%%%%%%%%%%%%%%%%
%%%%%%%%%%%%%%%%%%%%%%%%%%%%%%%%%%%%%%%%%%%%%%%%%%%%%%%%%%%%%%%%%%%%
\subhead\newsubsection
Duality and Hopf $\ast$-algebras
\endsubhead
A pairing between two vector spaces $A$ and $U$ is a bilinear 
mapping $A\times U\to\C$, $(a,u)\mapsto \langle a,u\rangle$. 
We say that
the pairing is non-degenerate, or better, doubly non-degenerate,
if $\langle a,u\rangle=0$ for all $u\in U$ implies $a=0$ and
$\langle a,u\rangle =0$ for all $a\in A$ implies $u=0$. Such a 
pairing can be
extended to a pairing of $A\otimes A$ and $U\otimes U$ by
$\langle a\otimes b, u\otimes v\rangle =
\langle a,u\rangle \langle b,v\rangle$.

\proclaim{Definition \theoremname{\defdualHopfalgebras}} Two
Hopf algebras $A$ and $U$ are said to be in duality if there
exists a non-degenerate pairing
$\langle\cdot,\cdot\rangle \colon A\times U\to\C$ such that
$$
\gather
\langle a\otimes b, \De(u)\rangle = \langle ab, u\rangle, \quad
\langle \De(a),u\otimes v\rangle = \langle a,uv\rangle, \\
\langle 1, u\rangle = \ep(u),\quad \langle a, 1\rangle = 
\ep(a),\qquad \langle S(a),u\rangle = \langle a,S(u)\rangle.
\endgather
$$
\endproclaim

\demo{Example \theoremname{\exampairingGandUg}}
Recall that if $G$ is an (algebraic) subgroup of $SL(n,\C)$ and
its Lie algebra ${\frak g}$ is a subalgebra of ${\frak{sl}}(n,\C)$,
then we have a natural pairing of $\PG$ and ${\frak g}$ by
$$
\langle X, f\rangle = {d\over{dt}}\Big|_{t=0} f(\exp tX),
\qquad X\in{\frak g}, \ f\in \PG.
$$ From Definition \thmref{\defdualHopfalgebras} and
Examples \thmref{\exaHopfalgebraone} and \thmref{\exaHopfalgebratwo}
we find the familiar expression
$$
\langle X_1\ldots X_k, f\rangle = {{\partial^k}\over
{\partial t_1\ldots\partial t_k}}\Big|_{t_1=\ldots=t_k=0}
f\bigl( (\exp t_1X_1)\ldots(\exp t_kX_k)\bigr)
$$
for $X_1,\dots,X_k\in{\frak g}$, $f\in \PG$.

Using this pairing we can define a left and right action of
${\frak g}$ on $\PG$ by left and right invariant first order 
differential operators; for $g\in G$, $X\in{\frak g}$, $f\in\PG$,
$$
X.f(g) = {d\over{dt}}\Big|_{t=0} f(g\exp tX), \qquad
f.X(g) = {d\over{dt}}\Big|_{t=0} f((\exp tX) g)
$$
and we extend this to $U({\frak g})$ by $(XY).f=X.(Y.f)$ and
$f.(XY)=(f.X).Y$. This gives the action of $U({\frak g})$
on smooth functions on $G$ by left or right invariant differential
operators. This can be done for arbitrary Hopf algebras in duality.
\enddemo

\proclaim{Proposition \theoremname{\propactionUonA}} Let $A$ and $U$
be two Hopf algebras in duality, then for $u\in U$ and $a\in A$
$$
u.a = (id \otimes \langle\cdot,u\rangle)\circ\De(a), \qquad
a.u = (\langle\cdot,u\rangle\otimes id)\circ\De(a)
$$
define a left and right action of $U$ on $A$.
\endproclaim

\demo{Proof} Using \thetag{\vglSweedlernot} we see that
$$
(vu).a = \sum_{(a)} \langle a_{(2)},vu\rangle a_{(1)} =
\sum_{(a)}a_{(1)} \langle a_{(2)},v\rangle\langle a_{(3)},u\rangle
=v.(u.a),
$$
and and similarly for the right action.
\qed\enddemo

Recall that an algebra $A$ is $\ast$-algebra if there exists
an antilinear antimultiplicative involution $a\mapsto a^\ast$.
So $(\la a+\mu b)^\ast=\bar\la a^\ast+\bar \mu b^\ast$,
$(ab)^\ast=b^\ast a^\ast$, $(a^\ast)^\ast=a$
for $a,b\in A$, $\la,\mu\in\C$. If $\phi\colon A\to B$
is an algebra homomorphism of two $\ast$-algebras $A$ and $B$,
then $\phi$ is a $\ast$-homomorphism if $\phi(a^\ast)=
\phi(a)^\ast$.

\proclaim{Definition \theoremname{\defHopfstaralg}} A
{\bf Hopf $\ast$-algebra} is a Hopf algebra $A$, such that
$A$ is a $\ast$-algebra and $\De$ and $\ep$ are
$\ast$-homomorphisms. A {\bf Hopf $\ast$-algebra morphism}
$\phi\colon A\to B$ of two Hopf $\ast$-algebras $A$ and $B$
is a Hopf algebra morphism such that
$\phi(a^\ast)=\bigl(\phi(a)\bigr)^\ast$.
\endproclaim

There is no requirement on the relation between the antipode
$S$ and $\ast$ in Definition \thmref{\defHopfstaralg}.

\proclaim{Proposition \theoremname{\propSandstar}} Let $A$ be a
Hopf $\ast$-algebra, then $(S\circ\ast)^2=id$. In particular,
$S$ is invertible.
\endproclaim

\demo{Proof} For an algebra $A$ we define the opposite algebra
$A_{\text{opp}}$ as the algebra with the same vector space
structure and multiplication and unit defined by
$$
m_{\text{opp}}=m\circ\si \colon A_{\text{opp}} \otimes 
A_{\text{opp}} \to A_{\text{opp}}, \qquad \et_{\text{opp}}=
\et\colon A_{\text{opp}} \to \C.
$$
Assume that $A$ is an Hopf algebra with invertible antipode, then
$A_{\text{opp}}$ is a Hopf algebra with unchanged comultiplication
and counit and antipode $S^{-1}$.

Let $A$ be a Hopf $\ast$-algebra.
Use \thetag{\vgldiagramS} on $a^\ast$, and apply again the 
$\ast$-operator
to get, using the notation \thetag{\vglSweedlernot},
$$
\sum_{(a)} a_{(2)}\bigl( S(a_{(1)}^\ast)\bigr)^\ast = \ep(a)1=
\sum_{(a)} \bigl( S(a_{(a)}^\ast)\bigr)^\ast a_{(2)}
$$
or $\ast\circ S\circ\ast$ is the antipode for $A_{\text{opp}}$.
Using Proposition \thmref{\proppropsofHopfalgebra}(i) and the
previous paragraph it follows that $S^{-1}=\ast\circ S\circ\ast$.
\qed\enddemo

\demo{Example \theoremname{\exaHopfalgebrathree}}
Let $G\subset SL(n,\C)$ be an algebraic group as
in Example \thmref{\exaHopfalgebraone}, and let $G_0$ be
a real connected group which is a real form of $G$.
E.g. if $G=SL(n,\C)$ we can take $G_0$ equal to
$SL(n,\R)$ ($n\times n$ matrices with real entries and
determinant one), which are the fixed points of
the involution which is entry-wise complex conjugation,
or to $SU(n)$ ($n\times n$ unitary
matrices with determinant one), which are the fixed points of
the involution which takes adjoints.
A polynomial $p$ on $G$
is then completely determined by its restriction to
the real form $G_0$. Suppose that for every $p\in\PG$
there exists a polynomial $p^\ast\in\PG$ such that
$p^\ast(g_0)=\overline{p(g_0)}$ for all $g_0\in G_0$,
then $\PG$ is a Hopf $\ast$-algebra. Conversely,
if $\PG$ is a Hopf $\ast$-algebra, then
$G_0=\{ g\in G\mid p^\ast(g)=\overline{p(g)}\, \forall p\in\PG\}$
defines a real form of $G$.

So $\PG$ as a Hopf algebra carries the properties of a complex
group, and $\PG$ considered as a Hopf $\ast$-algebra carries
the properties of a real group.
\enddemo

\proclaim{Definition \theoremname{\defHopfastalgduality}}
Two Hopf $\ast$-algebras $A$ and $U$ are in duality (as Hopf 
$\ast$-algebras), if they are in duality as Hopf algebras and 
$\langle a^\ast,u\rangle =\overline{\langle a, (S(u))^\ast\rangle}$.
\endproclaim

\demo{Example \theoremname{\exaHopfalgebrafour}} From
Example \thetag{\exaHopfalgebratwo} we know that $U({\frak g})$
is a Hopf algebra, which is in duality with the Hopf algebra
$\PG$ if the Lie algebra of $G$ is ${\frak g}$. Suppose
$\PG$ is a Hopf $\ast$-algebra and let $G_0$ be the corresponding
real form of $G$ and ${\frak g}_0$ the corresponding real form
of ${\frak g}={\frak g}_0+i{\frak g}_0$. Then for $X\in {\frak g}_0$
we have
$$
\langle X, p\rangle = \overline{\langle X,(S(p))^{\ast}\rangle}
=\overline{ {d\over{dt}}\Big\vert_{t=0} {\overline{p(\exp (-tX))}}}
= -\langle X, p\rangle,
$$
so $X^\ast=-X$ for $X\in{\frak g}_0$. E.g. in case $G=SL(n,\C)$ and
$G_0=SU(n)$, then the $\ast$-operator on the Lie algebra level
is taking adjoints.
\enddemo

%%%%%%%%%%%%%%%%%%%%%%%%%%%%%%%%%%%%%%%%%%%%%%%%%%%%%%%%%%%%%%%%%%%%
%%N E W   S U B S E C T I O N%%%%%%%%%%%%%%%%%%%%%%%%%%%%%%%%%%%%%%%
%%%%%%%%%%%%%%%%%%%%%%%%%%%%%%%%%%%%%%%%%%%%%%%%%%%%%%%%%%%%%%%%%%%%
\subhead\newsubsection
Invariant functional
\endsubhead
Let $A$ be a Hopf algebra and assume that a linear functional
$h\colon A\to \C$ satisfies
$$
(id\otimes h)\circ\De = \et\circ h = (h\otimes id)\circ\De
\tag\eqname{\vgldefinvariantint}
$$
then we say that $h$ is an {\sl invariant functional}, or $h$ is a
{\sl Haar functional}.
If only the first, respectively last,
equation of \thetag{\vgldefinvariantint} holds,
then we say that $h$ is a left, respectively right, invariant 
functional. If it exists we normalise $h$ by $h(1)=1$.

\proclaim{Theorem \theoremname{\thmHaarfunctionalgen}} If the
Haar functional exists on a Hopf algebra $A$,
then it is unique up to a scalar multiple.
\endproclaim

\demo{Example \theoremname{\exaHopfalgebratfive}} Let $G$ be an 
algebraic subgroup of $SL(n,\C)$, then the left Haar measure 
$d\mu$ gives a left invariant functional on $\PG$ by $h(p) = 
\int_G p(g)\, d\mu(g)$, assuming that the polynomials are
integrable on $G$ with respect to the Haar measure.
Here we have used that
$$
(id \otimes h)\De(p)(g') = \int_G p(g'g)\, d\mu(g) = \int_G p(g)\,
d\mu(g) = h(p).
$$
If $d\mu$ is also right invariant, i.e. if $G$ is unimodular, then
$h$ is also a right invariant functional.
\enddemo

For the special case of the quantised function algebra on the 
compact group $SU(2)$ we give a proof of the existence of the 
invariant integral.

%%%%%%%%%%%%%%%%%%%%%%%%%%%%%%%%%%%%%%%%%%%%%%%%%%%%%%%%%%%%%%%%%%%%
%%N E W   S U B S E C T I O N%%%%%%%%%%%%%%%%%%%%%%%%%%%%%%%%%%%%%%%
%%%%%%%%%%%%%%%%%%%%%%%%%%%%%%%%%%%%%%%%%%%%%%%%%%%%%%%%%%%%%%%%%%%%
%NOTES AND REFERENCES%%%%%%%%%%%%%%%%%%%%%%%%%%%%%%%%%%%%%%%%%%%%%%%
\subhead Notes and references
\endsubhead
The standard references for the theory of Hopf algebras before
the introduction of quantum groups are Abe \cite{\Abe} and
Sweedler \cite{\Swee}. Since its introduction a wealth of
papers has appeared, and we have used Chari and Pressley
\cite{\CharP, Ch.~4} and especially Koornwinder \cite{\KoorTrento}.
Important concepts within the theory of Hopf algebras related
to quantum groups are quasitriangular Hopf algebras  and the
quantum double construction, both due to
Drinfeld \cite{\Drin}, see also \cite{\CharP}, \cite{\Maji}.
Also, we haven't mentioned corepresentations of a coalgebra,
see Exercise~1.6.

For Hopf $\ast$-algebras in duality we refer to Van Daele 
\cite{\VDaelBLMS}. The invariant functional is an important tool 
in the harmonic analysis on quantum groups and proofs of Theorem 
\thmref{\thmHaarfunctionalgen}
can be found for the $C^\ast$-algebra setting in Woronowicz's
influential fundamental paper \cite{\Woro}, see also Van Daele
\cite{\VDaelPAMS} for an up-to-date version. 
A purely algebraic proof
can be found in Dijkhuizen and Koornwinder \cite{\DijkK}, see als
\cite{\KoorTrento}.

%%%%%%%%%%%%%%%%%%%%%%%%%%%%%%%%%%%%%%%%%%%%%%%%%%%%%%%%%%%%%%%%%%%%
%%N E W   S U B S E C T I O N%%%%%%%%%%%%%%%%%%%%%%%%%%%%%%%%%%%%%%%
%%%%%%%%%%%%%%%%%%%%%%%%%%%%%%%%%%%%%%%%%%%%%%%%%%%%%%%%%%%%%%%%%%%%
%EXERCISES%%%%%%%%%%%%%%%%%%%%%%%%%%%%%%%%%%%%%%%%%%%%%%%%%%%%%%%%%%
\subhead Exercises
\endsubhead

\item{\the\sectionno.1} In \thetag{\vgldefmAtensorA} the 
multiplication $m_{A\otimes A}$ is defined using the flip 
automorphism $\si$. Now let
$\Psi\colon A\otimes A\to A\otimes A$ be a linear map and define
$m_{A\otimes A}=m\otimes m \circ id\otimes \Psi\otimes id$.
Suppose that $\Psi(a\otimes 1)=1\otimes a$, 
$\Psi(1\otimes a)=a\otimes 1$ for all $a\in A$, and that
$$
\align
\Psi\circ(m\otimes id)&=
(id\otimes m)\circ(\Psi\otimes id)\circ(id\otimes\Psi),\\
\Psi\circ(id\otimes m)&=
(m\otimes id)\circ(id\otimes \Psi)\circ(\Psi\otimes id).
\endalign
$$
Prove that $A\otimes A$ is an associative algebra with unit 
$1\otimes 1$. The map $\Psi$ is called the {\sl braiding}.

\item{\the\sectionno.2} Define the notion of coalgebra 
homomorphism and prove Remark \thmref{\remdefbialgebra}.

\item{\the\sectionno.3} Prove the last statement of Proposition
\thmref{\proppropsofHopfalgebra}(ii).

\item{\the\sectionno.4} Check that the examples in
Example \thmref{\exaHopfalgebraone} and \thmref{\exaHopfalgebratwo}
are Hopf algebras. (First prove that in Example
\thmref{\exaHopfalgebratwo} the mappings $\De$, $\ep$ and $S$ are
well-defined as follows; $U({\frak g}) = T({\frak g})/I$ where
$T({\frak g})$ is the tensor algebra for ${\frak g}$ and $I$ is
the ideal in $T({\frak g})$ generated by 
$X\otimes Y-Y\otimes X -[X,Y]$
for all $X,Y\in{\frak g}$. Then show that
$\De(I)\subset T({\frak g})\otimes I + I\otimes T({\frak g})$,
$\ep(I)=0$, $S(I)\subset I$. This means that $I\subset T({\frak g})$
is a Hopf ideal.)

\item{\the\sectionno.5} Prove that $A_{\text{opp}}$ is a Hopf 
algebra, see proof of Proposition \thmref{\propSandstar}. 
Similarly, if $A$ is a coalgebra
we define $A^{\text{opp}}$ as the coalgebra with comultiplication
$\De^{\text{opp}}=\si\circ\De$ and $\ep^{\text{opp}}=\ep$. 
Suppose that $A$ is a Hopf algebra with invertible antipode $S$,
prove that $A^{\text{opp}}$ is a Hopf algebra
with unchanged multiplication and unit and antipode $S^{-1}$.
Prove also that $A_{\text{opp}}^{\text{opp}}=
(A_{\text{opp}})^{\text{opp}}=(A^{\text{opp}})_{\text{opp}}$
is a Hopf algebra with antipode $S$.

\item{\the\sectionno.6} A (left) representation $\pi$ of an
algebra $A$ in the
linear space $V$ is a bilinear mapping $\pi\colon A\otimes V\to V$,
$\pi\colon a\otimes v\mapsto a\cdot v$, such that
$(ab)\cdot v=a\cdot(b\cdot v)$ and $1\cdot v=v$. Rephrase this in 
terms of commutative diagrams using the multiplication $m$ and 
unit $\et$. Define the notion of
a corepresentation of a coalgebra by reversing arrows.

%%%%%%%%%%%%%%%%%%%%%%%%%%%%%%%%%%%%%%%%%%%%%%%%%%%%%%%%%%%%%%%%%%%%
%%N E W   S E C T I O N%%%%%%%%%%%%%%%%%%%%%%%%%%%%%%%%%%%%%%%%%%%%%
%%%%%%%%%%%%%%%%%%%%%%%%%%%%%%%%%%%%%%%%%%%%%%%%%%%%%%%%%%%%%%%%%%%%
\newpage

\head\newsection The quantum $SL(2,\C)$ group\endhead

In this lecture we consider the fundamental example of a
non-commutative, non-co\-com\-mu\-tative Hopf algebra; the so-called
quantised universal enveloping algebra for ${\frak{sl}}(2,\C)$.

%%%%%%%%%%%%%%%%%%%%%%%%%%%%%%%%%%%%%%%%%%%%%%%%%%%%%%%%%%%%%%%%%%%%
%%N E W   S U B S E C T I O N%%%%%%%%%%%%%%%%%%%%%%%%%%%%%%%%%%%%%%%
%%%%%%%%%%%%%%%%%%%%%%%%%%%%%%%%%%%%%%%%%%%%%%%%%%%%%%%%%%%%%%%%%%%%
\subhead\newsubsection 
The Hopf algebra $\U$
\endsubhead
Let $\U$ be the complex unital
associative algebra generated by $A$, $B$, $C$, $D$ subject to 
the relations
$$
AD=1=DA, \quad AB=qBA,\quad AC=q^{-1}CA,\quad
BC-CB = {{A^2-D^2}\over{q-q^{-1}}}.
\tag\eqname{\vglcommrelUq}
$$
On the level of generators we define the comultiplication, counit 
and antipode by
$$
\gather
{\align
&\De(A)=A\otimes A,\quad \De(B)=A\otimes B+B\otimes D,\\
&\De(C) = A\otimes C+C\otimes D, \quad \De(D)=D\otimes D,
\endalign} \\
\ep(A)=\ep(D)=1,\qquad
\ep(C)=\ep(B)=0,
\tag\eqname{\vglHopfstructureUq}\\
S(A)=D,\quad  S(B)=-q^{-1}B,\quad S(C)=-qC,\quad S(D)=A.
\endgather
$$
Here $q$ is thought of as a deformation parameter, and at first
we take $q\in\C\backslash\{-1,0,1\}$. (It is also possible to view
$\U$ as an algebra over $\C(q)$.)
The case $q\to 1$ is considered in a moment. However, we will always
assume that $q$ is not a root of unity, i.e. $q^m\not= 1$ for all
$m\in\Zp$.
Observe that $S$ is invertible, but $S^2\not= 1$ since $q^2\not=1$.

\proclaim{Proposition \theoremname{\propUqisHopfalgebra}} Define
$\De$ and $\ep$ on $\U$ by \thetag{\vglHopfstructureUq} as
(unital) algebra homomorphisms and $S$ by 
\thetag{\vglHopfstructureUq} as
(unital) anti-algebra homomorphisms, then $\U$ is a Hopf algebra.
\endproclaim

\demo{Proof} We have to check that $\De$, $\ep$ and $S$ are 
well-defined
and that they satisfy the axioms of a Hopf algebra. These are
straighforward computations. E.g.
$$
\De(AB)= A^2 \otimes AB + AB\otimes 1 = q(A^2\otimes BA + 
BA\otimes 1)=q\De(BA)
$$
and
$$
m\circ (id\otimes S)\circ \De(B) = AS(B) + BS(D)= 
-q^{-1}AB+BA=0=\ep(B).
$$
Continuing in this way proves the proposition.
\qed\enddemo

The element
$$
\Omega = {{q^{-1}A^2+qD^2-2}\over{(q^{-1}-q)^2}} + BC
 = {{qA^2+q^{-1}D^2-2}\over{(q^{-1}-q)^2}} + CB
\tag\eqname{\vgldefCasimir}
$$
is the Casimir element of the quantised universal enveloping algebra
$\U$. $\Omega$ belongs to the centre of $\U$.

To justify the name for this Hopf algebra, we replace $A$ by
$\exp( (q-1)H/2)$, and hence $D$ by $\exp( (1-q)H/2)$ and we
let $q\uparrow 1$. Then we can deduce from \thetag{\vglcommrelUq}
that in the limit we get
$$
[H,B]=2B, \qquad [H,C]=-2C, \qquad [B,C]=H.
$$
To see the first relation, use $A=\exp( (q-1)H/2)$ in $AB=qBA$
 to get
$$
\gather
B+ {1\over 2}(q-1)HB = qB +{1\over 2}q(q-1)BH + {\Cal O}
\bigl( (q-1)^2\bigr)\\
\Longrightarrow {1\over 2}(q-1)\bigl( HB-qBH\bigr) = (q-1)B +
{\Cal O}\bigl( (q-1)^2\bigr).
\endgather
$$
Divide both sides by $q-1$ and let $q\to 1$ to get the first 
relation. The other relations are obtained similarly.
With
$$
H=\pmatrix 1&0 \\ 0&-1\endpmatrix, \qquad
B=\pmatrix 0&1 \\ 0&0 \endpmatrix , \qquad
C=\pmatrix 0&0 \\ 1&0 \endpmatrix
$$
we see that $\{ H,B,C\}$ forms a basis for the three-dimensional
Lie algebra ${\frak{sl}}(2,\C)$.
Moreover, we see that \thetag{\vglHopfstructureUq} tends to the
standard Hopf algebra structure on $U({\frak{sl}}(2,\C)$ as in
Example \thmref{\exaHopfalgebratwo} as $q\to 1$.

For the universal enveloping algebra the Poincar\'e-Birkhoff-Witt
theorem gives a basis for the underlying linear space. Here we have
a similar result, but the proof is not a straightforward 
generalisation of the PBW-theorem.

\proclaim{Lemma \theoremname{\lemPBWforU}} A linear basis for
$\U$ is given by $D^lC^kB^m$ for $k,m\in\Zp$, $l\in\Z$ with
the convention $D^{-l}=A^l$ for $l\in\Zp$.
\endproclaim

%%%%%%%%%%%%%%%%%%%%%%%%%%%%%%%%%%%%%%%%%%%%%%%%%%%%%%%%%%%%%%%%%%%%
%%N E W   S U B S E C T I O N%%%%%%%%%%%%%%%%%%%%%%%%%%%%%%%%%%%%%%%
%%%%%%%%%%%%%%%%%%%%%%%%%%%%%%%%%%%%%%%%%%%%%%%%%%%%%%%%%%%%%%%%%%%%
\subhead\newsubsection 
Finite dimensional representations of $\U$
\endsubhead
We first concentrate on the algebra structure of $\U$. Let
$H$ be a finite dimensional vector space, then ${\Cal B}(H)$
is the (unital associative) algebra of linear operators of $H$ into
itself. A representation $t$ of $\U$ in $H$ is an algebra
morphism $t\colon\U\to {\Cal B}(H)$.
A linear subspace $V\subset H$ is called
invariant if $t(X)V\subset V$ for all $X\in\U$, and then the 
restriction
$t\vert_V\colon \U\to {\Cal B}(V)$, $t\vert_V(X)=t(X)\vert_V$ is
a subrepresentation of $t$. Obviously, $\{ 0\}$ and $H$ are
invariant subspaces of $H$, and if there are no more invariant
subspaces we say that the representation $t$ in $H$ is irreducible.
Two irreducible representations of $\U$, say 
$t\colon\U\to{\Cal B}(H)$
and $s\colon\U\to{\Cal B}(V)$, are said to be equivalent if there
exists a linear bijection $T\colon H\to V$ such that
$Tt(X)=s(X)T$ for all $X\in\U$.

\proclaim{Theorem \theoremname{\thmfindimreprU}} For each
dimension $N+1$, $N\in\Zp$, there are four inequivalent irreducible
representations. Explicitly, there exists a basis
$\{ e_0,\ldots,e_N\}$ of $\C^{N+1}$ such that they are given by
$t(A)\, e_k= \la q^{N/2-k}\, e_k$, $t(C)\, e_k = e_{k+1}$, and
$$
t(B)\, e_k=
{{q^{N+1}\la^2(1-q^{-2k})+q^{-1-N}\la^{-2}(1-q^{2k})}\over{ 
(q-q^{-1})^2}}\, e_{k-1}
$$
for $\la^4=1$ with the convention $e_{-1}=0=e_{N+1}$.
\endproclaim

\demo{Proof} Let $t$ be an irreducible representation of $\U$ in 
$H$, and let $\la$ be an eigenvalue of the operator $t(A)$ for the
eigenvector $v$. Then $t(B)^k v$ is an eigenvector for
$t(A)$ for the eigenvalue $\la q^k$. Since $H$ is finite dimensional
$t(B)^kv=0$ for $k$ large enough, so we may assume that
$t(B)v=0$. Define $e_0=v$ and $e_k=t(C)^ke_0$ for $k\in\Zp$.
Let $N$ be the smallest integer such that $e_{N+1}=0$ and
$e_N\not=0$, then $\text{span}(e_0,\ldots,e_N)$ is a non-zero
invariant subspace and hence equals $H$. Then
$t(C)e_k=e_{k+1}$ with the convention $e_{N+1}=0$ and
$t(A)e_k=\la q^{N/2-k}e_k$ after rescaling $\la$.
Since $t(B)e_k$ is an eigenvector
of $t(A)$ for the eigenvalue $\la q^{N/2+1-k}$ we have
$t(B)e_k=\mu_k e_{k-1}$ for some constant $\mu_k$.
The commutation relation for $B$ and $C$ imply
$$
\mu_{k+1}-\mu_{k} = {{\la^2q^{N-2k}-\la^{-2}  q^{2k-N}}
\over{q-q^{-1}}}\ \Longrightarrow \mu_k=
{{q^{N+1}\la^2(1-q^{-2k})+q^{-1-N}\la^{-2}(1-q^{2k})}
\over{ (q-q^{-1})^2}}.
$$
This calculation takes into account the initial condition $\mu_0=0$.
Since we also have the end condition
$\mu_N=(\la^{-2}q^N-\la^2q^{-N})/(q-q^{-1})$ these two expressions
have to be equal, and this leads to the condition $\la^4=1$.

The representations obtained in this way give all finite dimensional
representations of $\U$. For different dimensions they are mutually
inequivalent, and for the same dimension we see that they are 
mutually inequivalent since the spectrum of the operator 
corresponding to $A$ is different.
\qed\enddemo

In the sequel we only use the representations with
the spectrum of $A$ contained in $q^{\hZp}$, or $\la=1$.

%%%%%%%%%%%%%%%%%%%%%%%%%%%%%%%%%%%%%%%%%%%%%%%%%%%%%%%%%%%%%%%%%%%%
%%N E W   S U B S E C T I O N%%%%%%%%%%%%%%%%%%%%%%%%%%%%%%%%%%%%%%%
%%%%%%%%%%%%%%%%%%%%%%%%%%%%%%%%%%%%%%%%%%%%%%%%%%%%%%%%%%%%%%%%%%%%
\subhead\newsubsection 
$\ast$-structures on $\U$
\endsubhead
We have seen that $\U$ is a Hopf algebra, and we now consider
how the Hopf algebra structure depends on $q$ and whether $\U$
can be made into a Hopf $\ast$-algebra. For this purpose we
need some special elements of $\U$. We call $0\not=X\in\U$ 
{\sl group
like} if $\De(X)=X\otimes X$. Note that this implies $\ep(X)=1$.
Let $X$ be group like, then we say that $Y\in\U$ is {\sl twisted
primitive} with respect to $X$ if 
$\De(Y)=X\otimes Y + Y\otimes S(X)$.
Note that this implies $\ep(Y)=0$ and that the twisted primitive
elements form a linear subspace of $\U$.
In case we consider the Hopf
algebra $U({\frak g})$, cf. Example \thmref{\exaHopfalgebratwo},
the unit element $1$ is the only group like
element and the elements of ${\frak g}\subset U({\frak g})$
are the only twisted primitive elements.

\proclaim{Proposition \theoremname{\propclassgroupprimelts}}
In the Hopf algebra $\U$ we have

\noindent
{\rm (i)} $A^m$, $m\in\Z$, are the group like elements;

\noindent
{\rm (ii)} $\C (B,C,A-D)$ is the space of twisted
primitive elements with respect to $A$;

\noindent
{\rm (iii)} $\C(A^m-D^m)$ is the space of twisted
primitive elements with respect to $A^m$ for $m\not= 1$.
\endproclaim

For the proof of Proposition \thmref{\propclassgroupprimelts}
we need the following lemma describing a $q$-analogue
of Newton's binomial formula. For $a\in\C$ we define
the $q$-shifted factorial $(a;q)_n=\prod_{i=0}^{n-1}(1-aq^i)$
with the empty product equal to $1$. The $q$-binomial
coefficient is defined by
$$
\left[ {n\atop k}\right]_q = {{(q;q)_n}\over{(q;q)_k(q;q)_{n-k}}}
= {{(q^n;q^{-1})_k}\over{(q;q)_k}} = \left[ {n\atop{n-k}}\right]_q.
\tag\eqname{\vgldefqbinomialcoef}
$$

\proclaim{Lemma \theoremname{\lemqbinomiallemma}} Let $x$, $y$
be elements of an associative algebra satisfying $xy=qyx$, then
$$
(x+y)^n = \sum_{k=0}^n \left[ {n\atop k}\right]_{q^{-1}}x^ky^{n-k}
= \sum_{k=0}^n \left[ {n\atop k}\right]_q y^kx^{n-k}.
$$
\endproclaim

\demo{Proof} Use the recurrence
$$
\left[ {{n+1}\atop k}\right]_q = q^k\left[ {n\atop k}\right]_q +
\left[ {n\atop{k-1}}\right]_q
\tag\eqname{\vglrecurqbinomialcoef}
$$
and complete induction with respect to $n$.
\qed\enddemo

\demo{Proof of Proposition \thmref{\propclassgroupprimelts}} Use
Lemma \thmref{\lemqbinomiallemma} to see that
$$
\De(C^m)=(A\otimes C+C\otimes D)^m =
\sum_{i=0}^m \left[ {m\atop i}\right]_{q^{-2}} C^iA^{m-i}
\otimes D^iC^{m-i}
$$
and an expression for $\De(B^k)$ can be obtained in this way too. 
So using Lemma \thmref{\lemPBWforU} we consider general
$X=\sum c_{klm}D^lC^mB^k$ and rewriting each factor in the tensor 
product into the basis of Lemma \thmref{\lemPBWforU} we obtain
$$
\multline
\De(X) = \sum_{klm} \sum_{j=0}^k \sum_{i=0}^m c_{klm}
\left[{k\atop j}\right]_{q^2}\left[{m\atop i}\right]_{q^{-2}}
q^{(i-j)(k-j+m-i)}\\
\times D^{l-m+i-k+j}C^iB^j\otimes D^{l+i+j}C^{m-i}B^{k-j}.
\endmultline
\tag\eqname{\vglDeltaongeneralX}
$$

To prove (i) we compare \thetag{\vglDeltaongeneralX} with
$$
X\otimes X = \sum_{uvw} \sum_{abc} c_{uvw}c_{abc} D^vC^wB^u
\otimes D^bC^cB^a.
\tag\eqname{\vglXotimesX}
$$
So $w+c=m$, $a+u=k$, $b-v=m+k$. Suppose $m>0$, then
\thetag{\vglXotimesX} has a non-zero term of the form $D^vC^mB^u
\otimes D^bC^mB^a$, which cannot occur in 
\thetag{\vglDeltaongeneralX}. Hence, $m=0$ and similarly $k=0$. 
So $X=\sum_l c_lD^l$ and
$$
\sum_{l,k} c_lc_k D^l\otimes D^k = \sum_l c_l D^l\otimes D^l
$$
implying $c_l=0$ except for one $l\in\Z$ and the non-zero
$c_l$ has to satisfy $c_l^2=c_l$ or $c_l=1$.

In the proofs of (ii) and (iii) we have to compare
\thetag{\vglDeltaongeneralX} with the appropriate
expressions. These proofs are left as exercises.
\qed\enddemo

Proposition \thmref{\propclassgroupprimelts}(ii) shows that only
for the group like element $A$ the corresponding twisted
primitive elements leads to a three-dimensional space. So from
now on we consider twisted primitive elements only with respect
to $A$ and we consider these elements as the proper analogues of the
three-dimensional Lie algebra ${\frak{sl}}(2,\C)$.

\proclaim{Theorem \theoremname{\thmisomUasHopfalg}}
$U_q({\frak{sl}}(2,\C))\cong U_p({\frak{sl}}(2,\C))$ as Hopf 
algebras if and only if $p=q$ or $p=q^{-1}$.
\endproclaim

\demo{Proof} Let
$\phi\colon U_q({\frak{sl}}(2,\C))\to U_p({\frak{sl}}(2,\C))$ be the
Hopf algebra isomorphism, then $\phi$ maps group like elements onto
group like elements, so $\phi(D^l)=D^{\ta(l)}$. Since $\phi$
is an algebra homomorphism, $\ta$ is an automorphism of $\Z$.
Thus, $\phi(D)=D$ or $\phi(D)=D^{-1}=A$. In the last situation
Proposition \thmref{\propclassgroupprimelts}(ii), (iii) shows that
$\phi$ maps the three-dimensional space of twisted primitive 
elements with respect to $A$ to the one-dimensional space of 
twisted primitive elements with respect to $D$, contradicting 
$\phi$ being an isomorphism.

So $\phi(A)=A$, $\phi(D)=D$ and $\phi$ maps the twisted primitive
elements onto the twisted primitive elements (both with respect 
to $A$). Now the map $X\mapsto AXD$ maps the space of twisted 
primitive elements into itself and has eigenvalues $1$, $q$ and 
$q^{-1}$ for respectively the eigenvectors $A-D$, $B$ and $C$. 
Hence, $\{ 1,q,q^{-1}\}=\{ 1,p,p^{-1}\}$ and thus $p=q$ or 
$p=q^{-1}$.

{}From \thetag{\vglcommrelUq} we see that in case $p=q$ we can 
define $\phi(B)=\la B$, $\phi(C)=\la^{-1}C$ for some non-zero 
$\la\in\C$ and that in case $p=q^{-1}$ we can define
$\phi(B)=\la C$, $\phi(C)= \la^{-1}B$ for some non-zero $\la\in\C$.
It is straightforward to check that $\phi$ defined in this way
gives a Hopf algebra isomorphism.
\qed\enddemo

So we can assume without loss of generality that $|q|\leq 1$.

We say that two $\ast$-structures on a Hopf algebra are
equivalent if there exists a Hopf $\ast$-algebra isomorphism
of the Hopf algebra onto itself intertwining the two
$\ast$-structures. Otherwise, they are inequivalent 
$\ast$-structures.

\proclaim{Theorem \theoremname{\thmstarstructureonU}} The
list of mutually inequivalent $\ast$-structures on $\U$ is

\noindent
{\rm (i)} $|q|=1$; $A^\ast=A$, $B^\ast=-B$, $C^\ast=-C$ and
corresponding real form $U_q({\frak{sl}}(2,\R))$,

\noindent
{\rm (ii)} $-1<q<1$, $q\not=0$; $A^\ast=A$, $B^\ast=C$, 
$C^\ast=B$ and
corresponding real form $U_q({\frak{su}}(2))$,

\noindent
{\rm (iii)} $-1<q<1$, $q\not=0$; $A^\ast=A$, $B^\ast=-C$, 
$C^\ast=-B$ and
corresponding real form $U_q({\frak{su}}(1,1))$.
\endproclaim

\demo{Proof} Suppose $\U$ is a Hopf $\ast$-algebra, then we 
can think
of $\ast$ as an antilinear Hopf algebra isomorphism of $\U$ onto
$\bigl( U_{\bar q}({\frak{sl}}(2,\C))\bigr)_{\text{opp}} \cong
 U_{\bar q^{-1}}({\frak{sl}}(2,\C)$ as Hopf algebras. The 
Hopf algebra isomorphism is just interchanging $B$ and $C$. 
So we must have $q=\bar q$ or $q=\bar q^{-1}$, so $q\in\R$ or 
$|q|=1$. Using the proof of Theorem \thmref{\thmisomUasHopfalg} 
we see that in the
first case we have to have $B^\ast=\la C$, $C^\ast=\la^{-1} B$
for some non-zero $\la\in\R$ and in the second case
we have to have $B^\ast=\la B$, $C^\ast=\la^{-1} C$ for $|\la|=1$.
The condition on the $\la$ follows from the $\ast$-operator being
an involution.
Now use Theorem \thmref{\thmisomUasHopfalg} again to pick out the
inequivalent $\ast$-structures.
\qed\enddemo

\demo{Remark \theoremname{\remthmstarstructureonU}} The names 
for these real forms are motivated by the fact that in case 
$q\uparrow 1$ the $-1$-eigenspace of the corresponding 
$\ast$-operator in the Lie algebra,
cf. Example \thmref{\exaHopfalgebrafour}, are ${\frak{sl}}(2,\R)$,
${\frak{su}}(2)$ and ${\frak{su}}(1,1)$. Note that
${\frak{su}}(1,1)=C\,{\frak{sl}}(2,\R)\, C^{-1}$ with
$C=\pmatrix 1&-i\\1&i\endpmatrix$, so
${\frak{sl}}(2,\R)$ and ${\frak{su}}(1,1)$
are conjugate and hence have the
same representation theory. Theorem \thmref{\thmstarstructureonU}
shows that this is no longer true in the quantum case.
The $\ast$-operator from Theorem \thmref{\thmstarstructureonU}(ii)
and corresponding Hopf $\ast$-algebra is sometimes called the
compact real form, since $SU(2)$ is a compact group.
\enddemo

{}From now on, when considering the real forms $\Usu$
or $U_q({\frak{su}}(1,1))$ we always take $0<q<1$. This can be done
without much loss of generality, since it turns out that the
special functions associated to these Hopf $\ast$-algebras 
essentially depend on $q^2$.

It is straightforward to see that for $U_q({\frak{su}}(2))$
the finite dimensional representations of Theorem 
\thmref{\thmfindimreprU}
are unitarisable for $\la =\pm 1$. We only consider
$\la=1$ and we redefine these representations in the following 
theorem.

\proclaim{Theorem \theoremname{\thmunitaryrtepsofUsu}}
For each spin $l\in\hZp$ there
exists a unique $(2l+1)$-dimensional $\ast$-re\-pre\-sen\-tation
of $\Usu$ such that the spectrum of $A$ is contained in $q^{\hZp}$.
Equip $\C^{2l+1}$ with orthonormal basis $\{ e^l_n\}$, 
$n=-l,-l+1,\ldots,l$ and denote the representation by $t^l$. The
action of the generators is given by
$$
\aligned
&t^l(A)\, e^l_n=q^{-n}e^l_n,\qquad t^l(D)\, e^l_n=q^n\, e^l_n, \\
&t^l(B)\, e^l_n= 
{{\sqrt{(q^{-l+n-1}-q^{l-n+1})(q^{-l-n}-q^{l+n})}}\over
{q^{-1}-q}} \,  e^l_{n-1} \\
&t^l(C)\, e^l_n= 
{{\sqrt{(q^{-l+n}-q^{l-n})(q^{-l-n-1}-q^{l+n+1})}}\over
{q^{-1}-q}}\,  e^l_{n+1},
\endaligned
\tag\eqname{\vgldefreprU}
$$
where $e^l_{l+1}=0=e^l_{-l-1}$.
\endproclaim

%%%%%%%%%%%%%%%%%%%%%%%%%%%%%%%%%%%%%%%%%%%%%%%%%%%%%%%%%%%%%%%%%%%%
%%N E W   S U B S E C T I O N%%%%%%%%%%%%%%%%%%%%%%%%%%%%%%%%%%%%%%%
%%%%%%%%%%%%%%%%%%%%%%%%%%%%%%%%%%%%%%%%%%%%%%%%%%%%%%%%%%%%%%%%%%%%
\subhead\newsubsection 
Dual Hopf algebra
\endsubhead
Let us now consider the fundamental two-dimensional representation
of $\U$ for spin $l=1/2$. The matrix elements of $t^{1/2}$ give four
linear functionals on $\U$;
$$
\cases t^{1/2}(X)\, e_{-1/2}^{1/2} = \al(X)\, e^{1/2}_{-1/2}
+ \be(X)\ e^{1/2}_{1/2} \\
t^{1/2}(X)\, e_{1/2}^{1/2} = \ga(X)\, e^{1/2}_{-1/2}
+ \de(X)\ e^{1/2}_{1/2}
\endcases
\Longleftrightarrow
\quad t^{1/2}(X) = \pmatrix \al(X) & \be(X)\\ \ga(X) &\de(X) 
\endpmatrix.
$$
It follows from \thetag{\vgldefreprU} that on the basis of Lemma
\thmref{\lemPBWforU} the linear functionals are given by
$$
\gathered
\al(D^lC^mB^n)=\de_{n0}\de_{m0}q^{-l/2}, \qquad
\be(D^lC^mB^n)=\de_{n1}\de_{m0}q^{-l/2}, \\
\ga(D^lC^mB^n)=\de_{n0}\de_{m1}q^{l/2}, \qquad
\de(D^lC^mB^n)=(\de_{n0}\de_{m0}+\de_{n1}\de_{m1})q^{l/2}.
\endgathered
\tag\eqname{\vglabgdonbasis}
$$

\proclaim{Theorem \theoremname{\thmdescripAasdualHAofU}} Let
$\A$ be the complex unital associative subalgebra of
the linear dual of $\U$ generated by $\al$, $\be$, $\ga$, $\de$.
Then the following relations hold;
$$
\gathered
\al\be =q\be\al ,\quad \al\ga = q\ga\al ,\quad \be\de = q\de\be ,
\quad \ga\de = q\de\ga ,\\
\be\ga =\ga\be ,\quad \al\de -q\be\ga = \de\al - q^{-1}\be\ga =1 .
\endgathered
\tag\eqname{\vglcommrelAq}
$$
Then $\A$ is Hopf algebra.
The comultiplication $\De$, the counit $\ep$ and the antipode $S$
given on the generators by
$$
\gather
\eqalign{
&\De(\al)=\al\otimes\al + \be\otimes\ga ,\quad
\De(\be )=\al\otimes\be + \be\otimes\de ,\cr
&\De(\ga )=\ga\otimes\al + \de\otimes\ga ,\quad
\De(\de )=\ga\otimes\be + \de\otimes\de , \cr}
\tag\eqname{\vgldefDeltaonAq}\\
\ep\pmatrix \al &\be\\ \ga&\de\endpmatrix = \pmatrix
1&0\\0&1\endpmatrix,\qquad
S\pmatrix \al &\be\\ \ga&\de\endpmatrix = \pmatrix \de&-q^{-1}\be \\
-q\ga&\al\endpmatrix.
\tag\eqname{\vglSenepsilonopAq}
\endgather
$$
$\A$ has a linear basis formed by the matrix elements
$t^l_{m,n}\colon X\mapsto \langle t^l(X)e^l_n,e^l_m\rangle$.
\endproclaim

For the proof we need a lemma, which is of interest on its own.
First we discuss the notion of
a tensor product representation of $\U$. If $t$ and $s$ are
representations of $\U$ in $V$ and $W$, then we define the
tensor product representation $t\otimes s$ of $\U$ in $V\otimes W$
by using the comultiplication $\De$, cf. \thetag{\vglSweedlernot},
$$
\bigl(t\otimes s\bigr)(X) v\otimes w =
\sum_{(X)} t(X_{(1)})\, v\otimes s(X_{(2)})\, w, \qquad X\in\U.
$$
Then we have the following Clebsch-Gordan decomposition.

\proclaim{Lemma \theoremname{\lemmaCGCforUsu}}
$t^{l_1}\otimes t^{l_2} \cong \bigoplus_{l=|l_1-l_2|}^{l_1+l_2} t^l$
\endproclaim

\demo{Proof} Although Lemma \thmref{\lemmaCGCforUsu} is not 
concerned with $\ast$-structures we use the $\ast$-structure from 
$\Usu$. Then $t^{l_1}\otimes t^{l_2}$ is a finite-dimensional 
unitary representation of $\Usu$ and hence completely reducible. 
Since the spectrum of $\bigl(t^{l_1}\otimes t^{l_2}\bigr)(A)$ is 
contained in $q^{\hZp}$ we have a decomposition of the form 
$t^{l_1}\otimes t^{l_2}=\bigoplus_l
m_l\ t^l$ for certain multiplicities $m_l$. Since
$e_n^{l_1}\otimes e_m^{l_2}$ is an eigenvector of the action for $A$
for the eigenvalue $q^{-n-m}$ we can read off the multiplicities
from Theorem \thmref{\thmunitaryrtepsofUsu} by counting eigenvalues
for $A$.
\qed\enddemo

\demo{Proof of Theorem \thmref{\thmdescripAasdualHAofU}}
$\A$ is automatically a Hopf algebra by
Definition \thmref{\defdualHopfalgebras}
if we can show that comultiplication maps into the (algebraic)
tensor product $\A\otimes\A$. The definition of
comultiplication, counit and antipode follow from
Definition \thmref{\defdualHopfalgebras}. Using
$t^{1/2}(XY)=t^{1/2}(X)t^{1/2}(Y)$ we find
that $\langle \De(\al),X\otimes Y\rangle =\al(XY) = \al(X)\al(Y) +
\be(X)\ga(Y)=\langle\al\otimes\al+\be\otimes\ga,X\otimes Y\rangle$
by inspecting the upper left entry. Inspection of the other entries
leads to the action of the comultiplication on the other generators.
In particular, $\A$ is a Hopf algebra.

The action of the counit follows by taking $l=m=n=0$ in
\thetag{\vglabgdonbasis}.
The action of $S$ can be calculated by $S(D^lC^mB^n)=
(-q)^{m-n}B^nC^mA^l$, so that $\langle S(\al),D^lC^mB^n\rangle
=(-q)^{m-n}\al(B^nC^mA^l)=(\de_{n0}\de_{m0}+\de_{n1}\de_{m1})q^{l/2}
=\de(D^lC^mB^n)$ by \thetag{\vglabgdonbasis}, and similarly for the
other generators. From \thetag{\vgldiagramS} we see that
$$
\pmatrix \al &\be\\ \ga&\de\endpmatrix  \pmatrix \de&-q^{-1}\be \\
-q\ga&\al\endpmatrix =\pmatrix 1&0\\0&1\endpmatrix =
\pmatrix \de&-q^{-1}\be \\ -q\ga&\al\endpmatrix
\pmatrix \al &\be\\ \ga&\de\endpmatrix
$$
which implies the relations \thetag{\vglcommrelAq}.

Consider the matrix element $(t^{l_1}\otimes t^{l_2})_{ij;lm}=
t^{l_1}_{il}t^{l_2}_{jm}$. Take $l_1=1/2$ and use induction with
respect to $l_2$ and Lemma \thmref{\lemmaCGCforUsu} to see that
the matrix elements of all spin representations are contained in
$\A$. Iterating Lemma \thmref{\lemmaCGCforUsu} shows that each
matrix element $t^l_{mn}$ can be written in terms of the generators.
In order to prove the linear independence we need the following
lemma.

\proclaim{Lemma \theoremname{\lemHaarfunconAone}} Define 
$h\colon\A\to\C$
by $h(1)=1$, $h(t^l_{mn})=0$ for $l>0$, then $h$ is an invariant
functional on $\A$ and the Schur orthogonality relations hold;
$$
h\bigl( S(t^l_{mn}) t^k_{ij}\bigr) = \de_{lk}\de_{mj}\de_{in}
q^{2(l-n)} {{1-q^2}\over{1-q^{4l+2}}}.
$$
\endproclaim

\demo{Proof of Lemma \thmref{\lemHaarfunconAone}} First observe that
$1=t^0_{00}$ and that $1$ cannot be written as a linear combination
of matrix elements $t^l_{mn}$ for $l>0$, since otherwise the
trivial representation would occur as a subrepresentation of $t^l$
contradicting its irreducibility. So $h$ is well-defined on the 
whole of $\A$. Since 
$\De(t^l_{mn})=\sum_{k=-l}^l t^l_{mk}\otimes t^l_{kn}$
by $t^l(XY)=t^l(X)t^l(Y)$, we immediately see that $h$ is an
invariant functional.

In order to prove the orthogonality relations we first observe
that
$$
\sum_{j=-k}^k h\bigl( S(t^l_{mn}) t^k_{ij}\bigr) \, t^k_{jp}
= \sum_{j=-l}^l t^l_{mj}\, h\bigl( S(t^l_{jn}) t^k_{ip}\bigr).
\tag\eqname{\vglhasintertwiner}
$$
This follows from the right invariance of $h$;
$$
h\bigl( S(t^l_{mn}) t^k_{ij}\bigr)\, 1 =
(h\otimes id)\bigl(\De(S(t^l_{mn}) t^k_{ij})\bigr) =
\sum_{r=-l}^l\sum_{s=-k}^k h\bigl( S(t^l_{rn}) t^k_{is}) \,
S(t^l_{mr}) t^k_{sj}
$$
by Proposition \thmref{\proppropsofHopfalgebra}(ii).
{}From \thetag{\vgldiagramS} we obtain $\sum_{m=-l}^l
t^l_{km}\, S(t^l_{mr}) = \ep(t^l_{kr})=\de_{kr}$ and this leads
to \thetag{\vglhasintertwiner}.
So the matrix $T^{(n,i)}_{mj} = h\bigl( S(t^l_{mn}) t^k_{ij}\bigr)$
intertwines $t^l$ and $t^k$. Hence, it is zero for $k\not=l$ and
if $k=l$ we have $T^{(n,i)}_{mj}=c^{(n,i)}\de_{mj}$ for some
constant $c^{(n,i)}\in\C$.

For a representation $t$ of $\U$ in $V$ the contragredient
representation $t^c$ of $\U$ in $V$ is defined by
$\langle t^c(X)e_n,e_m\rangle = 
\langle t\bigl(S(X)\bigl)e_m,e_n\rangle$,
or in terms of the antipode $S$ of $\A$,
$t^c_{mn}=S(t_{nm})$ for a representation $t$ of $\U$ such that 
the spectrum of $A$ is contained in $q^{\hZp}$.
Using Theorem \thmref{\thmunitaryrtepsofUsu}
we see that the contragredient of $t^l$ is equivalent to $t^l$, and
consequently $h(S(\xi))=h(\xi)$ for all $\xi\in\A$. So we get
$T^{(n,i)}_{mj}= h\bigl(  S(t^k_{ij}) S^2(t^l_{mn})\bigr)$.
Now $S^2(t^l_{mn})$ are the matrix coefficients of the double
contragredient representation of $t^l$ and hence there exists
an invertible intertwining operator $F\in End(\C^{2l+1})$ such
that $(t^l)^{cc}(X)F=Ft^l(X)$ for all $X\in\U$ and hence
$S^2(t^l_{mn})=\sum_{p,r=-l}^l F_{mp}t^l_{pr} (F^{-1})_{rn}$.
Thus,
$$
\de_{kl}\de_{mj} c^{(n,i)} = T^{(n,i)}_{mj}= \sum_{p,r=-l}^l F_{mp}
h\bigl(  S(t^k_{ij})t^l_{pr}\bigr) (F^{-1})_{rn} =
\de_{kl} (F^{-1})_{in} \sum_{p=-l}^l F_{mp} c^{(j,p)},
$$
or $c^{(n,i)}=c(F^{-1})_{in}$ for some constant $c\in\C$.
To determine the constant we observe that
$\sum_{m=-l}^l S(t^l_{km}) t^l_{mr}=\de_{kr}$ implies
$1 = \sum_{m=-l}^l c^{(m,m)}$, so $c^{-1} =\text{tr}(F^{-1})$.
Summarising,
$$
h\bigl( S(t^l_{mn})\, t^k_{ij}\bigr) =\de_{kl}\de_{mj}
{{ (F^{-1})_{in}}\over{\text{tr}(F^{-1})}}.
\tag\eqname{\vglgeneralSchurorth}
$$
In this case we easily show that $Fe^l_k=q^{2k}e^l_k$ is such
an intertwiner. Then \thetag{\vglgeneralSchurorth}
and a calculation finish the proof.
\qed\enddemo

Finally, assume that $\sum_{lmn} c_{lmn}t^l_{mn}=0$, then multiply
this element from the left
by $S(t^l_{nm})$
and apply the invariant functional $h$.
We get $c_{lmn} h\bigl( S(t^l_{mn}) t^l_{mn}\bigr)=0$
from Lemma \thmref{\lemHaarfunconAone} and this shows $c_{lmn}=0$.
\qed\enddemo

\demo{Remark \theoremname{\remlimtqtooneofA}} If we let 
$q\uparrow 1$
in the definition of $\A$ we obtain a commutative Hopf algebra
which is nothing but $\text{Pol}(SL(2,\C))$ by identifying
$\al$, $\be$, $\ga$ and $\de$ with the coordinate functions.
So for $g\in SL(2,\C)$ we let $\al(g)=g_{11}$, $\be(g)=g_{12}$,
$\ga(g)=g_{21}$ and $\de(g)=g_{22}$.
\enddemo

%%%%%%%%%%%%%%%%%%%%%%%%%%%%%%%%%%%%%%%%%%%%%%%%%%%%%%%%%%%%%%%%%%%%
%%N E W   S U B S E C T I O N%%%%%%%%%%%%%%%%%%%%%%%%%%%%%%%%%%%%%%%
%%%%%%%%%%%%%%%%%%%%%%%%%%%%%%%%%%%%%%%%%%%%%%%%%%%%%%%%%%%%%%%%%%%%
\subhead\newsubsection 
More on the dual Hopf algebra
\endsubhead
Do the relations in \thetag{\vglcommrelAq} describe all the
relations between the generators of $\A$? We show that the
answer is yes.

\proclaim{Lemma \theoremname{\lemabstractalgsrtuc}} Let $B$
be the algebra generated by $\al$, $\be$, $\ga$ and $\de$
subject to the relations of \thetag{\vglcommrelAq}, then
a linear basis for $B$ is given by
$\de^l\ga^m\be^n$, $l,m,n\in\Zp$, and
$\al^l\ga^m\be^n$, $l\in\N$, $m,n\in\Zp$.
\endproclaim

\demo{Proof} From
\thetag{\vglcommrelAq} it follows that these elements span $B$.
We have to show that they are linearly independent.

Consider the infinite dimensional representation
of $B$ in $\ell^2(\Zp\times \Z)$ with standard orthonormal
basis $e_{n,k}$, $n\in\Zp$, $k\in\Z$. The action of the
generators is given by
$$
\gather
\al\, e_{n,k} = \sqrt{1-q^{2n}}\, e_{n-1,k}, \quad
\be\, e_{n,k} = -q^{n+1}\, e_{n,k-1}, \\
\ga\, e_{n,k} = q^n\, e_{n,k+1}, \quad
\de\, e_{n,k} = \sqrt{1-q^{2n+2}}\, e_{n+1,k}.
\endgather
$$
Using this representation we show that a non-trivial linear
combination of such elements cannot give the zero operator.
\qed\enddemo

To show that the same elements also give a basis of $\A$
we proceed by explicitly calculating the duality for
such elements. The exercises contain some hints on how
to prove the following theorem.

\proclaim{Theorem \theoremname{\thmexplicitdualityonbasis}}
Define
$$
C^{L,M,N}_{l,m,n} = q^{l(L+M-N)/2} q^{-L(m+n)/2}
q^{-m(m-1)/2}q^{-n(n-1)/2}
{{(q^2;q^2)_n(q^2;q^2)_m}\over{(1-q^2)^{m+n}}},
$$
then
$$
\langle \al^L\ga^M\be^N,D^lC^mB^n\rangle = \de_{Mm}\de_{Nn}
C^{-L,M,N}_{l,m,n}
$$
and
$$
\langle \de^L\ga^M\be^N,D^lC^mB^n\rangle =
\cases {\dsize{ q^{(m-M)^2} \left[ {L\atop{m-M}}\right]_{q^2}
C^{L,M,N}_{l,m,n} }}, &\text{if $0\leq m-M=n-N\leq L$,}\\
0&\text{otherwise}.
\endcases
$$
\endproclaim

\proclaim{Corollary \theoremname{\corthmexplicitdualityonbasisone}}
A linear basis for $\A$ is given by
$\de^l\ga^m\be^n$, $l,m,n\in\Zp$, and
$\al^l\ga^m\be^n$, $l\in\N$, $m,n\in\Zp$.
In particular, $\A$ is isomorphic to $B$ defined in Lemma
\thmref{\lemabstractalgsrtuc}, and \thetag{\vglcommrelAq}
are the only relations in $\A$.
\endproclaim

\demo{Proof} Suppose that in $\A$ we
have
$$
0= \sum c_{LMN} \de^L\ga^M\be^N + \sum c_{-LMN} \al^L\ga^M\be^N.
$$
Let $m$ be the minimal $M$ such that $c_{LMN}\not=0$, and $n$
the minimal $N$ such that $c_{LmN}\not=0$.
Testing against $D^lC^mB^n$ shows
$$
0=\sum_{L\in\Z} c_{Lmn} C(m,n) q^{-L(m+n)/2} q^{lL/2}, 
\qquad \forall
\, l\in\Z,
$$
for some non-zero constant $C(m,n)$ by Theorem
\thmref{\thmexplicitdualityonbasis}. So $c_{Lmn}=0$ for all $L$.
\qed\enddemo

\proclaim{Corollary \theoremname{\corthmexplicitdualityonbasistwo}}
The duality between $\U$ and $\A$, considered as an
abstract algebra generated by $\al$, $\be$, $\ga$ and $\de$,
subject to the relations \thetag{\vglcommrelAq},
defined on the generators by
$$
\aligned
\Bigl\langle A, \pmatrix \al&\be\\ \ga&\de\endpmatrix \Bigr\rangle =
\pmatrix q^{1/2}&0 \\ 0&q^{-1/2}\endpmatrix , \quad
\Bigl\langle B, \pmatrix \al&\be\\ \ga&\de\endpmatrix \Bigr\rangle =
\pmatrix 0&1 \\ 0&0 \endpmatrix , \\
\Bigl\langle C, \pmatrix \al&\be\\ \ga&\de\endpmatrix \Bigr\rangle =
\pmatrix 0&0 \\ 1&0 \endpmatrix , \quad
\Bigl\langle D, \pmatrix \al&\be\\ \ga&\de\endpmatrix \Bigr\rangle =
\pmatrix q^{-1/2}&0 \\ 0&q^{1/2}\endpmatrix .
\endaligned
\tag\eqname{\vglexpldualityongenAU}
$$
and extended as Hopf algebra duality, is doubly non-degenerate.
\endproclaim

%%%%%%%%%%%%%%%%%%%%%%%%%%%%%%%%%%%%%%%%%%%%%%%%%%%%%%%%%%%%%%%%%%%%
%%N E W   S U B S E C T I O N%%%%%%%%%%%%%%%%%%%%%%%%%%%%%%%%%%%%%%%
%%%%%%%%%%%%%%%%%%%%%%%%%%%%%%%%%%%%%%%%%%%%%%%%%%%%%%%%%%%%%%%%%%%%
\subhead\newsubsection 
Action of $\U$ on $\A$ and $\ast$-structures on $\A$
\endsubhead
Now that we have established $\U$ and $\A$ in duality as Hopf
algebras, we may consider the action of $\U$ on $\A$ as defined
in Proposition \thmref{\propactionUonA}. It is a simple calculation
to give the action on the level of generators using
\thetag{\vgldefDeltaonAq} and \thetag{\vglexpldualityongenAU}.
We get in an obvious notation
$$
\gathered
A. \pmatrix \al&\be\\ \ga&\de\endpmatrix =
\pmatrix q^{1/2}\al&q^{-1/2}\be\\ q^{1/2}\ga&q^{-1/2}\de\endpmatrix,
\quad
\pmatrix \al&\be\\ \ga&\de\endpmatrix . A =
\pmatrix q^{1/2}\al&q^{1/2}\be\\ q^{-1/2}\ga&q^{-1/2}
\de\endpmatrix, \\
B. \pmatrix \al&\be\\ \ga&\de\endpmatrix =
\pmatrix 0 &\al\\ 0& \ga\endpmatrix,
\quad
\pmatrix \al&\be\\ \ga&\de\endpmatrix . B =
\pmatrix \ga &\de\\ 0& 0\endpmatrix, \\
C. \pmatrix \al&\be\\ \ga&\de\endpmatrix =
\pmatrix \be&0\\ \de& 0\endpmatrix,
\quad
\pmatrix \al&\be\\ \ga&\de\endpmatrix . C =
\pmatrix 0&0\\ \al& \be\endpmatrix.
\endgathered
\tag\eqname{\vglexplactioUonAongen}
$$
For $q\uparrow 1$ \thetag{\vglexplactioUonAongen} corresponds
to the action between $U(\frak{sl}(2,\C))$ and
$\text{Pol}(SL(2,\C))$ as described in Example
\thmref{\exampairingGandUg}.

Corresponding to the three inequivalent $\ast$-structures on
$\U$ described in Theorem \thmref{\thmstarstructureonU}
we obtain three $\ast$-structures on $\A$ making it into
a Hopf $\ast$-algebra by transposing the $\ast$-operator
from $\U$ to $\A$ by Definition \thmref{\defdualHopfalgebras}.

\proclaim{Theorem \theoremname{\thmstarstructureonA}} The
list of mutually inequivalent $\ast$-structures on $\A$ is

\noindent
{\rm (i)} $|q|=1$; $\al^\ast=\al$, $\be^\ast=q^{-1}\be$, 
$\ga^\ast=q\ga$,
$\de^\ast=\de$ and
corresponding real form $A_q(SL(2,\R))$,

\noindent
{\rm (ii)} $-1<q<1$, $q\not=0$; $\al^\ast=\de$, $\be^\ast=-q\ga$,
$\ga^\ast=-q^{-1}\be$, $\de^\ast=\al$ and
corresponding real form $\Asu$,

\noindent
{\rm (iii)} $-1<q<1$, $q\not=0$; $\al^\ast=\de$, $\be^\ast=q\ga$,
$\ga^\ast=q^{-1}\be$, $\de^\ast=\al$ and
corresponding real form $A_q(SU(1,1))$.
\endproclaim

\proclaim{Proposition \theoremname{\propSchurorthoforSU}} 
The invariant functional $h$ on $\Asu$ satisfies
$$
h\bigl( (t^l_{nm})^\ast t^k_{ij}\bigr) = \de_{lk}\de_{mj}\de_{in}
q^{2(l-n)} {{1-q^2}\over{1-q^{4l+2}}}.
$$
In particular, $h\colon\Asu\to\C$ is a positive linear functional,
$h(\xi^\ast\xi)>0$ for $0\not=\xi\in\Asu$.
\endproclaim

\demo{Proof} Since $q$ is real,
Lemma \thmref{\lemHaarfunconAone} shows that it suffices
to prove $S(t^l_{ij})=(t^l_{ji})^\ast$. By Theorem
\thmref{\thmunitaryrtepsofUsu} we know that $t^l$ is a
$\ast$-representation of $\Usu$, so
$$
\langle (t^l_{ij})^\ast, X\rangle =
\overline{\langle t^l_{ij}, S(X)^\ast\rangle} =
\langle e^l_i, t^l\bigl( S(X)^\ast\bigr)e^l_j\rangle =
\langle t^l\bigl( S(X)\bigr)e^l_i,e^l_j\rangle =
\langle S(t^l_{ji}), X\rangle
$$
for all $X\in\Usu$.
\qed\enddemo

Proposition \thmref{\propSchurorthoforSU} shows that
the Hopf $\ast$-algebra $\Asu$ has the proper $q$-analogue
of the Schur orthogonality relations. It is also a nice
$\ast$-structure since we can give a complete list of
mutually inequivalent $\ast$-representations of $\Asu$.
It is straightforward to check that the $\ast$-representations
defined in the next theorem are indeed representations.

\proclaim{Theorem \theoremname{\thmirredstarrepsofAsu}}
The following is a complete list of irreducible inequivalent
$\ast$-represen\-ta\-tions
of $\Asu$.

\noindent
{\rm (i)}
The one-dimensional $\ast$-representations $\pi_\th$ defined by
$\pi_\th(\al)=e^{i\theta}$ and $\pi_\th(\be)=0$
for $\th\in[0,2\pi)$.

\noindent
{\rm (ii)}
Infinite dimensional $\ast$-representations
acting in the Hilbert space $\ell^2(\Zp)$. For an orthonormal
basis $\{ e_n\}_{n\in\Zp}$ the action of the generators is given by
$$
\pi^\infty_\th(\al)\, e_n = \sqrt{1-q^{2n}}\, e_{n-1}, \quad
\pi^\infty_\th(\ga)\, e_n = e^{i\theta}q^n\, e_n,
$$
with the convention $e_{-1}=0$.
\endproclaim

\demo{Sketch of Proof} It is straightforward to see that
the representations given in Theorem 
\thmref{\thmirredstarrepsofAsu} are irreducible 
$\ast$-representations of $\Asu$. Conversely, let $\pi$
be an irreducible $\ast$-representation of $\Asu$, then it
follows from the commutation relations 
\thetag{\vglcommrelAq} that the kernel of $\pi(\ga)$ is
an irreducible subspace. Hence, $\pi(\ga)=0$, leading
to the one-dimensional representations, or $\text{ker}(\pi(\ga))$
is trivial. 

In the last case we can use the Spectral Theorem
for the normal operator $\pi(\ga)$ to show that the spectrum
is of the form $\la q^k$, $k\in\Zp$, for some $\la\in\C$. From
the commutation relations \thetag{\vglcommrelAq} it follows
that $|\la|=1$, and that $\pi(\al)$ and $\pi(\de)$ are
acting as shift operators on the eigenvectors. We then find
that $\pi$ is equivalent to an infinite dimensional 
representation as in (ii), since the eigenspaces of $\pi(\ga)$
are one-dimensional.
\qed\enddemo

\demo{Remark \theoremname{\remdefinnonunitaryonedrepsofAsu}}
For $\lambda\not= 0$
we also define one-dimensional representations of $\A$ given by
$$
\ta_\la(\al)=\la,\quad \ta_\la(\be)=0=\ta_\la(\ga),\quad
\ta_\la(\de)=\la^{-1}.
$$
Note that $\ta_\la$ is a $\ast$-representation of $\Asu$ if 
and only if
$\la=e^{i\th}$, or $\ta_\la=\pi_\th$, for some $\th\in [0,2\pi)$.
The counit
$\ep$ coincides with the special case $\ta_1=\pi_0$.
\enddemo

%%%%%%%%%%%%%%%%%%%%%%%%%%%%%%%%%%%%%%%%%%%%%%%%%%%%%%%%%%%%%%%%%%%%
%%N E W   S U B S E C T I O N%%%%%%%%%%%%%%%%%%%%%%%%%%%%%%%%%%%%%%%
%%%%%%%%%%%%%%%%%%%%%%%%%%%%%%%%%%%%%%%%%%%%%%%%%%%%%%%%%%%%%%%%%%%%
%NOTES AND REFERENCES%%%%%%%%%%%%%%%%%%%%%%%%%%%%%%%%%%%%%%%%%%%%%%%
\subhead Notes and references
\endsubhead
The quantised universal enveloping algebra for ${\frak{sl}}(2,\C)$
is the simplest case of a series of quantised universal enveloping 
algebras. In fact, there is a canonical way to associate to any
simple Lie algebra $\frak g$ a quantised universal enveloping
algebra $U_q({\frak g})$. The PBW-basis in Lemma \
thmref{\lemPBWforU} can be proved in various ways, such as by use 
of the representations in Theorem \thmref{\thmfindimreprU} or by 
using the so-called diamond lemma.
As is the case for ${\frak g}={\frak{sl}}(2,\C)$ the representation
theory of simple $\frak g$ and $U_q({\frak g})$ is very similar,
see \cite{\CharP} and references therein.

The approach to $\U$ and its dual algebra presented here
is much inspired by the
paper \cite{\MasuMNNSU} by Masuda et al.
The existence of the invariant functional and
the corresponding
orthogonality relations for the invariant functional have several
proofs, see \cite{\DijkK}, \cite{\KoorTrento},
\cite{\VDaelPAMS}, \cite{\Woro}.

Theorem \thmref{\thmirredstarrepsofAsu} is due to Vaksman and
Soibelman \cite{\VaksS} and can be used to complete $\Asu$ into
a $C^\ast$-algebra, and so making the connection with
Woronowicz's \cite{\Woro} approach to compact quantum groups.
See also \cite{\KoelComp} for the details of the proof of Theorem 
\thmref{\thmirredstarrepsofAsu}.
For general simple compact quantum groups
the $\ast$-representations have been classified by
Soibelman, see references in \cite{\CharP}.

%%%%%%%%%%%%%%%%%%%%%%%%%%%%%%%%%%%%%%%%%%%%%%%%%%%%%%%%%%%%%%%%%%%%
%%N E W   S U B S E C T I O N%%%%%%%%%%%%%%%%%%%%%%%%%%%%%%%%%%%%%%%
%%%%%%%%%%%%%%%%%%%%%%%%%%%%%%%%%%%%%%%%%%%%%%%%%%%%%%%%%%%%%%%%%%%%
%EXERCISES%%%%%%%%%%%%%%%%%%%%%%%%%%%%%%%%%%%%%%%%%%%%%%%%%%%%%%%%%%
\subhead Exercises
\endsubhead

\item{\the\sectionno.1} Finish the proof of Proposition
\thmref{\propUqisHopfalgebra}.

\item{\the\sectionno.2} Verify \thetag{\vglrecurqbinomialcoef}.

\item{\the\sectionno.3} Prove (ii) and (iii) of Proposition
\thmref{\propclassgroupprimelts}.

\item{\the\sectionno.4} Replace in $\Usu$ its generators $A$, 
$B$, $C$
and $D$ by $A$, $\rho^{-1}B$, $\rho^{-1}C$ and $D$. Let 
$\rho\downarrow 0$
and determine the Hopf $\ast$-algebra structure in the limit case.
This Hopf $\ast$-algebra is denoted by $U_q({\frak{m}}(2))$, and is
related to the Lie algebra for the group of orientation and distance
preserving motions of the Euclidean plane.

\item{\the\sectionno.5} The proof of Theorem
\thmref{\thmexplicitdualityonbasis} has to be done in the following
stages, which each can be proved using induction and
Definition \thmref{\defdualHopfalgebras}.
\itemitem{$\bullet$} First show that for $X\in\U$
$$
\align
\langle\de^L\ga^M\be^N,D^lX\rangle &= q^{l(L+M-N)/2}
\langle \de^L\ga^M\be^N,X\rangle, \\
\langle\al^L\ga^M\be^N,D^lX\rangle&= q^{-l(L+M-N)/2}
\langle\al^L\ga^M\be^N,X\rangle.
\endalign
$$
\itemitem{$\bullet$} Show that
$$
\align
\langle\de^L,C^mB^n\rangle&= \cases q^{m(1-L)}
{\dsize{  {{(q^{2L};q^{-2})_m(q^2;q^2)_m}\over{(1-q^2)^{2m}}}   ,}}
&\text{if $m=n\leq L$,}\\ 0,&\text{otherwise,} \endcases \\
\langle\al^L,C^mB^n&\rangle= \de_{m0}\de_{n0}.
\endalign
$$
\itemitem{$\bullet$} Show that
$$
\langle\ga^M\be^N,C^mB^n\rangle=\de_{mM}
\de_{nN}q^{-(n(n-1)/2}q^{-m(m-1)/2}
{{(q^2;q^2)_m(q^2;q^2)_n}\over{(1-q^2)^{m+n}}}
$$
and finish the proof of Theorem
\thmref{\thmexplicitdualityonbasis}.

\item{\the\sectionno.6} Let $k>0$ and let 
$\{e_n\}_{n=0}^\infty$ be the
standard orthonormal basis for $\ell^2(\Zp)$. Prove that there 
exists a unitary representation of $U_q({\frak{su}}(1,1))$ in 
$\ell^2(\Zp)$ such that $A\cdot e_n=q^{k+n}\, e_n$, and 
$C\cdot e_0=0$. Give an explicit expression for the action of 
$B$ and $C$ on a basis vector.
These representations of $U_q({\frak{su}}(1,1))$ are {\sl
positive discrete series} representations. (Be aware that $D$
is represented by an unbounded operator, since $-1<q<1$.)

\item{\the\sectionno.7} Prove that the representation in the proof 
of Lemma \thmref{\lemabstractalgsrtuc} is also a 
$\ast$-representation of $\Asu$. Decompose this 
$\ast$-representation in terms of the
irreducible $\ast$-represen\-ta\-tions of
Theorem \thmref{\thmirredstarrepsofAsu}.

%%%%%%%%%%%%%%%%%%%%%%%%%%%%%%%%%%%%%%%%%%%%%%%%%%%%%%%%%%%%%%%%%%%%
%%N E W   S E C T I O N%%%%%%%%%%%%%%%%%%%%%%%%%%%%%%%%%%%%%%%%%%%%%
%%%%%%%%%%%%%%%%%%%%%%%%%%%%%%%%%%%%%%%%%%%%%%%%%%%%%%%%%%%%%%%%%%%%
\newpage

\head\newsection Orthogonal polynomials and basic
hypergeometric series\endhead

In this lecture we consider first some general theory of
orthogonal polynomials and next some explicit examples
which can be written in terms of basic hypergeometric series,
in particular the famous Askey-Wilson polynomials.

%%%%%%%%%%%%%%%%%%%%%%%%%%%%%%%%%%%%%%%%%%%%%%%%%%%%%%%%%%%%%%%%%%%%
%%N E W   S U B S E C T I O N%%%%%%%%%%%%%%%%%%%%%%%%%%%%%%%%%%%%%%%
%%%%%%%%%%%%%%%%%%%%%%%%%%%%%%%%%%%%%%%%%%%%%%%%%%%%%%%%%%%%%%%%%%%%
\subhead\newsubsection 
Orthogonal polynomials on the real line
\endsubhead
Let $\mu$ be a non-negative Borel measure on $\R$ such that all
moments exist, i.e. $\int x^k\, d\mu(x)<\infty$ for all $k\in\Zp$,
and such that $\text{supp}(\mu)$ contains at least a countably
infinite number of points.
The polynomials $\{ p_n(x)\}_{n=0}^\infty$ with 
$\text{degree}(p_n)=n$ and real coefficients
are said to be orthogonal polynomials with repect to $\mu$
if
$$
\int_{-\infty}^\infty p_n(x)p_m(x)\, d\mu(x) = \de_{nm}h_n,
\qquad h_n>0.
\tag\eqname{\vgldeforthopols}
$$
This definition is equivalent to
$$
\int_{-\infty}^\infty p_n(x) x^m\, d\mu(x) = \de_{nm}g_n,
\qquad g_n\not=0, \quad 0\leq m\leq n,
\tag\eqname{\vglequivdeforthopols}
$$
or with $x^m$ replaced by any other polynomial of degree $m$.

The orthogonal polynomials are uniquely determined up to
a scalar depending on the degree $n$. There are two canonical
choices; (i) such that the leading coefficient is $1$, and then
we speak of monic orthogonal polynomials,
(ii) such that $h_n$ in \thetag{\vgldeforthopols} is
independent of $n$ and that the leading coefficient
is positive and then we assume $\mu$ normalised
by $m_0=\mu(\R)=1$ so that $h_n=1$ and we speak of orthonormal
polynomials.

Orthogonal polynomials satisfy a fundamental three-term
recurrence relation.

\proclaim{Theorem \theoremname{\thmthreetermrecurrfororthopol}}
Let $p_n$ be orthogonal polynomials, then
$$
x\, p_n(x) = A_n\, p_{n+1}(x) + B_n\, p_n(x) + C_n\, p_{n-1}(x)
$$
with $A_n,B_n,C_n\in\R$ and $A_{n-1}C_n>0$ for $n\geq 1$.
\endproclaim

\demo{Proof} $xp_n(x)$ is a polynomial of degree $n+1$,
so $xp_n(x)=\sum_{k=0}^{n+1} c_kp_k(x)$. The orthogonality
relations \thetag{\vgldeforthopols} show that
$$
h_k\, c_k = \int_{-\infty}^\infty p_n(x)\, x p_k(x)\, d\mu(x).
$$
Now $h_k>0$ and the integral is real-valued and
zero for $k= 0,\ldots,n-2$
by \thetag{\vglequivdeforthopols}. We also see that
$A_{n-1}h_n = \int p_n(x) x p_{n-1}(x) d\mu(x) = C_n h_{n-1}$,
so that $A_{n-1}C_n\geq 0$ since $h_k>0$.
It cannot be zero since $A_n\not= 0$.
\qed\enddemo

On the other hand, the three-term recursion relation together
with the initial conditions $p_{-1}(x)=0$, $p_0(x)=1$ defines
polynomial $p_n$ as polynomials of degree $n$. Then
the converse also holds, and this is commonly called
Favard's Theorem.

\proclaim{Theorem \theoremname{\thmFavard}}
{\rm (Favard)} Define polynomials
$p_n$ of degree $n$ by
$$
x\, p_n(x) = A_n\, p_{n+1}(x) + B_n\, p_n(x) + C_n\, p_{n-1}(x)
$$
with initial conditions $p_{-1}(x)=0$ and $p_0(x)=1$. Assume
$A_n,B_n,C_n\in\R$ and $A_{n-1}C_n>0$ for $n\geq 1$, then there
exists a non-negative Borel measure $\mu$ on $\R$ such that
$\{ p_n\}_{n=0}^\infty$ are orthogonal polynomials
with respect to $\mu$.
\endproclaim

\demo{Remark \theoremname{\remthreetermrecrel}} (i) Observe that
for monic polynomials $A_n=1$ and that for orthonormal
polynomials $A_{n-1}=C_n$.

(ii) Note that the orthogonality measure from Favard's Theorem
is not unique. However, under some extra conditions on the
coefficients, e.g. the polynomials are orthonormal, i.e. 
$A_{n-1}=C_n$
and the coefficients $A_n$, $B_n$ are bounded for $n\to\infty$,
the measure $\mu$ is unique. In the case mentioned we even have
that $\text{supp}(\mu)$ is compact.
\enddemo

\demo{Sketch of Proof} By rescaling we may suppose that we deal with
$$
x\, p_n(x) = a_{n+1}\, p_{n+1}(x) + b_n\, p_n(x) + a_n\, p_{n-1}(x)
$$
with initial conditions $p_{-1}(x)=0$ and $p_0(x)=1$. In the
Hilbert space $\Hi$ with orthonormal basis $\{ e_n\}_{n=0}^\infty$
we define a linear operator $J$ by
$$
J\, e_n = a_{n+1}\, e_{n+1} + b_n\, e_n + a_n\, e_{n-1},
$$
i.e. $J$ is a tridiagonal operator which are also called Jacobi
matrices. We now prove Favard's Theorem under the extra assumption
that $a_n$ and $b_n$ are bounded. This implies that $J$ is a bounded
self-adjoint operator on $\Hi$, hence the Spectral Theorem
applies and $J=\int \la dE(\la)$ for some projection valued
Borel measure $E$ on $\R$.
Define $\mu(B)=E_{e_0,e_0}(B)=\langle E(B)e_0,e_0\rangle$,
then $\mu$ is a non-negative Borel measure on $\R$. Moreover,
$\text{supp}(\mu)$ is compact, since
$\text{supp}(\mu)\subset \text{supp}(E)\subset \si(J)$ and
the spectrum of $J$ is compact.

Define the mapping
$$
\Lambda\colon \Hi\to L^2(\mu),\qquad e_n\mapsto p_n,
$$
where $L^2(\mu)=\{f\colon \R\to\R\mid \int |f(x)|^2d\mu(x)<\infty\}$
is a weighted $L^2$-space. The three-term recurrence relation
then implies $\Lambda\circ J = M\circ \Lambda$, where
$Mf(x)=xf(x)$ is the multiplication operator. Moreover, $\Lambda$
is a unitary mapping. Indeed, the polynomials are dense in
$L^2(\mu)$ since $\text{supp}(\mu)$ is compact. And from
$\Lambda J^ne_0= x^n$ we get
$$
\langle J^ne_0,J^me_0\rangle =\int \la^{n+m} \, dE_{e_0,e_0}(\la)
=\int \la^{n+m} \, d\mu(\la) = \langle x^n,x^m\rangle_{L^2(\mu)}.
$$
The unitarity implies
$$
\de_{nm}=\langle e_n,e_m\rangle = \langle \Lambda e_n,
\Lambda e_m \rangle_{L^2(\mu)} = \int p_n(x)p_m(x)\, d\mu(x).
$$

A similar, but more delicate, construction works if $J$ is
unbounded. The deficiency indices are $(0,0)$ or $(1,1)$,
in the first case the measure is still unique but might
cease to have compact support, and in the second case
the orthogonality measure depends on the self-adjoint
extension and is no longer unique.
\qed\enddemo

%%%%%%%%%%%%%%%%%%%%%%%%%%%%%%%%%%%%%%%%%%%%%%%%%%%%%%%%%%%%%%%%%%%%
%%N E W   S U B S E C T I O N%%%%%%%%%%%%%%%%%%%%%%%%%%%%%%%%%%%%%%%
%%%%%%%%%%%%%%%%%%%%%%%%%%%%%%%%%%%%%%%%%%%%%%%%%%%%%%%%%%%%%%%%%%%%
\subhead\newsubsection 
Basic hypergeometric series
\endsubhead
The series $\sum c_k$ is a {\sl basic hypergeometric series}
if $c_{k+1}/c_k$ is a rational function of $q^k$ for a base $q$.
So
$$
{{c_{k+1}}\over{c_k}} = {{(1-a_1q^k)\ldots (1-a_rq^k)}
\over{(1-b_1q^k)
\ldots (1-b_sq^k)}} {{(-q^k)^{1+s-r}}\over{1-q^{k+1}}} z,
$$
and we have the following form for basic hypergeometric series, also
known as $q$-hyper\-geo\-me\-tric series,
$$
\aligned
&{}_r\vp_s(a_1,\ldots,a_r;b_1,\dots, b_s;q,z)=
{}_r\vp_s\left( {{a_1,\ldots,a_r}\atop{b_1,\dots, b_s}};
q,z\right)\\
&= \sum_{k=0}^\infty {{(a_1,\ldots,a_r;q)_k}\over{(b_1,\ldots, 
b_s;q)_k}}
{{\bigl( (-1)^k q^{k(k-1)/2}\bigr)^{1+s-r} z^k}\over{(q;q)_k}},
\endaligned
\tag\eqname{\vgldefbasichypgeomser}
$$
with the notation for $q$-shifted factorials
$$
(a_1,\ldots,a_r;q)_k = (a_1;q)_k\ldots (a_r;q)_k,\qquad
(a;q)_k = \prod_{i=0}^{k-1} (1-aq^i).
$$
When dealing with $q$-hypergeometric series our standard assumption
on $q$ is $0<q<1$. Then the $q$-shifted factorials are also 
well-defined for $k\to\infty$, 
$(a;q)_\infty = \lim_{k\to\infty}(a;q)_k$.

Since $(q^{-n};q)_k=0$ for $k>n$ we see that the series
\thetag{\vgldefbasichypgeomser} terminates if one of the upper
parameters $a_i$ equals $q^{-n}$ for $n\in\Zp$. If one of the
lower parameters equals $q^{-N}$ for $N\in\Zp$ the series
\thetag{\vgldefbasichypgeomser} is not well-defined, unless
one of the upper parameters equals $q^{-n}$ for some
$n\in\{ 0,1,\ldots,N\}$. In case $n=N$ we follow the convention
that the series consists of the first $N+1$ terms.

The ratio test shows that for generic values of the parameters
the radius of convergence is $\infty$, $1$ or $0$ for
$r<s+1$, $r=s+1$ or $r>s+1$.

\proclaim{Theorem \theoremname{\thmqbinomialthm}}
{\rm ($q$-binomial theorem)}
${}_1\vp_0(a;-;q,z)= {\dsize{ 
{{(az;q)_\infty}\over{(z;q)_\infty}} }}$
for $|z|<1$.
\endproclaim

\demo{Proof} Let $h_a(z)={}_1\vp_0(a;-;q,z)$, which is analytic in
the unit disc. Simple calculations show
$h_a(z)-h_{aq}(z)=-azh_{aq}(z)$ and
$h_a(z)-h_a(qz)=z(1-a)h_{aq}(z)$. Eliminating $h_{aq}(z)$ leads
to
$$
h_a(z)= {{1-az}\over{1-z}}\, h_a(qz) =
{{(az;q)_n}\over{(z;q)_n}}\, h_a(q^nz).
$$
Finally, let $n\to\infty$ and use that $h_a(z)$ is continuous at
$z=0$ and $h_a(0)=1$.
\qed\enddemo

\proclaim{Corollary \theoremname{\corstothmqbinomialthm}}
{\rm (i)} ${}_1\vp_0(q^{-n};-;q,q^nz) = (z;q)_n$ for $z\in\C$,

\noindent
{\rm (ii)} $e_q(z) = {}_1\vp_0(0;-;q,z) = (z;q)_
\infty^{-1}$ for $|z|<1$,

\noindent
{\rm (iii)} $E_q(z) = {}_0\vp_0(-;-;q,-z) = (-z;q)_\infty$ 
for $z\in\C$.
\endproclaim

\demo{Remark \theoremname{\remcorstothmqbinomialthm}} We call $e_q$
and $E_q$ $q$-exponential functions. To motivate this terminology
we note that $\lim_{q\uparrow 1} e_q(z(1-q)) = e^z =
\lim_{q\uparrow 1} E_q(z(1-q))$, which follows formally from the
power series representation for $e_q$ and $E_q$.
\enddemo

\demo{Proof} Case (i) and (ii) follow from
Theorem \thmref{\thmqbinomialthm} by specialisation of $a$. 
Case (iii) follows from replacing $z$ by $z/a$ in Theorem 
\thmref{\thmqbinomialthm} and letting $a\to\infty$.
\qed\enddemo

Identities for $q$-hypergeometric series can be obtained by
playing around with the $q$-binomial Theorem \
thmref{\thmqbinomialthm}. As an example we derive Heine's 
transformation formulae for the ${}_2\vp_1$-series.

\proclaim{Theorem \theoremname{\thmHeinetransformtwophione}}
{\rm (Heine)}
$$
\align
{}_2\vp_1(a,b;c;q,z) &= {{(b,az;q)_\infty}\over{(c,z;q)_\infty}}
\, {}_2\vp_1(c/b,z;az;q,b), \quad |z|, |b|<1, \\
{}_2\vp_1(a,b;c;q,z) &= {{(c/b,bz;q)_\infty}\over{(c,z;q)_\infty}}
\, {}_2\vp_1(abz/c,b;bz;q,c/b), \quad |z|, |c/b|<1, \\
{}_2\vp_1(a,b;c;q,z) &= {{(abz/c;q)_\infty}\over{(z;q)_\infty}}
\, {}_2\vp_1(c/a,c/b;c;q,abz/c), \quad |z|, |abz/c|<1.
\endalign
$$
\endproclaim

\demo{Proof} The second equality follows from the first by applying
it twice. The third equality follows from the second on the first,
or by a threefold application of the first equality. Hence it 
suffices to prove the first equality. Write
$$
\align
{}_2\vp_1(a,b;c;q,z) &= {{(b;q)_\infty}\over{(c;q)_\infty}}
\sum_{n=0}^\infty {{(a;q)_n(cq^n;q)_\infty}\over{
(q;q)_n(bq^n;q)_\infty}} z^n \\
& = {{(b;q)_\infty}\over{(c;q)_\infty}}
\sum_{n=0}^\infty {{(a;q)_n}\over{
(q;q)_n}} z^n \sum_{m=0}^\infty {{(c/b;q)_m}\over{(q;q)_m}}
(bq^n)^m \\
\intertext{by the $q$-binomial Theorem \thmref{\thmqbinomialthm}}
&={{(b;q)_\infty}\over{(c;q)_\infty}}
\sum_{m=0}^\infty {{(c/b;q)_m}\over{(q;q)_m}} b^m
\sum_{n=0}^\infty {{(a;q)_n}\over{ (q;q)_n}} (zq^m)^n \\
&= {{(b;q)_\infty}\over{(c;q)_\infty}}
\sum_{m=0}^\infty {{(c/b;q)_m}\over{(q;q)_m}} b^m
{{(azq^m;q)_\infty}\over{(zq^m;q)_\infty}} \\
\intertext{by the $q$-binomial Theorem \thmref{\thmqbinomialthm} 
again}
&= {{(b,az;q)_\infty}\over{(c,z;q)_\infty}}
\, {}_2\vp_1(c/b,z;az;q,b).
\endalign
$$
Interchanging summations is allowed for $|z|,|b|<1$, since then
the double sum is estimated easily by a product of two
absolutely convergent series.
\qed\enddemo

\proclaim{Corollary \theoremname{\corthmHeinetransformtwophione}}
{\rm ($q$-Saalsch\"utz summation)} For $n\in\Zp$
$$
{}_3\vp_2\left( {{q^{-n}, a,b}\atop{c, q^{1-n}ab/c}};q,q\right)
= {{(c/a,c/b;q)_n}\over{(c,c/(ab);q)_n}}.
$$
\endproclaim

\demo{Proof} Consider the last of Heine's transformation
formulae in Theorem \thmref{\thmHeinetransformtwophione}.
Use the $q$-binomial Theorem \thmref{\thmqbinomialthm} for the
quotient of the infinite $q$-shifted factorials to write it as
a power series of $z$. Now compare the coefficients of $z$
on both sides. Finally, replace $a$, $b$ by $c/a$, $c/b$.
\qed\enddemo

\demo{Remark \theoremname{\rembalancedseries}} The series
$_{r+1}\vp_r(a_1,\ldots,a_{r+1};b_1,\ldots,b_r;q,z)$
is called balanced, or Saal\-sch\"utzian, if $b_1\ldots b_r =
qa_1\ldots a_{r+1}$ and $z=q$. Corollary 
\thmref{\corthmHeinetransformtwophione} shows that any
balanced ${}_3\vp_2$-series is summable.
\enddemo

%%%%%%%%%%%%%%%%%%%%%%%%%%%%%%%%%%%%%%%%%%%%%%%%%%%%%%%%%%%%%%%%%%%%
%%N E W   S U B S E C T I O N%%%%%%%%%%%%%%%%%%%%%%%%%%%%%%%%%%%%%%%
%%%%%%%%%%%%%%%%%%%%%%%%%%%%%%%%%%%%%%%%%%%%%%%%%%%%%%%%%%%%%%%%%%%%
\subhead\newsubsection 
Askey-Wilson polynomials
\endsubhead
To formulate the celebrated Askey-Wilson integral we have to
introduce the following weight function;
$$
w\bigl( {1\over 2}(z+z^{-1})\bigr) =
{{(z^2,z^{-2};q)_\infty}\over{
(az,a/z,bz,b/z,cz,c/z,dz,d/z;q)_\infty}}
$$
and we use $w(x)=w(x;a,b,c,d|q)$ to stress the dependence on the
parameters when needed.

\proclaim{Theorem \theoremname{\thmAskeyWilsonintegral}}
Assume $|a|,|b|,|c|,|d|< 1$, then
$$
\int_0^\pi w(\cos\th;a,b,c,d|q)\, d\th =
{{2\pi (abcd;q)_\infty}\over{(q,ab,ac,ad,bc,bd,cd;q)_\infty}}.
$$
\endproclaim

Theorem \thmref{\thmAskeyWilsonintegral} is the key to the
Askey-Wilson polynomials. To motivate its introduction
define the monomial function $m_k(\cos\th;a) =
(ae^{i\th},ae^{-\th};q)_k$ and then we consider
$$
\multline
{1\over{2\pi}} \int_0^\pi m_j(\cos\th;a)m_k(\cos\th;b)
w(\cos\th;a,b,c,d|q)\, d\th =
{1\over{2\pi}} \int_0^\pi
w(\cos\th;aq^j,b^k,c,d|q)\, d\th \\
= {{(ab;q)_{k+j} (ac,ad;q)_j (bc,bd;q)_k}\over{(abcd;q)_{k+j}}}
{{(abcd;q)_\infty}\over{(q,ab,ac,ad,bc,bd,cd;q)_\infty}}.
\endmultline
$$
So we can try to find orthogonal polynomials
$r_m(x)=\sum_{j=0}^m c_j m_j(a;x)$ such that
$r_m$ is orthogonal to each $m_k(x;b)$ for $0\leq k<m$. Hence
the coefficients $c_j$ have to satisfy
$$
\sum_{j=0}^m c_j {{(abq^k;q)_j}\over{(abcdq^k;q)_j}}
(ac,ad;q)_j = \de_{m,k} g_m \qquad \text{for}\ 0\leq k \leq m.
\tag\eqname{\vglAskeyWilsonorthoneeded}
$$
Take $c_j=(q^{-m},abcdq^{m-1};q)_jq^j/(q,ab,ac,ad;q)_j$, then
\thetag{\vglAskeyWilsonorthoneeded} equals
$$
{}_3\vp_2\left( {{q^{-m}, abcdq^{m-1}, abq^k}\atop{ 
ab,\  abcdq^k}};q,q
\right) =
{{(q^{k-m+1}, cd;q)_m}\over{(abcdq^k,q^{1-m}/(ab);q)_m}}
$$
by the $q$-Saalsch\"utz formula of Corollary
\thmref{\corthmHeinetransformtwophione}. This is zero
for $0\leq k<m$ and non-zero for $k=m$, and thus
$$
r_m(\cos\th;a,b,c,d|q) = {}_4\vp_3
\left( {{q^{-m},abcdq^{m-1},ae^{i\th},ae^{-i\th}}
\atop{ab,\ ac,\ ad}}; q,q\right)
\tag\eqname{\vglAskeyWilsonasfourphithreeser}
$$
are the orthogonal polynomials with respect to the
Askey-Wilson integral of Theorem \thmref{\thmAskeyWilsonintegral}.
Since the value of $g_m$ also follows from this calculation, and
since the leading coefficients of $r_m$ and $m_m(x)$ are
easily calculable, we can also calculate the norm of $r_m$
with respect to the Askey-Wilson measure.

\proclaim{Theorem \theoremname{\thmAWpolsareorthogonal}} Define
the Askey-Wilson polynomials as
$$
p_m(\cos\th;a,b,c,d|q) = a^{-m} (ab,ac,ad;q)_m\, {}_4\vp_3
\left( {{q^{-m},abcdq^{m-1},ae^{i\th},ae^{-i\th}}
\atop{ab,\ ac,\ ad}}; q,q\right),
%\tag\eqname{\vgldefAskeyWilsonpols}
$$
and assume $|a|, |b|,|c|,|d|<1$,
then $p_m(x)=p_m(x;a,b,c,d|q)$ satisfies
$$
{1\over{2\pi}} \int_0^\pi p_n(\cos\th)p_m(\cos\th)\, 
w(\cos\th)\, d\th = \de_{nm} h_n
$$
with
$$
\align
h_n &= {{ (1-q^{n-1}abcd)}\over{(1-q^{2n-1}abcd)}}
{{(q,ab,ac,ad,bc,bd,cd;q)_n}\over{(abcd;q)_n}}\, h_0 , \\
h_0 &= {{(abcd;q)_\infty}\over{(q,ab,ac,ad,bc,bd,cd;q)_\infty}}.
\endalign
$$
\endproclaim

\demo{Remark \theoremname{\remthmAWpolsareorthogonal}} The factor
in front of the ${}_4\vp_3$-series in
Theorem \thmref{\thmAWpolsareorthogonal}
is chosen such that the squared
norm is symmetric in the parameters. Since the weight function
possesses this symmetry as well, the Askey-Wilson polynomials
are symmetric in its parameters. Using the symmetry in $a$
and $b$ leads to a transformation for ${}_4\vp_3$-series,
which is Sears's transformation formula.
\enddemo

For more general values of the parameters the Askey-Wilson
polynomials defined in Theorem \thmref{\thmAWpolsareorthogonal}
are still orthogonal. The more general result can be obtained
by working with contour integration and next  use contour shifts
and residue calculus.

\proclaim{Proposition \theoremname{\propAWorthogonality}}
Let $a$, $b$, $c$ and $d$ be real and let all the pairwise 
products of $a$, $b$, $c$ and $d$ be less than $1$. Then the 
Askey-Wilson polynomials $p_n(x)= p_n(x;a,b,c,d| q)$ satisfy 
the orthogonality relations
$$
{1\over {2\pi}} \int_0^\pi p_n(\cos\theta)p_m(\cos\theta)
w(\cos\theta)\,d\theta
+\sum_k p_n(x_k)p_m(x_k)w_k = \delta_{n,m} h_n.
$$
The points $x_k$ are of the form ${1\over 2}(eq^k+e^{-1}q^{-k})$ 
for $e$ any of the parameters $a$, $b$, $c$ or $d$ with absolute 
value greater than $1$; the sum is over $k\in\Zp$ such that 
$\vert eq^k\vert >1$ and $w_k$ is the residue of 
$z\mapsto w\bigl({1\over 2}(z+z^{-1})\bigr)$ at
$z=eq^k$ minus the residue at $z=e^{-1}q^{-k}$.
\endproclaim

The orthogonality relations remain valid for complex parameters 
$a$, $b$, $c$ and $d$, if they occur in conjugate pairs. If all 
parameters have absolute value less than $1$, the Askey-Wilson 
orthogonality measure is absolutely continuous, i.e. we are in 
the situation of Theorem \thmref{\thmAWpolsareorthogonal}.

We use the notation $dm(x)=dm(x;a,b,c,d| q)$ for the normalised
orthogonality measure. So for any polynomial $p$
$$
\int_\R p(x) \, dm(x) = {1\over{h_0}} \Bigl( {1\over{2\pi}}
\int_{-1}^1 p(x)w(x) {{dx}\over{\sqrt{1-x^2}}} + \sum_k p(x_k)w_k
\Bigr).
\tag\eqname{\vglnormalisedAWmeasure}
$$

%%%%%%%%%%%%%%%%%%%%%%%%%%%%%%%%%%%%%%%%%%%%%%%%%%%%%%%%%%%%%%%%%%%%
%%N E W   S U B S E C T I O N%%%%%%%%%%%%%%%%%%%%%%%%%%%%%%%%%%%%%%%
%%%%%%%%%%%%%%%%%%%%%%%%%%%%%%%%%%%%%%%%%%%%%%%%%%%%%%%%%%%%%%%%%%%%
%NOTES AND REFERENCES%%%%%%%%%%%%%%%%%%%%%%%%%%%%%%%%%%%%%%%%%%%%%%%
\subhead Notes and references
\endsubhead
There is an enormous amount of literature available on general
orthogonal polynomials, and we only have given a very small
portion of the available results. Further introductions  can be 
found in e.g. Chihara \cite{\Chih}, Temme \cite{\Temm, Ch.~6},
Szeg\H o \cite{\Szeg}. More details on this proof of
Favard's Theorem \thmref{\thmFavard} and the
relation between orthogonal polynomials and functional analysis
can be found in Berezanski\u\i\ \cite{\Bere, Ch.~7.1}, see
also Dombrowski \cite{\Domb}. Chihara \cite{\Chih} gives another
proof of Favard's Theorem. In case the moment problem is
not determined, i.e. there exist more than one orthogonality
measure for the corresponding orthogonal polynomials,
the analysis becomes much more delicate, see Berg \cite{\Berg} and 
references given there.

The basic hypergeometric ${}_2\vp_1$-series was introduced
in 1846 by Heine. Since then research has been going on, and
a very good account of properties of basic hypergeometric series
can be found in the book \cite{\GaspR} on this subject
by Gasper and Rahman. A number of connections with other fields,
such as quantum groups, Lie algebras, number theory, statistical
mechanics and other areas in physics are known, see e.g.
Andrews \cite{\Andr} for a nice account and references in
\cite{\GaspR}.

The Askey-Wilson integral of Theorem 
\thmref{\thmAskeyWilsonintegral} is due to Askey and Wilson 
\cite{\AskeW}, who evaluated the integral by calculating residues 
and a number of summation formulas.
See also \cite{\GaspR, \S 6.1} for another proof and further
references. A nice proof of Theorem \thmref{\thmAskeyWilsonintegral}
is given by Kalnins and Miller \cite{\KalnMRMJM}
using the symmetry in $a$, $b$,
$c$ and $d$ and a suitable iteration, see \cite{\KoorTrento},
and Exercise~3.7,
for an adaptation of this method. Another recent elegant
proof of the Askey-Wilson measure is given by Berg and Ismail
\cite{\BergI}.

The Askey-Wilson polynomials, together with their finite
discrete counterpart, the so-called $q$-Racah polynomials, form
the top level of the $q$-analogue of the Askey-scheme.
A very useful compendium of properties of orthogonal polynomials
both from the Askey scheme and its $q$-analogue is given
by Koekoek and Swarttouw \cite{\KoekS}.

%%%%%%%%%%%%%%%%%%%%%%%%%%%%%%%%%%%%%%%%%%%%%%%%%%%%%%%%%%%%%%%%%%%%
%%N E W   S U B S E C T I O N%%%%%%%%%%%%%%%%%%%%%%%%%%%%%%%%%%%%%%%
%%%%%%%%%%%%%%%%%%%%%%%%%%%%%%%%%%%%%%%%%%%%%%%%%%%%%%%%%%%%%%%%%%%%
%EXERCISES%%%%%%%%%%%%%%%%%%%%%%%%%%%%%%%%%%%%%%%%%%%%%%%%%%%%%%%%%%
\subhead Exercises
\endsubhead

\item{\the\sectionno.1} Let $p_n$ be orthonormal polynomials, 
prove the Christoffel-Darboux formula
$$
\sum_{m=0}^n p_m(x)p_m(y)= {{A_n}\over{x-y}}
\Bigl( p_{n+1}(x)p_n(y)-p_n(x)p_{n+1}(y)\Bigr).
$$
What is the resulting identity for $y\to x$?

\item{\the\sectionno.2} Let $p_n$ be orthogonal polynomials. 
Prove that $p_n$ has $n$ real simple roots in the convex hull 
of $\text{supp}(\mu)$.
(Hint: Use Exercise {\the\sectionno.1} with $y=x$.)

\item{\the\sectionno.3} Suppose that $\mu$ is a finite discrete
non-negative measure with the support containing $N+1$
points. Show that we can still define a finite collection
of orthogonal polynomials $\{ p_n\}_{n=0}^N$. In this case
the Jacobi matrix is a self-adjoint $(N+1)\times(N+1)$-matrix
whose eigenvalues are the mass points.
As an example, consider the $q$-Krawtchouk polynomials
$$
K_n(x;q^\si,N;q) = {}_3\vp_2 \left( {{q^{-n},x,-q^{n-N-\si}}\atop
{q^{-N},0}};q,q\right)
$$
and prove the orthogonality relations
$$
{{q^{N+\si}}\over{(-q^\si;q)_{N+1}}} \sum_{x=0}^N
\bigl( K_nK_m\bigr) (q^{-x};q^\si,N;q) w_x(q^\si,N)
= \de_{n,m} \bigl( h_n(q^\si,N)\bigr)^{-1}
$$
with
$$
\gather
w_x(q^\si,N) = (-q^{N+\si})^x {{(q^{-N};q)_x}\over{(q;q)_x}}, \\
h_n(q^\si,N) = {{(1+q^{2n-N-\si})}\over{(-q^{n-2N-\si})^n}}
{{(-q^{-N-\si},q^{-N};q)_n}\over{(q,-q^{1-\si};q)_n}}.
\endgather
$$
Deduce from this that the dual $q$-Krawtchouk polynomials
$$
R_n(q^{-x}-q^{x-N-\si};q^\si,N;q) =
{}_3\vp_2 \left( {{q^{-n},q^{-x},-q^{x-N-\si}}\atop
{q^{-N},0}};q,q\right) = K_x(q^{-n};q^\si,N;q)
$$
satisfy the orthogonality relations
$$
{{q^{N+\si}}\over{(-q^\si;q)_{N+1}}} \sum_{x=0}^N
\bigl( R_nR_m\bigr) (q^{-x}-q^{x-N-\si};q^\si,N;q)
h_x(q^\si,N) = \de_{n,m} \bigl( w_n(q^\si,N)\bigr)^{-1}.
$$

\item{\the\sectionno.4} Derive the three-term recurrence relation 
for the Al-Salam and Chihara polynomials $s_n$, which
are Askey-Wilson polynomials with $c=d=0$, so
$s_n(x)=s_n(x;a,b|q)= p_n(x;a,b,0,0|q)$;
$$
x\, s_n(x) = s_{n+1}(x) + q^n(a+b)\, s_n(x) +
(1-q^n)(1-abq^{n-1})\, s_{n-1}(x).
$$

\item{\the\sectionno.5} Define the $q$-integral by
$$
\int_0^b f(x) \,d_qx = (1-q)b\sum_{k=0}^\infty f(bq^k)q^k,\qquad
\int_a^b f(x)\, d_qx = \int_0^b f(x)\, d_qx - \int_0^a f(x)\, d_qx,
$$
whenever $f$ is such that the series involved are convergent.
Show that $\int_0^b f(x) \,d_qx$ tends to $\int_0^b f(x) \,dx$
as $q\uparrow 1$ for Riemann integrable $f$. Show that
the inverse operator $D_q$, the $q$-derivative, is given by
$D_qf(x)=\bigl( f(x)-f(qx)\bigr)/(1-q)x$ for $x\not= 0$.

\item{\the\sectionno.6} Prove Gau\ss's summation formula;
$$
{}_2\vp_1(a,b;c;q,c/ab) = {{(c/a,c/b;q)_\infty}
\over{(c,c/ab;q)_\infty}}, \qquad |c/ab|<1.
$$
(Hint: use the first of Heine's transformation formulae.) Then 
prove Jackson's transformation formula
$$
{}_2\vp_1(a,b;c;q,z) = {{(az;q)_\infty}\over{(z;q)_\infty}}
\, {}_2\vp_2(a,c/b;c,az;q,bz), \qquad |z|<1.
$$
(Hint: develop $(b;q)_k/(c;q)_k$ in the summand on the left hand 
side using the terminating Gau\ss\ formula, i.e. $a=q^{-k}$.)

\item{\the\sectionno.7} In this exercise we sketch a proof
of the Askey-Wilson integral in Theorem 
\thmref{\thmAskeyWilsonintegral},
which is taken from \cite{\KoorTrento} which in turn is motivated
by \cite{\KalnMRMJM}. So $|a|,|b|,|c|,|d|<1$ and consider
$w_{a,b,c,d}(z) = w\bigl((z+z^{-1})/2\bigr)$ as in \S 3.3. Define
the contour integral
$$
I_{a,b,c,d} = {1\over{2\pi i}}\oint_{|z|=1} w_{a,b,c,d}(z)\, 
{{dz}\over z}.
$$
Prove that $I_{a,b,c,d}=(1-abcd)I_{aq,b,c,d}/(1-ab)(1-ac)(1-ad)$. 
Do this by giving two expressions for
$$
\oint_{|z|=1} {{w_{aq^{1/2},bq^{1/2},cq^{1/2},dq^{1/2}}}\over
{z-z^{-1}}}\, {{dz}\over z}
$$
by shifting the contour to $|z|=q^{\pm 1/2}$ and scaling back to 
$|z|=1$, and subtracting the result. Then prove as a consequence
$$
I_{a,b,c,d}={{(abcd;q)_\infty}\over{(ab,ac,ad.bc,bd,cd;q)_\infty}}
 I_{0,0,0,0}.
$$
Then prove $I_{1,q^{1/2},-1,-q^{1/2}}=1$. (Why is this choice okay?)
Now finish the proof of Theorem \thmref{\thmAskeyWilsonintegral}.

\item{\the\sectionno.8} Formulate Sears's transformation formula, 
cf. Remark \thmref{\remthmAWpolsareorthogonal}, and prove
$$
{}_3\vp_2 \left( {{q^{-n},b,c}\atop{d,\, e}};q,q\right) =
{{(de/bc;q)_n}\over{(e;q)_n}} (bc/d)^n
\, {}_3\vp_2\left( {{q^{-n},d/b,d/c}\atop{d,\, de/bc}};q,q\right).
$$

\item{\the\sectionno.9} Define
$\Gamma_q(x)= (q;q)_\infty (1-q)^{1-x}/(q^x;q)_\infty$, show that 
the first of Heine's transformation formulae of
Theorem \thmref{\thmHeinetransformtwophione} can be written as
$$
{}_2\vp_1(q^a,q^b;q^c;q,z) = {{\Gamma_q(c)}\over{
\Gamma_q(b)\Gamma_q(c-b)}} \, \int_0^1 t^{b-1} 
{{(tzq^a,tq;q)_\infty}\over{(tz,tq^{c-b};q)_\infty}}\, d_qt.
$$
Assuming $\Gamma_q(x)\to \Gamma_x$ as $q\uparrow 1$, what is the
corresponding limit as $q\uparrow 1$? See the $q$-binomial theorem
as the analogue of the binomial theorem
$\sum_{k=0}^\infty (a)_k z^k/k! = (1-z)^{-a}$,
where $(a)_k=a(a+1)\ldots (a+k-1)=\Gamma(a+k)/\Gamma(a)$ is a
shifted factorial, or Pochhammer symbol.

%%%%%%%%%%%%%%%%%%%%%%%%%%%%%%%%%%%%%%%%%%%%%%%%%%%%%%%%%%%%%%%%%%%%
%%N E W  S E C T I O N%%%%%%%%%%%%%%%%%%%%%%%%%%%%%%%%%%%%%%%%%%%%%%
%%%%%%%%%%%%%%%%%%%%%%%%%%%%%%%%%%%%%%%%%%%%%%%%%%%%%%%%%%%%%%%%%%%%
\newpage

\head\newsection Quantum subgroups and the Haar functional\endhead

We discuss a way to consider the analogues of bi-$K$-invariant 
functions on $SU(2)$ for arbitrary one-parameter subgroups $K$. 
In the group case it is not of much importance which $K$ is 
chosen, since these groups are all conjugated. However, in the 
quantum group case it makes
a difference. Of particular interest is the Haar functional on
bi-$K$-invariant functions. We also start investigating
the analogues of functions
on $SU(2)$ behaving like $f(kgh)= \chi_K(k)\chi_H(h)f(g)$ for
possibly different one-parameter subgroups $K$ and $H$ of $SU(2)$
having characters $\chi_K$ respectively $\chi_H$. For
$\chi_K\equiv 1\equiv \chi_H$ the functions live on 
$K\backslash SU(2)/H$.

%%%%%%%%%%%%%%%%%%%%%%%%%%%%%%%%%%%%%%%%%%%%%%%%%%%%%%%%%%%%%%%%%%%%
%%N E W   S U B S E C T I O N%%%%%%%%%%%%%%%%%%%%%%%%%%%%%%%%%%%%%%%
%%%%%%%%%%%%%%%%%%%%%%%%%%%%%%%%%%%%%%%%%%%%%%%%%%%%%%%%%%%%%%%%%%%%
\subhead\newsubsection 
Generalised matrix elements
\endsubhead
Recall from Proposition \thmref{\propclassgroupprimelts} that
the twisted primitive elements are the analogues of the Lie algebra.
In the group case we have that $f\in\PG$ is
left and right $K$-invariant for the one-parameter group
$K=\exp(tX)$ if and only if $X.f=0=f.X$ with the action of
$U({\frak g})$ on $\PG$ as in Example \thmref{\exampairingGandUg}.
Such functions are called spherical.

For the quantum $SL(2,\C)$ group we already know
that
$$
(D-A). t^l_{ij} = \sum_{k=-l}^l t^l_{ik} \, t^l_{kj}(D-A)
= (q^j-q^{-j})t^l_{ij}
$$
and similarly $t^l_{ij}.(D-A)=(q^i-q^{-i})t^l_{ij}$. So
for $l\in\Zp$ we may consider $t^l_{00}$ as the spherical
functions on $\A$ with respect to the subgroup `generated'
by $D-A$. We now want to do this for more general twisted
primitive elements.
The $\ast$-structure of the real form $\Usu$ and $\Asu$
is needed.

For $\si\in\R$ we define
$$
X_\si = iq^{1\over 2}B-iq^{-{1\over 2}}C - {{q^\si-q^{-\si}}
\over{q-q^{-1}}}(A-D) \in\Usu .
\tag\eqname{\vgldsefXsigma}
$$
We also define
$$
X_\infty = D-A = \lim_{\si\to\infty} (q^{-1}-q)q^\si X_\si =
\lim_{\si\to -\infty} (q-q^{-1})q^{-\si} X_\si,
\tag\eqname{\vgldefXoneindig}
$$
recovering the case of the previous paragraph as a special case.
{}From \thetag{\vglHopfstructureUq} and Theorem
\thmref{\thmstarstructureonA}(ii)
$$
\Delta(X_\si) = A\otimes X_\si + X_\si\otimes D, \qquad 
S(X_\si)=-X_\si, \qquad (X_\si A)^\ast = X_\si A.
\tag\eqname{\vglpropXsigma}
$$
Note that $X_\si A$ is self-adjoint, so we consider it as the
analogue of the element $\pmatrix \si&i\\-i&-\si\endpmatrix$
in $i{\frak{su}}(2)\subset {\frak{sl}}(2,\C)$ instead of an
element of the real form ${\frak{su}}(2)$.
Instead of considering $X_\si$ we can take a self-adjoint
element in the space of primitive elements, but we need this
definition in order to be able to prove certain multiplicative
properties in the dual algebra later.

Fundamental for what follows is that we can explicitly calculate
the eigenvectors and eigenvalues of each self-adjoint
matrix $t^l(X_\si A)$.

\proclaim{Proposition \theoremname{\propspectrumdualqK}}
The self-adjoint operator $t^l(X_\si A)$ has an orthonormal basis 
of eigenvectors $v^{l,j}(\si)=\sum_{n=-l}^l v^{l,j}_n(\si)e^l_n$
corresponding to the eigenvalue
$$
\la_j(\si ) = {{q^{-2j-\si}-q^{\si+2j}+q^\si - q^{-\si}}
\over{q-q^{-1}}}, \quad j=-l,-l+1,\ldots,l.
$$
The coefficients $v^{l,j}_n(\si)$ are explicitly known by
$$
\multline
v^{l,j}_n(\si) = \, C^{l,j}(\si ) i^{n-l} q^{\si (l-n)}
q^{{1\over 2}(l-n)(l-n-1)}
\left( {{(q^{4l};q^{-2})_{l-n}}\over{(q^2;q^2)_{l-n}}}
\right)^{1/2} \\
\times R_{l-n}(q^{2j-2l}-q^{-2j-2l-2\si};q^{2\si},2l;q^2),
\endmultline
$$
where $R_{l-n}$ is a dual $q$-Krawtchouk polynomial,
cf. Exercise ~3.3, and the constant
is given by
$$
C^{l,j}(\si) = q^{l+j} \left[ {{2l}\atop{l-j}}\right]_{q^2}^{1/2}
\left( {{1+q^{-4j-2\si}}\over{1+q^{-2\si}}}\right)^{1/2}
\Bigl( (-q^{2-2\si};q^2)_{l-j} (-q^{2+2\si};q^2)_{l+j}\Bigr)^{-1/2}.
$$
\endproclaim

\demo{Proof} The relation $t^l(X_\si A)\sum_{m=-l}^l c_me_m^l=
\la \sum_{m=-l}^l c_me_m^l$ leads to a (finite) three-term 
recurrence relation for the coefficients $c_m$ by
\thetag{\vgldsefXsigma} and Theorem \thmref{\thmunitaryrtepsofUsu}.
Comparison with the three-term recurrence relation for the
dual $q$-Krawtchouk polynomials $R_n(y)=R_n(y;q^\si,N;q)$
$$
y\, R_n(y) = (1-q^{n-N})\, R_{n+1}(y) + (q^{-N}-q^{-N-\si})q^n\, 
R_n(y) -(1-q^n)q^{-N-\si}\, R_{n-1}(y)
$$
for $y=q^{-x}-q^{x-N-\si}$, $x\in\{0,\ldots,N\}$, gives the value
for the coefficients and for the eigenvalues. The normalisation
follows from the orthogonality relations for the $q$-Krawtchouk
polynomials, cf. Exercise~3.3.
\qed\enddemo

Since the vectors $v^{l,j}(\si)$ are orthonormal in a finite 
dimensional space, we have the orthogonality relations
$$
\sum_{n=-l}^l v^{l,j}_n(\si) \overline{v^{l,i}_n(\si )} = 
\delta_{i,j}, \qquad
\sum_{i=-l}^l v^{l,i}_m(\si) \overline{v^{l,i}_n(\si )} = 
\delta_{n,m}.
\tag\eqname{\vglorthorelvector}
$$
The first part of \thetag{\vglorthorelvector} is equivalent to the
orthogonality relations
for the $q$-Kraw\-tchouk polynomials, and
the second part of \thetag{\vglorthorelvector} is equivalent to the
orthogonality relations
for the dual $q$-Krawtchouk polynomials, cf. Exercise~3.3.

With this orthonormal basis for $\C^{2l+1}$ we can define
generalised matrix elements.

\proclaim{Lemma \theoremname{\lemintrogenmatrixelts}} In $\Asu$ we
define
$$
a^l_{i,j}(\ta,\si) \bigl(X\bigr) = \langle t^l(X)v^{l,j}(\si),
v^{l,i}(\ta)\rangle, \qquad \si,\ta\in\R,\ i,j=-l,-l+1,\ldots,l,
$$
then
$$
a^l_{i,j}(\ta,\si) = \sum_{n,m=-l}^l v^{l,j}_m(\si) 
\overline{v^{l,i}_n(\ta)} t^l_{n,m} \in\Asu,
$$
and
$$
(X_\si A).a^l_{i,j}(\ta,\si)=\la_j(\si)\, a^l_{i,j}(\ta,\si)
\quad\text{and}\quad
a^l_{i,j}(\ta,\si).(X_\ta A) = \la_i(\ta)\, a^l_{i,j}(\ta,\si).
$$
\endproclaim

\demo{Proof} The first statement is obvious from
Proposition \thmref{\propspectrumdualqK}
and to prove the last we observe that $\langle Y,X.\xi\rangle=
\langle YX,\xi\rangle$, so that for all $Y\in\Usu$
$$
\multline
\langle Y,(X_\si A).a^l_{i,j}(\ta,\si)\rangle =
\langle t^l(Y)t^l(X_\si A)v^{l,j}(\si),
v^{l,i}(\ta)\rangle \\
= \la_j(\si)
\langle t^l(Y)v^{l,j}(\si),
v^{l,i}(\ta)\rangle = \la_j(\si)
\langle Y,a^l_{i,j}(\ta,\si)\rangle.
\endmultline
$$
The other case is treated similarly using $(X_\ta A)^\ast=X_\ta A$.
\qed\enddemo

We want the action of $X_\si$ and $X_\ta$ to be more
symmetric than in Lemma \thmref{\lemintrogenmatrixelts}.
For this we define $b^l_{i,j}(\ta,\si)=A.a^l_{i,j}(\ta,\si)$, then
$$
X_\si .b^l_{i,j}(\ta,\si)=\la_j(\si)\, D.b^l_{i,j}(\ta,\si)
\quad\text{and}\quad
b^l_{i,j}(\ta,\si).X_\ta = \la_i(\ta)\, b^l_{i,j}(\ta,\si).D
\tag\eqname{\vglbiinvariancedefgenmatelt}
$$
since $(X.\xi).Y=X.(\xi.Y)$ by the coassociativity of $\De$.
Explicitly,
$$
b^l_{i,j}(\ta,\si) = \sum_{n,m=-l}^l v^{l,j}_m(\si) 
\overline{v^{l,i}_n(\ta)} q^{-m}\, t^l_{n,m}\ \in \Asu .
\tag\eqname{\vglexprdefgenmateltinstanmatelt}
$$
The case $l=1/2$ is of particular interest and we write
$$
b^{1/2}(\ta,\si) =
{1\over{\sqrt{(1+q^{2\si})(1+q^{2\ta})} }}
\pmatrix \al_{\ta,\si},& \be_{\ta,\si}\\
\ga_{\ta,\si},&\de_{\ta,\si} \endpmatrix,
$$
or explicitly
$$
\align
\al_{\ta,\si} &= q^{1/2}\al-iq^{\si-1/2}\be +iq^{\ta+1/2}\ga +
q^{\si+\ta-1/2}\de, \\
\be_{\ta,\si} &= -q^{\si+1/2}\al-iq^{-1/2}\be -iq^{\si+\ta+1/2}
\ga + q^{\ta-1/2}\de, \\
\ga_{\ta,\si} &= -q^{\ta+1/2}\al+iq^{\ta+\si-1/2}\be +iq^{1/2}\ga +
q^{\si-1/2}\de, \\
\de_{\ta,\si} &= q^{\ta+\si+1/2}\al+iq^{\ta-1/2}\be -
iq^{\si+1/2}\ga + q^{-1/2}\de.
\endalign
$$

The following two propositions are fundamental in giving
explicit formulas for $b^l_{i,j}(\ta,\si)$ and hence for
the generalised matrix elements $a^l_{i,j}(\ta,\si)$.

\proclaim{Definition \theoremname{\defsitasphericalelts}}
$\xi\in\Asu$ is a {\bf $(\ta,\si)$-spherical element} if
$$
X_\si .\xi=0
\quad\text{and}\quad
\xi.X_\ta = 0.
$$
\endproclaim

\proclaim{Proposition \theoremname{\propeigdefmatelt}}
{\rm (i)} Let $\xi\in\Asu$ be a $(\ta,\si)$-spherical element
and let $\et\in\Asu$ satisfy
$$
X_\si .\et=\la\, D.\eta
\quad\text{and}\quad
\et.X_\ta = \mu\, \eta.D.
\tag\eqname{\vglgenbiinvariancet}
$$
for $\la,\mu\in\C$. Then $\et\xi$ satisfies
\thetag{\vglgenbiinvariancet} for the same $\la$, $\mu$.
Moreover, if $\la,\mu\in\R$,
then $\et^\ast\et$ is a $(\ta,\si)$-spherical element.
In particular, the space of $(\ta,\si)$-spherical element
forms a $\ast$-subalgebra of $\Asu$.

\noindent
{\rm (ii)} If $\et\in\text{span}(t^l_{nm})$ satisfies
\thetag{\vglgenbiinvariancet} for
arbitrary $\la, \mu\in\C$ and $\et$ is non-zero, then
$\la= \la_j(\si)$, $\mu=\la_i(\ta)$ for some
$i,j\in\{-l,-l+1,\ldots,l\}$ and $\et$ is a multiple of
$b^l_{i,j}(\ta,\si)$.
\endproclaim

\demo{Proof} To prove (i) we first note that by
Definition \thmref{\defdualHopfalgebras} and Proposition
\thmref{\propactionUonA} we have in general for Hopf algebras
in duality
$$
u.(ab) = \sum_{(u)} \bigl( u_{(1)}.a\bigr)\bigl( u_{(2)}.b\bigr),
\qquad \De(u) = \sum_{(u)} u_{(1)}\otimes u_{(2)},
\tag\eqname{\vglactionUonproductinA}
$$
with the notation of \thetag{\vglSweedlernot}.

First consider $\et\xi$, then
$$
X_\si .(\eta\xi) = (A.\eta)(X_\si .\xi)+(X_\si .\eta)(D.\xi) =
\la (D.\eta)(D.\xi) = \la\, D.(\eta\xi),
$$
since $X_\si$ is twisted primitive and $D$ group like.
Similarly we prove $(\et\xi).X_\ta = \mu\, (\eta\xi).D$.

To prove the other statement of (i) we observe that in general
for Hopf $\ast$-algebras in duality the left action satisfies
$u.a^\ast = (S(u)^\ast.a)^\ast$. Now
proceed as before to obtain
$$
X_\si .(\et^\ast\et)= (A.\et^\ast)(X_\si .\et) + (X_\si .\et^\ast)
(D.\et) = (\la -\bar\la) (D.\et)^\ast (D.\et),
$$
since $S(A)^\ast=D$ and $S(X_\si)^\ast=-X_\si$.
This yields zero for $\la\in\R$. Similarly we prove
$(\et^\ast\et).X_\ta=0$ for real $\mu$.

The proof of (ii) follows immediately from the linear basis
of $\Asu$ given in Theorem \thmref{\thmdescripAasdualHAofU}
and the multiplicity of each eigenvalue of $t^l(X_\si A)$
being one.
\qed\enddemo

\proclaim{Proposition \theoremname{\propshifteigenval}}
Let $\et$ satisfy \thetag{\vglgenbiinvariancet} with 
$\la=\la_j(\si)$ and
$\mu=\la_i(\ta)$, then

\noindent
{\rm (i)} $\al_{\ta+2i,\si+2j}\eta$ satisfies 
\thetag{\vglgenbiinvariancet} with
$\la=\la_{j-1/2}(\si)$ and $\mu=\la_{i-1/2}(\ta)$,

\noindent
{\rm (ii)} $\be_{\ta+2i,\si+2j}\eta$ satisfies 
\thetag{\vglgenbiinvariancet} with
$\la=\la_{j+1/2}(\si)$ and $\mu=\la_{i-1/2}(\ta)$,

\noindent
{\rm (iii)} $\ga_{\ta+2i,\si+2j}\eta$ satisfies 
\thetag{\vglgenbiinvariancet} with
$\la=\la_{j-1/2}(\si)$ and $\mu=\la_{i+1/2}(\ta)$,

\noindent
{\rm (iv)} $\de_{\ta+2i,\si+2j}\eta$ satisfies 
\thetag{\vglgenbiinvariancet} with
$\la=\la_{j+1/2}(\si)$ and $\mu=\la_{i+1/2}(\ta)$.
\endproclaim

\demo{Remark \theoremname{\rempropshifteigenval}} Proposition
\thmref{\propshifteigenval} gives a way to express the corner
elements, i.e. $b^l_{ij}(\ta,\si)$ with $l=\max(|i|,|j|)$,
in terms of products of $\al_{\ta,\si}$, $\be_{\ta,\si}$,
$\ga_{\ta,\si}$, $\de_{\ta,\si}$, but with integer shifts in the
parameters $\ta$ and $\si$.
\enddemo

\demo{Sketch of proof} Take $\xi\in\Asu$ arbitrary for the moment, 
then as in the proof of Proposition \thmref{\propeigdefmatelt},
$$
X_\si.(\xi\et) =
\bigl[ \la\, A.\xi + X_\si.\xi\bigr] (D.\et), \quad
(\xi\et).X_\ta  =
\bigl[ \mu\, \xi.A + \xi.X_\ta\bigr] (\et.D).
$$
So, if $\xi$ satisfies
$$
\bigl[ \la\, A.\xi + X_\si.\xi\bigr] = \la_1\, D.\xi, \qquad
\bigl[ \mu\, \xi.A + \xi.X_\ta\bigr] = \mu_1\, \xi.D
\tag\eqname{\vglrequironxi}
$$
for some $\la_1,\mu_1\in\C$, then we get
$$
X_\si.(\xi\et)  = \la_1\, D.(\xi\et), \qquad \text{and} \qquad
(\xi\et).X_\ta = \mu_1\, (\xi\et).D.
$$
Next we assume
$\xi = a\al+b\be+c\ga+d\de$
for some complex constants $a$, $b$, $c$ and $d$.
Using \thetag{\vglexplactioUonAongen} we see that the first
requirement of \thetag{\vglrequironxi} leads to
$M \pmatrix a\\b\endpmatrix = \pmatrix 0\\0\endpmatrix$,
$M \pmatrix c\\d\endpmatrix = \pmatrix 0\\0\endpmatrix$
for a $2\times 2$-matrix depending on $\la$, $\la_1$ and $\si$.
{}From $\text{det}(M)=0$ it follows that for $\la=\la_j(\si)$
we have $\la_1=\la_{j\pm 1/2}$. Similarly, the second
requirement of \thetag{\vglrequironxi} leads to linear
equations of a similar form and again a determinant
condition implies $\mu=\la_i(\ta)$
and $\mu_1=\la_{i\pm1/2}(\ta)$. Combining each of the cases
leads to a solution for $a$, $b$, $c$ and $d$ up to
a scalar.
\qed\enddemo

Since $\la_i(\si)=0\Leftrightarrow i=0$
we see from Proposition \thmref{\propeigdefmatelt}(ii)
that there is no $(\ta,\si)$-spherical element in
$\text{span}(t^{1/2}_{nm})$ and that there is a one-dimensional
space of $(\ta,\si)$-spherical elements in
$\text{span}(t^{1}_{nm})$. Take such an element and keep only the
non-constant terms to find the following  $(\ta,\si)$-spherical 
element;
$$
\multline
\rts = {1\over 2}\bigl( \al^2+\de^2 +q\ga^2 +q^{-1}\be^2
+ i(q^{-\si}-q^\si )(q\de\ga +\be\al )\\
 -i(q^{-\ta}-q^\ta)(\de\be +q\ga\al)
+ (q^{-\si}-q^\si)(q^{-\ta}-q^\ta)\be\ga \bigr).
\endmultline
$$

\proclaim{Proposition \theoremname{\propsphericaleltsalgebra}}
The $\ast$-subalgebra of $\Asu$ of $(\ta,\si)$-spherical elements 
is generated by the self-adjoint element $\rts$. Moreover,
$$
\align
\be_{\ta+1,\si-1}\ga_{\ta,\si}&=
2q^{\ta+\si}\rts -q^{2\si-1}-q^{2\ta+1}, \\
\ga_{\ta-1,\si+1}\be_{\ta,\si}&=
2q^{\ta+\si}\rts-q^{2\si+1}-q^{2\ta-1}, \\
\al_{\ta+1,\si+1}\de_{\ta,\si}&=
2q^{\ta+\si+1}\rts+1+q^{2\si+2\ta+2}, \\
\de_{\ta-1,\si-1}\al_{\ta,\si}&=
2q^{\ta+\si-1}\rts+1+q^{2\si+2\ta-2}.
\endalign
$$
\endproclaim

\demo{Proof} Indeed, $\rts^\ast=\rts$ for the $\ast$-operator
of Theorem \thmref{\thmstarstructureonA}(ii). From Proposition
\thmref{\propeigdefmatelt}(ii) and $\la_i(\si)=0$ if and only
if $i=0$, it follows that a linear basis for the
space of $(\ta,\si)$-spherical elements is given
by $b^l_{00}(\ta,\si)$, $\l\in\Zp$. From
Lemma \thmref{\lemmaCGCforUsu} we see that $\rts^l$ has a non-zero
component in $\text{span}(t^l_{nm})$. Since $\text{span}(t^l_{nm})$
is invariant under the left action of $X_\si$ and under the
right action of $X_\ta$ it follows by induction on $l$ that
$b^l_{00}(\ta,\si)$ is a polynomial of degree $l$ in $\rts$.

{}From Proposition \thmref{\propshifteigenval} it follows that
each product on the left hand side is a
$(\ta,\si)$-spherical element, and by considering the degrees
it has to be a polynomial of degree $1$ in $\rts$ by the first
part. It remains to calculate the leading and constant coefficient
by comparing the coefficients of $\al^2$ and the unit $1$ on 
both sides.
\qed\enddemo

\proclaim{Corollary \theoremname{\corcommrst}}
The following relations hold in $\A$;
$$
\gather
\al_{\ta,\si}\rts=\rho_{\ta-1,\si-1}\al_{\ta,\si}, \qquad
\be_{\ta,\si}\rts=\rho_{\ta-1,\si+1}\be_{\ta,\si}, \\
\ga_{\ta,\si}\rts=\rho_{\ta+1,\si-1}\ga_{\ta,\si}, \qquad
\de_{\ta,\si}\rts=\rho_{\ta+1,\si+1}\de_{\ta,\si}.
\endgather
$$
\endproclaim

\demo{Proof} The
proofs of these statements are all similar. To prove the first, 
multiply the last equation of Proposition 
\thmref{\propsphericaleltsalgebra} by $\al_{\ta,\si}$
and use the third equation in the left hand side. Cancelling 
terms proves the first statement.
\qed\enddemo

\demo{Remark \theoremname{\rempropsphericaleltsalgebra}} We
define $\rti$, $\ris$ and $\rii$ as limit cases of $\rts$ as 
follows;
$$
\ris = \lim_{\ta\to\infty} 2q^{\si+\ta-1}\rts, \qquad
\rti = \lim_{\si\to\infty} 2q^{\si+\ta-1}\rts,
$$
and $\rii=\lim_{\si\to\infty} \ris = \lim_{\ta\to\infty} \rti =
-\ga\ga^\ast$. Since we take $0<q<1$, this is just replacing 
$q^\infty=0$.
The corresponding factorisations of Proposition
\thmref{\propsphericaleltsalgebra} remain valid under these limit
transitions, since we can take limits in $\al_{\ta,\si}$, etcetera.
In case $\si=\ta\to\infty$, Proposition
\thmref{\propsphericaleltsalgebra} reduces to part of the relations
\thetag{\vglcommrelAq} in $\Asu$.
\enddemo

%%%%%%%%%%%%%%%%%%%%%%%%%%%%%%%%%%%%%%%%%%%%%%%%%%%%%%%%%%%%%%%%%%%%
%%N E W   S U B S E C T I O N%%%%%%%%%%%%%%%%%%%%%%%%%%%%%%%%%%%%%%%
%%%%%%%%%%%%%%%%%%%%%%%%%%%%%%%%%%%%%%%%%%%%%%%%%%%%%%%%%%%%%%%%%%%%
\subhead\newsubsection 
The Haar functional restricted to $(\ta,\si)$-spherical elements
\endsubhead
Recall that the invariant functional or Haar functional is a 
positive linear functional on $\Asu$, which is defined by 
$h(t^l_{mn})=\de_{l0}$.

\proclaim{Lemma \theoremname{\lemHaaronbasis}}
$h(\de^l\ga^m\be^n)=h(\al^l\ga^m\be^n)=0$ unless $l=0$ and $n=m$
and
$$
h(\ga^m\be^m)= (-q)^m {{1-q^2}\over{1-q^{2m+2}}}.
$$
\endproclaim

\demo{Proof} The invariance of $h$ implies 
$h(a.X)=h(a)\ep(X)=h(X.a)$.
Use this with $X=A$, and $X=D$ (so that the left and right action 
is a homomorphism because $A$ and $D$ are group like),
\thetag{\vglexplactioUonAongen} and $a=\de^l\ga^m\be^n$ and
$a=\al^l\ga^m\be^n$ to see that the Haar functional applied to 
these basis elements gives zero unless $l=0$, $m=n$.

In this case we take $X=C$ and we calculate, using the analogue of
\thetag{\vglactionUonproductinA} for the right action,
$$
\align
\de\ga(-q^{-1}\be\ga)^n.C &=
\bigl(\de\ga(-q^{-1}\be\ga)^{n-1}.A\bigr)\bigl( 
-q^{-1}\be\ga.C\bigr) +
\bigl(\de\ga(-q^{-1}\be\ga)^{n-1}.C\bigr)\bigl( 
-q^{-1}\be\ga.D\bigr)\\
&=\bigl(q^{-1}\de\ga(-q^{-1}\be\ga)^{n-1}\bigr)\bigl( 
-q^{-1/2}\be\al\bigr)+
\bigl(\de\ga(-q^{-1}\be\ga)^{n-1}.C\bigr)\bigl( 
-q^{-1}\be\ga\bigr)\\
&=q^{-1/2-2n}(-q^{-1}\be\ga)^n\bigl(1+q^{-1}\be\ga\bigr) +
\bigl(\de\ga(-q^{-1}\be\ga)^{n-1}.C\bigr)\bigl( 
-q^{-1}\be\ga\bigr)\\
\intertext{and together with the initial condition for $n=0$ we get}
&={{q^{-1/2-2n}}\over{1-q^2}}\bigl( (1-q^{2n+2})(-q^{-1}\be\ga)^n
- (1-q^{2n+4})(-q^{-1}\be\ga)^{n+1}\bigr).
\endalign
$$
Apply $h$ and $\ep(C)=0$ to get a two-term recurrence relation, 
which is uniquely solved with the initial condition $h(1)=1$.
\qed\enddemo

\proclaim{Corollary \theoremname{\corlemHaaronbasis}}
Let $D$ be the self-adjoint positive
diagonal operator $\Hi$ with orthonormal
basis $\{e_n\}_{n=0}^\infty$ defined by $D\colon e_n\mapsto 
q^{2n}e_n$,
$$
h(a) = {{(1-q^2)}\over{2\pi}}
\int_0^{2\pi} \hbox{\rm{tr}}\bigl( D\pi^\infty_\phi(a)\bigr) 
\, d\phi,\qquad a\in\A
$$
\endproclaim

\demo{Proof} Check that the right hand side coincides with $h$
for $a=\de^l\ga^m\be^n$ or $a=\al^l\ga^m\be^n$.
\qed\enddemo

\demo{Remark \theoremname{\remcorlemHaaronbasis}}
The trace operation in Corollary \thmref{\corlemHaaronbasis}
is well-defined due to the
appearance of $D$. Since $D^{1/2}$ is a Hilbert-Schmidt operator, 
so is $\pi^\infty_\phi(a)D^{1/2}$, since $\pi^\infty_\phi(a)$ is 
a bounded operator.
The trace of the product of two Hilbert-Schmidt
operators is well-defined. Moreover,
$\hbox{\rm{tr}}\bigl( D\pi^\infty_\phi(a)\bigr)=
\hbox{\rm{tr}}\bigl(\pi^\infty_\phi(a)D\bigr)$ and
it is independent of the choice of
the basis. The trace can be estimated by the product of the
Hilbert-Schmidt norms of $D^{1/2}$ and $\pi^\infty_\phi(a)D^{1/2}$,
i.e. $\hbox{\rm{tr}}\bigl(D\pi^\infty_\phi(a)\bigr)|\leq
\parallel a\parallel /(1-q^2)$, so that
the function in Corollary \thmref{\corlemHaaronbasis}
is integrable. Here $\parallel a\parallel$, $a\in\Asu$, denotes
the supremum of the operator norm $\pi^\infty_\th(a)$ as
$\th$ varies, where we use
the notation of Theorem \thmref{\thmirredstarrepsofAsu}.
\enddemo

The next theorem is the key in determining the generalised
matrix elements of the previous section in terms of
$q$-special functions.

\proclaim{Theorem \theoremname{\thmHaarongeneralspherelt}}
The Haar functional
on the subalgebra generated by the self-adjoint element
$\rts$ is given by
$$
h(p(\rts)) = \int_{\R} p(x)\, dm(x;a,b,c,d| q^2)
$$
for any polynomial $p$. Here $a=-q^{\si +\ta +1}$,
$b=-q^{-\si - \ta +1}$, $c=q^{\si -\ta +1}$, $d=q^{-\si +\ta +1}$ 
and $dm(x;a,b,c,d| q^2)$ denotes the normalised Askey-Wilson 
measure.
\endproclaim

The proof is based on Corollary \thmref{\corlemHaaronbasis} and
the spectral theory of Jacobi matrices and their connection
with orthonormal polynomials. The proof of the general statement
is of some computational complexity. We first consider the
case $\si\to\infty$ of Theorem \thmref{\thmHaarongeneralspherelt}.
This limit case is also needed in the proof of the general
statement. We then illustrate the
procedure for the general case in a less computational case, namely
the Haar functional restricted to the co-central elements.
In these two cases the calculations are much simpler,
and the proof of the general case is contained in the exercises.

The limiting case $\si\to\infty$ of
Theorem \thmref{\thmHaarongeneralspherelt} is phrased in terms
of the $q$-integral introduced in Exercise~3.5.

\proclaim{Theorem \theoremname{\thmHaaronspecialspherelt}}
The Haar functional
on the subalgebra generated by the self-adjoint element
$\rti$ is given by the $q$-integral
$$
h\bigl(p(\rti)\bigr) = {1\over{1+q^{2\ta}}}
\int_{-1}^{q^{2\ta}} p(x) \, d_{q^2}x
$$
for any polynomial $p$.
\endproclaim

To prove Theorem \thmref{\thmHaaronspecialspherelt}
from a spectral analysis of $\pi^\infty_\th(\rti)$ we
use the following proposition.

\proclaim{Proposition \theoremname{\propeigvectpirti}}
$\Hi$ has an orthogonal basis of
eigenvectors $v_\la^\th$, where $\la=-q^{2k}$, $k\in\Zp$, and
$\la = q^{2\ta + 2k}$,
$k\in\Zp$, for the eigenvalue $\la$ of the self-adjoint operator
$\pi^\infty_\th(\rti)$. The squared norm is given by
$$
\align
\langle v_\la^\th,v_\la^\th\rangle &= q^{-2k} (q^2;q^2)_k
(-q^{2-2\ta};q^2)_k (-q^{2\ta};q^2)_\infty,\qquad \la=-q^{2k},\\
\langle v_\la^\th,v_\la^\th\rangle &= q^{-2k} (q^2;q^2)_k
(-q^{2+2\ta};q^2)_k (-q^{-2\ta};q^2)_\infty,\qquad \la=q^{2\ta+2k}.
\endalign
$$
Moreover, $v_\la^\th=\sum_{n=0}^\infty i^ne^{in\th} p_n(\la)\, e_n$
with the polynomial $p_n(\la)$ defined by
$$
\aligned
p_n(\la) &= {{q^{-n\ta} q^{n(n-1)/2}}\over{\sqrt{(q^2;q^2)_n}}}
\, {}_2\vp_1(q^{-2n}, q^{2\ta}/\la;0;q;-q^2\la) \\
&= {{(-1)^nq^{n\ta} q^{n(n-1)/2}}\over{\sqrt{(q^2;q^2)_n}}}
\, {}_2\vp_1(q^{-2n}, -1/\la;0;q;q^{2-2\ta}\la).
\endaligned
\tag\eqname{\vglexppnl}
$$
\endproclaim

\demo{Remark \theoremname{\rempropeigvectpirti}} The basis 
described in Proposition \thmref{\propeigvectpirti} induces an 
orthogonal decomposition of the representation space 
$\Hi = V_1^\th\oplus V_2^\th$, where
$V_1^\th$ is the subspace with basis $v^\th_{-q^{2k}}$, $k\in\Zp$, 
and $V_2^\th$ is the subspace with basis $v^\th_{q^{2\ta+2k}}$, 
$k\in\Zp$.
\enddemo

\demo{Sketch of Proof} Consider the Al-Salam and Carlitz 
polynomials
$$
\align
U_n^{(a)}(x;q) &= (-a)^n q^{n(n-1)/2} \, {}_2\vp_1
(q^{-n},x^{-1};0;q,qx/a) \\
&= (-1)^n q^{n(n-1)/2} \, {}_2\vp_1(q^{-n},a/x;0;q,qx).
\endalign
$$
The equality follows from the second of Heine's transformation 
formulae of Theorem \thmref{\thmHeinetransformtwophione}. 
Take $b=q^{-n}$ and reverse the order of summation on the right 
hand side before letting $c\to 0$. The Al-Salam and Carlitz 
polynomials are orthogonal polynomials satisfying the three-term 
recurrence relation
$$
x\, U^{(a)}_n(x;q) =  U^{(a)}_{n+1}(x;q) + (a+1)\,  U^{(a)}_n(x;q)
-aq^{n-1}(1-q^n)\,  U^{(a)}_{n-1}(x;q).
$$
Compare this with the action of $\rti$ in $\Hi$,
$$
\pi^\infty_\th(\rti)\, e_n = iq^{\ta+n}e^{i\th}\sqrt{1-q^{2n+2}}
\, e_{n+1} -(1-q^{2\ta})q^{2n}\, e_n -iq^{\ta+n-1}e^{-i\th}
\sqrt{1-q^{2n}}\, e_{n-1},
$$
to see that $v_\la^\th$ as defined in Proposition 
\thmref{\propeigvectpirti} is indeed formally an eigenvector for 
$\pi^\infty_\th(\rti)$ for the eigenvalue $\la$.
For $\la=-q^{2k}$ or $\la=q^{2\ta+2k}$, $k\in\Zp$, it is easy to 
show that $v_\la^\th\in\Hi$. The orthogonality follows, since 
$\pi^\infty_\th(\rti)$ is self-adjoint.

It remains to calculate the squared norm and to prove the 
completeness in $\Hi$. The squared norm can be calculated 
using the fact that we already have established orthogonality 
and some easy calculations using Theorem \thmref{\thmqbinomialthm} 
and Corollary \thmref{\corstothmqbinomialthm}(iii). The 
completeness follows, since the dual orthogonality relations also 
hold, which are in fact the orthogonality relations for the 
Al-Salam and Carlitz polynomials.
\qed\enddemo

\demo{Proof of Theorem \thmref{\thmHaaronspecialspherelt}}
We calculate the trace with respect to
the orthogonal basis of eigenvectors described in
Proposition \thmref{\propeigvectpirti};
$$
\text{tr}\Bigl( D\pi^\infty_\th\bigl(p(\rti)\bigr)\Bigr) =
\sum_{k=0}^\infty p(-q^{2k}) {{\langle D v_{-q^{2k}}^\th,
v_{-q^{2k}}^\th\rangle}\over{\langle v_{-q^{2k}}^\th,
v_{-q^{2k}}^\th\rangle}}  +
\sum_{k=0}^\infty p(q^{2\ta+2k}){{\langle D v_{q^{2\ta+2k}}^\th,
v_{q^{2\ta+2k}}^\th\rangle}\over{\langle v_{q^{2\ta+2k}}^\th,
v_{q^{2\ta+2k}}^\th\rangle}}.
$$
So it remains to calculate the matrix coefficients on the diagonal
of the operator $D$ with respect to this basis. This can be done
by a straightforward calculation, or by use of the so-called
$q$-Charlier polynomials, cf. Exercise~4.3. The result is
$$
\align
\langle Dv_{-q^{2k}}^\th,v_{-q^{2k}}^\th\rangle &=
(-q^{2\ta+2};q^2)_\infty (q^2;q^2)_k (-q^{2-2\ta};q^2)_k,\\
\langle Dv_{q^{2\ta+2k}}^\th,v_{q^{2\ta+2k}}^\th\rangle &=
(-q^{2-2\ta};q^2)_\infty (q^2;q^2)_k (-q^{2+2\ta};q^2)_k,
\endalign
$$

Thus
$$
\text{tr}\Bigl( D\pi^\infty_\th\bigl(p(\rti)\bigr)\Bigr)
= {1\over{1+q^{2\ta}}}
\biggl( \sum_{k=0}^\infty p(-q^{2k}) q^{2k} +
q^{2\ta} \sum_{k=0}^\infty p(q^{2\ta+2k}) q^{2k}\biggr).
$$
This expression is independent of $\th$, so that we obtain
from Corollary \thmref{\corlemHaaronbasis}
$$
h\bigl(p(\rti)\bigr) = {{1-q^2}\over{1+q^{2\ta}}} \biggl(
\sum_{k=0}^\infty p(-q^{2k}) q^{2k} +
q^{2\ta} \sum_{k=0}^\infty p(q^{2\ta+2k}) q^{2k}\biggr),
$$
which is precisely the definition of the $q$-integral.
\qed\enddemo

The proof of Theorem \thmref{\thmHaaronspecialspherelt}
is straightforward, since the representation space
$\Hi$ has an orthonormal basis of eigenvectors
of $\pi^\infty_\th(\rti)$ for each $\th$, so that the
spectrum of $\pi^\infty_\th(\rti)$ is purely discrete.

To prove the general case we have to invoke Jacobi matrices,
as in \S 3.1. We can show that $\pi^\infty_\th(\rts)$
on each of the components of the decomposition
$\Hi = V_1^\th\oplus V_2^\th$ of
Remark \thmref{\rempropeigvectpirti} is represented
by a Jacobi matrix, although $\pi^\infty_\th(\rts)$ in the
standard basis is represented by a five-term recurrence operator.
On each of the components we can determine the corresponding
orthogonality measure using Al-Salam and Chihara polynomials,
see Exercise~3.4.
Then the Poisson kernel for the Al-Salam and Chihara polynomials
comes into play, which is a (very-well poised)
${}_8\vp_7$-series, on each of these
components. Finally, the results have to be matched using
Bailey's formula for the sum of two of such ${}_8\vp_7$-series
of the correct form. In Exercises~4.4--7 more explicit hints
are given. In order to show how this works we consider
the case of the Haar functional on the cocentral elements,
in which all ingredients are contained but in which the
computations are much simpler.

\proclaim{Theorem \theoremname{\thmHaarcocentral}} The Haar 
functional on the subalgebra generated by the self-adjoint element
$\al+\al^\ast$ is given by the integral
$$
h\bigl(p((\al+\al^\ast)/2)\bigr) = {2\over \pi}\int_{-1}^1 p(x)
\sqrt{1-x^2}\, dx
$$
for any polynomial $p$.
\endproclaim

\demo{Proof} Consider
$$
2\pi^\infty_\th\big( (\al+\al^\ast)/2\bigr) e_n = \sqrt{1-q^{2n}}
e_{n-1} + \sqrt{1-q^{2n+2}} e_{n+1}.
\tag\eqname{\vglcocentralpi}
$$
So the operator $\pi^\infty_\th\big( (\al+\al^\ast)/2\bigr)$ is
represented by a Jacobi matrix with respect to the standard basis 
of $\Hi$.

Take $a=b=0$ in the Al-Salam and Chihara polynomials of 
Exercise~3.4 to obtain Rogers's continuous $q$-Hermite polynomials 
$H_n(x|q)$ satisfying
$$
2 x\,  H_n(x| q) = H_{n+1}(x|q) + (1-q^n)\, H_{n-1}(x| q), \quad
H_{-1}(x| q)=0,\ H_0(x| q)=1.
\tag\eqname{\vglttrqHermite}
$$
The continuous $q$-Hermite polynomials satisfy the orthogonality 
relations
$$
\int_0^{\pi} H_n(\cos\phi| q)H_m(\cos\phi| q)
w(\cos\phi| q) \, d\phi = \de_{nm} {{2\pi
(q;q)_n}\over{(q;q)_\infty}},
%\tag\eqname{\vglorthoqHermite}
$$
with $w(\cos\phi| q) = (e^{2i\phi},e^{-2i\phi};q)_\infty$, by
taking $a=b=c=d=0$ in the Askey-Wilson orthogonality measure
of Theorem \thmref{\thmAWpolsareorthogonal}.

Compare \thetag{\vglcocentralpi} with the three-term recurrence 
relation \thetag{\vglttrqHermite} to see that the orthonormal
polynomials associated to the Jacobi matrix
$\pi^\infty_\th\big( (\al+\al^\ast)/2\bigr)$ are
the orthonormal continuous $q$-Hermite polynomials $H_n(x|
q^2)/\sqrt{(q^2;q^2)_n}$.

We use the spectral theory of Jacobi matrices
as in the proof of Theorem \thmref{\thmFavard}
with $p_n(x) = H_n(x|q^2)/\sqrt{(q^2;q^2)_n}$ as the 
orthonormal polynomials and the absolutely continuous measure
$dm(x| q^2)=(2\pi)^{-1}(q^2;q^2)_\infty w(x| q^2)(1-x^2)^{-1/2}dx$ 
on $[-1,1]$ as the (normalised) orthogonality measure.
So we obtain a unitary mapping $\Lambda$,
intertwining $\pi^\infty_\th\bigl((\al+\al^\ast)/2\bigr)$ on $\Hi$
with the multiplication operator $M$ on $L^2([-1,1],dm(x|q^2))$.
Hence,
$$
\align
\text{tr}\Bigl( D\pi^\infty_\th\bigl(&p((\al+\al^\ast)/2)\bigr)
\Bigr) =
\sum_{n=0}^\infty q^{2n}
\langle \pi^\infty_\th\bigl(p((\al+\al^\ast)/2)\bigr)\, e_n, e_n
\rangle \\
&= \sum_{n=0}^\infty q^{2n} \langle
\Lambda \pi^\infty_\th\bigl(p((\al+\al^\ast)/2)\bigr)\, e_n,
\Lambda e_n \rangle \\
&= \sum_{n=0}^\infty q^{2n} \int_{-1}^1
p(x)\, \bigl\vert\bigl(\Lambda e_n\bigr)(x)\bigr\vert^2\, 
dm(x| q^2) \\
&= \sum_{n=0}^\infty q^{2n} \int_{-1}^1
p(x)\, {{\bigl(H_n(x|q^2)\bigr)^2}\over{(q^2;q^2)_n}}\, 
dm(x| q^2) \\
&= \int_{-1}^1 p(x) P_{q^2}(x,x| q^2) \, dm(x| q^2),
\endalign
%\tag\eqname{\vglcocentralone}
$$
where $P_t(x;y|q)$ is the Poisson kernel for the continuous
$q$-Hermite polynomials. Explicitly, $P$ is given by
$$
\multline
P_t(\cos\phi,\cos\psi| q) =
\sum_{n=0}^\infty {{H_n(\cos\phi| q)H_n(\cos\psi| q)
t^n}\over{(q;q)_n}} = \\{{(t^2;q)_\infty}\over{(te^{i\phi+i\psi},
te^{i\phi-i\psi},te^{-i\phi+i\psi},te^{-i\phi-i\psi};q)_\infty}},
\qquad |t|<1.
\endmultline
%\tag\eqname{\vglPoissonqHermite}
$$
Interchanging summation
and integration is easily justified.

{}From the explicit expression for the Poisson kernel we get
$$
w(x| q^2) P_{q^2}(x,x| q^2) =
{{4(1-x^2)}\over{(1-q^2) (q^2;q^2)_\infty}},
$$
so that
$$
\text{tr}\Bigl( D\pi^\infty_\th\bigl(p((\al+\al^\ast)/2)\bigr)
\Bigr) = {2\over{\pi(1-q^2)}} \int_{-1}^1 p(x)\sqrt{1-x^2}\, dx.
$$
Since this is independent of the parameter $\th$ of the infinite
dimensional representation, the proof follows from
Corollary \thmref{\corlemHaaronbasis}.
\qed\enddemo

%%%%%%%%%%%%%%%%%%%%%%%%%%%%%%%%%%%%%%%%%%%%%%%%%%%%%%%%%%%%%%%%%%%%
%%N E W   S U B S E C T I O N%%%%%%%%%%%%%%%%%%%%%%%%%%%%%%%%%%%%%%%
%%%%%%%%%%%%%%%%%%%%%%%%%%%%%%%%%%%%%%%%%%%%%%%%%%%%%%%%%%%%%%%%%%%%
%NOTES AND REFERENCES%%%%%%%%%%%%%%%%%%%%%%%%%%%%%%%%%%%%%%%%%%%%%%%
\subhead Notes and references
\endsubhead
The introduction of these infinitesimally generated `subgroups'
of the quantum $SU(2)$ group is due to Koornwinder \cite{\KoorZSE}.
This paper has been very influential for the development
of the relation between $q$-special functions and quantum groups.
The rest of \S 4.1 has been taken from an unpublished
announcement \cite{\NoumMunpub} by Noumi and Mimachi,
of which an even shorter announcement
\cite{\NoumMPJA} has appeared, see also \cite{\NoumMLNM},
and from \cite{\KoelSIAM}, \cite{\KoelAAM}, \cite{\KoelFIC}.
The (dual) $q$-Krawtchouk polynomials and the element
of $X_\si A$ also naturally occur when determining the
spherical functions on the Hecke algebra of type $B_n$
with respect to the parabolic subalgebra corresponding to $A_{n-1}$,
see \cite{\KoelHAlg}.

Lemma \thmref{\lemHaaronbasis} is due to Woronowicz \cite{\Woro},
but its proof is taken from Noumi and Mimachi \cite{\NoumMCompM},
and Corollary \thmref{\corlemHaaronbasis} is due to Vaksman
and Soibelman \cite{\VaksS}.
See Dunford and Schwartz \cite{\DunfS, Ch.~XI, \S 6} for more 
details on Hilbert-Schmidt operators and trace class operators.
In fact, $\parallel a\parallel$, $a\in\Asu$, is a
$C^\ast$-norm, and we can complete $\Asu$ into a $C^\ast$-algebra
in which $\Asu$ is a dense subalgebra. This corresponds to
Woronowicz's approach to the quantum $SU(2)$ group, \cite{\Woro},
\cite{\Worotwee}. The key Theorem \thmref{\thmHaarongeneralspherelt}
is due to Koornwinder \cite{\KoorZSE}, but the proof sketched here 
is taken from Koelink and Verding \cite{\KoelV}.
The details of the last part of the proof of Proposition
\thmref{\propeigvectpirti} can be found in \cite{\KoelCJM}.
Theorem \thmref{\thmHaarcocentral} is due to Woronowicz 
\cite{\Woro}, with a different proof. The Poisson kernel
is known from the work by Rogers (1894-6) on the continuous
$q$-Hermite polynomials, see \cite{\AskeI}, \cite{\Bres}.

%%%%%%%%%%%%%%%%%%%%%%%%%%%%%%%%%%%%%%%%%%%%%%%%%%%%%%%%%%%%%%%%%%%%
%%N E W   S U B S E C T I O N%%%%%%%%%%%%%%%%%%%%%%%%%%%%%%%%%%%%%%%
%%%%%%%%%%%%%%%%%%%%%%%%%%%%%%%%%%%%%%%%%%%%%%%%%%%%%%%%%%%%%%%%%%%%
%EXERCISES%%%%%%%%%%%%%%%%%%%%%%%%%%%%%%%%%%%%%%%%%%%%%%%%%%%%%%%%%%
\subhead Exercises
\endsubhead

\item{\the\sectionno.1} Fill in the proof of Proposition
\thmref{\propshifteigenval}.

\item{\the\sectionno.2} Fill in the gaps of the last paragraph
of the proof of Proposition \thmref{\propeigvectpirti}. Use the
orthogonality relations for the Al-Salam and Carlitz polynomials;
$$
\multline
\int_a^1 U_m^{(a)}(q;x)U_n^{(a)}(q;x) (qx,qx/a;q)_\infty\, d_qx =\\
\de_{nm} (-a)^n (1-q) q^{n(n-1)/2} (q;q)_n (q,a,1/a;q)_\infty
\endmultline
$$
assuming $a<0$ for the measure to be non-negative.

\item{\the\sectionno.3} Define the discrete orthogonality measure
$\mu$ by
$$
\int_\R p(x)\, d\mu(x) =
\sum_{n=0}^\infty {{q^{2n\ta}q^{n(n-1)}}\over{(q^2;q^2)_n}}
p(q^{-2n}),
$$
and define the $q$-Charlier polynomials $c_n(x;a;q) =
{}_2\vp_1(q^{-n},x;0;q,-q^{n+1}/a)$. Prove from Proposition
\thmref{\propeigvectpirti}
$$
\int_\R c_k(x;q^{2\ta};q^2) c_l(x;q^{2\ta};q^2) \, d\mu(x)=
\de_{k,l} q^{-2k} (q^2,-q^{2-2\ta};q^2)_k (-q^{2\ta};q^2)_\infty.
$$
Assuming $k\geq l$, prove that
$$
\align
\langle Dv_{-q^{2k}}^\th,v_{-q^{2l}}^\th\rangle & =
\int_\R {1\over x} c_k(x;q^{2\ta};q^2)c_l(x;q^{2\ta};q^2) 
\, d\mu(x) \\ & =
c_l(0;q^{2\ta};q^2) \int_\R {1\over x} c_k(x;q^{2\ta};q^2) 
\, d\mu(x),
\endalign
$$
and calculate the last integral explicitly using the summation
formulas of \S 3.

\item{\the\sectionno.4} Use the notation 
$v_\la^\th(q^\ta)=v_\la^\th$
for the orthogonal basis of Proposition \thmref{\propeigvectpirti}
in order to stress the $\ta$-dependence. Prove that
$$
\align
\pi^\infty_\th(\al_{\ta,\infty})v^\th_\la(q^\ta) &=
e^{i\th}iq^{1/2-\ta}(1+\la) v^\th_{\la/q^2}(q^{\ta-1}), \\
\pi^\infty_\th(\be_{\ta,\infty})v^\th_\la(q^\ta)& =
e^{-i\th}iq^{1/2} v^\th_\la (q^{\ta-1}),\\
\pi^\infty_\th(\ga_{\ta,\infty})v^\th_\la(q^\ta) &=
e^{i\th}iq^{1/2} (q^{2\ta}-\la)v^\th_{\la}(q^{\ta+1}),\\
\pi^\infty_\th(\de_{\ta,\infty})v^\th_\la(q^\ta) &=
-e^{-i\th}iq^{1/2+\ta} v^\th_{\la q^2}(q^{\ta+1})
\endalign
$$
and from this that
$$
2\pi^\infty_\th(\rts)v_\la^\th = qe^{-2i\th}v_{\la q^2}^\th + 
q^{-1}e^{2i\th}(1 - q^{-2\ta}\la)(1+\la)v_{\la/q^2}^\th
+ \la q^{1-\ta}(q^{-\si}-q^\si)v_\la^\th.
$$
Use Proposition \thmref{\propsphericaleltsalgebra} and that
$\al_{\ta,\si}$ can be written as a linear combination of
$\al_{\ta,\infty}$ and $\be_{\ta,\infty}$, and similarly for the
other elements.

\item{\the\sectionno.5} $\pi^\infty_\th(\rts)$ preserves the 
orthogonal decomposition of Remark \thmref{\rempropeigvectpirti}, 
and is given by a Jacobi matrix on each component. Let
$$
D=\pmatrix D^\th_{1,1}&D^\th_{1,2}\\ D^\th_{2,1}&D^\th_{2,2}
\endpmatrix
$$
be the corresponding decomposition of $D$. So that
$$
\text{tr}\bigl( D\pi^\infty_\th(p(\rts))\bigr) =
\text{tr}_{V_1^\th}\bigl( D_{11}^\th\pi^\infty_\th(p(\rts))\bigr) +
\text{tr}_{V_2^\th}\bigl( D_{22}^\th\pi^\infty_\th(p(\rts))\bigr).
$$
Denote by $w_m^\th$, $m\in\Zp$, the orthonormal basis of $V_1^\th$
obtained by normalising $v_{-q^{2m}}^\th$, $m\in\Zp$, and by
$u_m^\th$, $m\in\Zp$, the orthonormal basis of $V_2^\th$
obtained by normalising $v_{q^{2\ta + 2m}}^\th$, $m\in\Zp$.
Show that the corresponding orthonormal polynomials can be given 
in terms of orthonormal Al-Salam and Chihara polynomials,
$$
h_n(x;s,t|q) = (q,-q/s^2;q)^{-1/2}\, s_n(x;q^{1/2}t/s, 
-q^{1/2}/st|q),
$$
with the notation of Exercise~3.4. And prove that
$$
\multline
\text{tr}_{V_1^\th}\bigl( D_{11}^\th\pi^\infty_\th(p(\rts))\bigr) =
\int_{\R}p(x) \sum_{n=0}^{\infty} \sum_{m=0}^{\infty}
\langle  D w_n^\th,w_m^\th\rangle \\ \times h_n(x;q^\ta,q^\si| q^2)
h_m(x;q^\ta,q^\si| q^2)e^{2i(m-n)\th}\, dm(x;q^{1+\si-\ta},
-q^{1-\si-\ta},0,0| q^2)
\endmultline
$$
and similarly
$$
\multline
\text{tr}_{V_2^\th}\bigl( D_{22}^\th\pi^\infty_\th(p(\rts))\bigr) =
\int_{\R}p(x) \sum_{n=0}^{\infty} \sum_{m=0}^{\infty}
\langle  D u_n^\th,u_m^\th\rangle \\ \times h_n(x;q^{-\ta},
q^{-\si}| q^2)
h_m(x;q^{-\ta},q^{-\si}| q^2)e^{2i(m-n)\th}\,
dm(x;q^{1-\si+\ta},-q^{1+\si+\ta},0,0| q^2)
\endmultline
$$
using the notation for the normalised orthogonality measure for
Askey-Wilson polynomials, cf. \thetag{\vglnormalisedAWmeasure}.

\item{\the\sectionno.6} Prove that $h\bigl(p(\rts)\bigr) = $
$$
\multline
{{1-q^2}\over {1+q^{2\ta}}}
\int_\R p(x) P_{q^2}(x,x;q^{1+\si-\ta},-q^{1-\si-\ta}\vert q^2)
\, dm(x;q^{1+\si-\ta},-q^{1-\si-\ta},0,0| q^2)
\\ + {{1-q^2}\over {1+q^{-2\ta}}}
\int_\R p(x) P_{q^2}(x,x;q^{1-\si+\ta},-q^{1+\si+\ta}\vert q^2)
\,dm(x;q^{1-\si+\ta},-q^{1+\si+\ta},0,0| q^2) ,
\endmultline
$$
where
$$
P_t(x,y;a,b\vert q) =
\sum_{k=0}^\infty t^k {{s_k(x;a,b\vert q)s_k(y;a,b
\vert q)} \over {(q,ab;q)_k}}
$$
is the Poisson kernel for the Al-Salam--Chihara polynomials.

\item{\the\sectionno.7} Using the standard notation for very-well
poised ${}_8\vp_7$-series, cf. \cite{\GaspR, Ch.~2},
$$
{}_8W_7(a;b,c,d,e,f;q,z) = {}_8\vp_7 \left(
{{a,\, q\sqrt{a},\, -q\sqrt{a},\, b,\, c,\, d,\, e,\, f}\atop
{\sqrt{a},-\sqrt{a},qa/b,qa/c,qa/d,qa/e,qa/f}};q,z\right),
$$
the Poisson kernel for the Al-Salam and Chihara polynomials has 
been evaluated by Askey, Rahman and Suslov, see 
\cite{\AskeRS, (14.8)},
$$
\multline
P_t(\cos\theta,\cos\psi;a,b\vert q) =
{{(ate^{i\theta},ate^{-i\theta},bte^{i\psi},
bte^{-i\psi},t;q)_\infty }\over
{(te^{i\theta +i\psi},te^{i\theta -i\psi},te^{i\psi -i\theta},
te^{-i\psi -i\theta} ,abt;q)_\infty}} \\
\times {}_8W_7({{abt}\over q};t,be^{i\theta},be^{-i\theta},
ae^{i\psi} ,ae^{-i\psi};q,t).
\endmultline
$$
Use Bailey's formula, cf. \cite{\GaspR, (2.11.7)},
$$
\multline
{1\over {(b/a;q)_\infty}} {}_8W_7(a;b,c,d,e,f ;q,q) +
{{(a q, c ,d ,e ;q)_\infty}\over
{(a/b ;q)_\infty (aq/c, aq/d, aq/e, aq/f;q)_\infty}} \\
\times {{(f ,bq/ c , bq/d ,bq/e, bq/f;q)_\infty}\over
{(bc/a, bd/a, be/a, bf/a,
b^2q/a;q)_\infty}} \,  {}_8W_7({b^{ 2}\over a };b ,
{{b c }\over a },{{b d }\over a }, {{b e }\over a },{{b f }
\over a } ;q,q) \\
= {{(aq, aq/(cd), aq/(ce), aq/(cf), aq/(de), aq/(df), 
aq/(ef);q)_\infty}\over
{(aq/c, aq/d, aq/e, aq/f, bc/a, bd/a, be/a, bf/a; q)_\infty}}
\endmultline
$$
to finish the proof of Theorem 
\thmref\thmref{\thmHaarongeneralspherelt}.
The absolutely continuous part of the measure has to be treated
differently from possible discrete mass points.

%%%%%%%%%%%%%%%%%%%%%%%%%%%%%%%%%%%%%%%%%%%%%%%%%%%%%%%%%%%%%%%%%%%%
%%N E W   S E C T I O N%%%%%%%%%%%%%%%%%%%%%%%%%%%%%%%%%%%%%%%%%%%%%
%%%%%%%%%%%%%%%%%%%%%%%%%%%%%%%%%%%%%%%%%%%%%%%%%%%%%%%%%%%%%%%%%%%%
\newpage

\head\newsection Askey-Wilson polynomials and
generalised matrix elements\endhead

We have now developed all the necessary ingredients for the
interpretation of Askey-Wilson polynomials on the quantum $SU(2)$ 
group. In this section we first show that orthogonal polynomials 
are of importance in describing generalised matrix elements. 
Next these polynomials are explicitly calculated in terms of 
Askey-Wilson polynomials.

%%%%%%%%%%%%%%%%%%%%%%%%%%%%%%%%%%%%%%%%%%%%%%%%%%%%%%%%%%%%%%%%%%%%
%%N E W   S U B S E C T I O N%%%%%%%%%%%%%%%%%%%%%%%%%%%%%%%%%%%%%%%
%%%%%%%%%%%%%%%%%%%%%%%%%%%%%%%%%%%%%%%%%%%%%%%%%%%%%%%%%%%%%%%%%%%%
\subhead\newsubsection 
Generalised matrix elements and orthogonal polynomials
\endsubhead
First we establish a relation between generalised matrix elements
and orthogonal polynomials.

\proclaim{Theorem \theoremname{\thmabstrorthpolgenmatelt}}
For fixed $i,j\in\hZp$ such that $i-j\in\Z$, there exists a system 
of orthogonal polynomials $(p_k)_{k\in\Zp}$ of degree $k$ such 
that for $l\geq m =\max (\vert i\vert,\vert j\vert)$, $l-m\in\Zp$,
$$
b^l_{i,j}(\ta,\si) = b^m_{i,j}(\ta,\si) p_{l-m}(\rts ).
$$
\endproclaim

\demo{Proof} We first prove that an expression of this form
exists. Consider for any polynomial
$s_{l-m}$ of degree $l-m$ the expression
$b^m_{i,j}(\ta,\si) s_{l-m}(\rts )$.
If we decompose this product with respect to the linear basis
the decomposition of $\Asu$ of Theorem
\thmref{\thmdescripAasdualHAofU}  we get
$$
b^m_{i,j}(\ta,\si) s_{l-m}(\rts ) = \sum_{k=\vert 2m-l\vert}^l b^k,
\qquad b^k \in \text{span}(t^k_{nm}),
$$
where the upper and lower bound follow from Lemma 
\thmref{\lemmaCGCforUsu}. The mappings $X.$ and $.X$ preserve
$\text{span}(t^k_{nm})$. Proposition \thmref{\propeigdefmatelt}(i)
shows that the left
hand side satisfies
\thetag{\vglgenbiinvariancet} with $\la=\la_j(\si)$ and
$\mu=\la_i(\ta)$. Consequently, each $b^k$ has to satisfy
\thetag{\vglgenbiinvariancet} with $\la=\la_j(\si)$ and
$\mu=\la_i(\ta)$. Proposition \thmref{\propeigdefmatelt}(ii) 
implies that $b^k=0$ for $k<m$ and $b^k=c_kb^k_{i,j}(\ta,\si)$ 
for $k\geq m$ and some constants $c_k$. Hence,
$$
b^m_{i,j}(\ta,\si) s_{l-m}(\rts ) = \sum_{k=m}^l c^k b^k_{i,j}
(\ta,\si).
$$
Since both sides contain the same degree of freedom,
the existence of such polynomials follows once we know that the 
mapping $s_{l-m}\mapsto b^m_{i,j}(\ta,\si)s_{l-m}(\rts )$ is 
injective. This can be seen by applying the one-dimensional 
$\ast$-representation $\pi_\th$ of $\Asu$, cf.
Theorem \thmref{\thmirredstarrepsofAsu},
and use of the explicit expression of $b^m_{i,j}(\ta,\si)$, cf.
\thetag{\vglexprdefgenmateltinstanmatelt}. Then use 
$\pi_\th(\rts)=\cos2\th$ and 
$\pi_\th(t^l_{nm})=\de_{nm}e^{-2in\th}$, which follows by 
observing that 
$\pi_\th(\xi)=\langle A^{2i\th/\ln q},\xi\rangle$ holds (formally)
for $\xi\in\Asu$.

For $l,k\geq m$, $l-m,k-m\in\Zp$ we have
$b^l_{i,j}(\ta,\si)\in\text{span}(t^l_{nm})$,
$b^k_{i,j}(\ta,\si)\in\text{span}(t^k_{nm})$, so that the Schur 
orthogonality relations for the Haar functional $h$, cf.
Theorem \thmref{\propSchurorthoforSU}, imply
$$
h\Bigl( \bigl(b^l_{i,j}(\ta,\si)\bigr)^\ast b^k_{i,j}(\ta,\si)
\Bigr) = \de_{k,l}h_l, \qquad h_l>0 .
\tag\eqname{\vgldefinitionhsubl}
$$
Now $\bigl(b^m_{i,j}(\ta,\si)\bigr)^\ast b^m_{i,j}(\ta,\si)=
w_m(\rts)$ for some
polynomial $w_m$ of degree $2m$,
by Proposition \thmref{\propeigdefmatelt}(i) and
\thmref{\propsphericaleltsalgebra}, and the
Clebsch-Gordan series of Lemma \thmref{\lemmaCGCforUsu}.
Hence,
$$
h\bigm( \bar p_{l-m}(\rts)p_{k-m}(\rts) w_m(\rts)\bigr) = 
\delta_{l,k} h_l, \qquad h_l >0,
$$
since we have already established the existence of such
polynomials and since $\rts^\ast=\rts$. Consequently,
by Theorem \thmref{\thmHaarongeneralspherelt} the polynomials
$p_{l-m}$, $l-m\in\Zp$, form a system of orthogonal polynomials with
respect to the measure
$w_m(x) \, dm(x;a,b,c,d\mid q^2)$ on $\R$,
where $dm(x;a,b,c,d\mid q^2)$ is the measure described in
Theorem \thmref{\thmHaarongeneralspherelt}.
\qed\enddemo

%%%%%%%%%%%%%%%%%%%%%%%%%%%%%%%%%%%%%%%%%%%%%%%%%%%%%%%%%%%%%%%%%%%%
%%N E W   S U B S E C T I O N%%%%%%%%%%%%%%%%%%%%%%%%%%%%%%%%%%%%%%%
%%%%%%%%%%%%%%%%%%%%%%%%%%%%%%%%%%%%%%%%%%%%%%%%%%%%%%%%%%%%%%%%%%%%
\subhead\newsubsection 
Generalised matrix elements and Askey-Wilson polynomials
\endsubhead
Combining Theorem \thmref{\thmabstrorthpolgenmatelt},
and in particular the description of the orthogonality measure,
with Theorem \thmref{\thmHaarongeneralspherelt} shows that it 
suffices to calculate
$$
w_m(\cos\theta)= \bigl\vert \pi_{\th/2}
\bigl( b^m_{i,j}(\ta,\si)\bigr)\bigr\vert^2, \qquad
m=\max(|i|,|j|),
$$
to find the orthogonality measure for the orthogonal polynomials
from Theorem \thmref{\thmabstrorthpolgenmatelt}. The first step is
the following proposition, where we use the one-dimensional
representations $\ta_\la$ defined in
Remark \thmref{\remdefinnonunitaryonedrepsofAsu}.

\proclaim{Proposition 
\theoremname{\propexplvalueonedimonmingenmatelt}}
For $i,j\in\hZp$, $i-j\in\Z$ and $m=\max 
(\vert i\vert,\vert j\vert)$ we have

\noindent
{\rm (i)} In case $m=i$ or $-i\leq j\leq i$
$$
\ta_\la \bigl( b^i_{i,j}(\ta,\si)\bigr) = C^{i,j}(\si)C^{i,i}(\ta)
q^{-i} \la^{-2i} (\la^2q^{1+\ta-\si};q^2)_{i-j}
(-\la^2q^{1+\ta+\si};q^2)_{i+j}.
$$

\noindent
{\rm (ii)} In case $m=j$ or $-j\leq i\leq j$
$$
\ta_\la \bigl( b^j_{i,j}(\ta,\si)\bigr) = C^{j,i}(\ta)C^{j,j}(\si)
q^{-j} \la^{-2j} (\la^2q^{1+\si-\ta};q^2)_{j-i}
(-\la^2q^{1+\ta+\si};q^2)_{i+j}.
$$

\noindent
{\rm (iii)} In case $m=-i$ or $i\leq j\leq -i$
$$
\ta_\la \bigl( b^{-i}_{i,j}(\ta,\si)\bigr) = C^{-i,-j}(-\si)
C^{-i,-i}(-\ta) q^i \la^{2i} (\la^2q^{1-\ta+\si};q^2)_{j-i}
(-\la^2q^{1-\ta-\si};q^2)_{-i-j}.
$$

\noindent
{\rm (iv)} In case $m=-j$ or $j\leq i\leq -j$
$$
\ta_\la \bigl( b^{-j}_{i,j}(\ta,\si)\bigr) = C^{-j,-i}(-\ta)
C^{-j,-j}(-\si) q^j \la^{2j} (\la^2q^{1-\si+\ta};q^2)_{i-j}
(-\la^2q^{1-\ta-\si};q^2)_{-i-j}.
$$
\endproclaim

\demo{Proof} First we observe that
the function $\ta_\la \bigl( b^l_{i,j}(\ta,\si)\bigr)$ of $\la$ 
satisfies the symmetry relations
$$
\aligned
\ta_\la \bigl( b^l_{i,j}(\ta,\si)\bigr) &=
\overline{ \ta_{\bar \la} \bigl( b^l_{j,i}(\si,\ta)\bigr)} \\
&= \ta_\la \bigl( b^l_{-j,-i}(-\si,-\ta)\bigr)\\
&= \overline{ \ta_{\bar \la} \bigl( b^l_{-i,-j}(-\ta,-\si)\bigr)} .
\endaligned
\tag\eqname{\vglsymmetriesgenmatelt}
$$
This follows from $\ta_\la(t^l_{nm})=\de_{nm}\la^{-2n}$ and
$v_n^{l,j}(\si)=\overline{v^{l,-j}_n(-\si)}$, which in turn
follows from $C^{l,j}(\si)=C^{l,-j}(-\si)$ and
Exercise~3.8.

So it suffices to prove the first statement.
Use Proposition \thmref{\propshifteigenval} to find
$$
\multline
b^l_{l,m}(\ta,\si)= C^{l,m}(\si) C^{l,l}(\ta) q^{\si (m-l)}
q^{{1\over 2}(l-m)(l-m-1)} \\
\times
\prod_{k=0}^{l+m-1} \de_{\ta+2l-1-k,\si+2m-1-k}
\prod_{j=0}^{l-m-1} \ga_{\ta+l-m-1-j,\si-l+m+1+j}.
\endmultline
$$
Apply $\ta_\la$ to find
$$
\multline
\ta_\la \bigl( b^l_{l,m}(\ta,\si)\bigr)  =
C^{l,m}(\si) C^{l,l}(\ta) q^{\si (m-l)}
q^{{1\over 2}(l-m)(l-m-1)} \\
\times
\prod_{k=0}^{l+m-1} (q^{\ta+\si+2l+2m-2-2k+1/2}\la + 
q^{-1/2}\la^{-1}) \\ \times
\prod_{j=0}^{l-m-1} (-q^{\ta+l-m-1-j+1/2}\la + 
q^{\si-l+m+1+j-1/2}\la^{-1}) \\
= C^{l,m}(\si)C^{l,l}(\ta) q^{-l} \la^{-2l} 
(\la^2q^{1+\ta-\si};q^2)_{l-m} (-\la^2q^{1+\ta+\si};q^2)_{l+m}
\endmultline
$$
and this proves (i).
\qed\enddemo

\demo{Remark \theoremname{\rempropexplvalueonedimonmingenmatelt}}
We can lift the symmetry relations of
\thetag{\vglsymmetriesgenmatelt} to the algebra $\Asu$,
as follows. Let $\Psi\colon\Asu\to\Asu$ be the algebra isomorphism
obtained by interchanging $\be$ and $\ga$. It follows
from Corollary \thmref{\corthmexplicitdualityonbasisone} that 
this is well-defined. From
\thetag{\vgldefDeltaonAq} we see that $\Psi$ is an anti-coalgebra
isomorphism. Hence, $\Psi$ gives a coalgebra-isomorphism and
anti-isomorphism of $\U$, which is just interchanging $B$ and $C$.
Using this, Theorem \thmref{\thmunitaryrtepsofUsu} and
Lemma \thmref{\lemPBWforU}, then implies $\Psi(t^l_{nm})=t^l_{mn}$.

Similarly,
let $\Phi\colon\Asu\to\Asu$ be the anti-algebra isomorphism
obtained by interchanging $\al$ and $\de$. Again
Corollary \thmref{\corthmexplicitdualityonbasisone} implies that
it is well-defined. From
\thetag{\vgldefDeltaonAq} we see that $\Phi$ is an anti-coalgebra
isomorphism. Hence, $\Phi$ gives a anti-coalgebra isomorphism and
anti-isomorphism of $\U$, which is just interchanging $A$ and $D$.
Using this, Theorem \thmref{\thmunitaryrtepsofUsu} and
Lemma \thmref{\lemPBWforU}, then implies $\Phi(t^l_{nm})=
t^l_{-m,-n}$. Combining gives 
$\Phi\circ\Psi(t^l_{nm})=t^l_{-n,-m}$, where
$\Phi\circ\Psi$ is interchanging $\be$ and $\ga$ and $\al$ and 
$\de$.
\enddemo

We can now prove the main result of this section in which we relate
Askey-Wilson polynomials to generalised matrix elements. The 
Askey-Wilson polynomials involve four continuous parameters. 
In the following theorem we establish an interpretation of the 
Askey-Wilson polynomials with two continuous and two discrete 
parameters. We rewrite the Askey-Wilson polynomials using the 
following notation;
$$
p_n^{(\al,\be)}(x;s,t|q) = p_n(x;q^{1/2}t/s,q^{1/2+\al}s/t,
-q^{1/2}/(st), -stq^{1/2+\be}| q).
$$
In this way the Askey-Wilson polynomials are $q$-analogues of the
Jacobi polynomials. The case $s=t=1$ goes under the name of
continuous $q$-Jacobi polynomials. Note also that the Al-Salam
and Chihara polynomials with the parametrisation as in Exercise~4.5
can be considered as the corresponding Hermite polynomial of
$p_n^{(\al,\be)}(x;s,t|q)$, i.e. they can be obtained by letting
$\al=\be\to\infty$.

\proclaim{Theorem \theoremname{\thmdefgenmateltAskeyWpol}}
For $i,j\in\hZp$, $i-j\in\Z$ and $l-m\in\Zp$, 
$m=\max (\vert i\vert,\vert j\vert)$ we have

\noindent
{\rm (i)} In case $m=i$ or $-i\leq j\leq i$:
$$
b^l_{i,j}(\ta,\si) = d^l_{i,j}(\ta,\si) b^i_{i,j}(\ta,\si)
p^{(i-j,i+j)}_{l-i}(\rts;q^\ta, q^\si \mid q^2).
$$

\noindent
{\rm (ii)} In case $m=j$ or $-j\leq i\leq j$:
$$
b^l_{i,j}(\ta,\si) = d^l_{j,i}(\si,\ta) b^j_{i,j}(\ta,\si)
p^{(j-i,i+j)}_{l-j}(\rts;q^\si, q^\ta \mid q^2).
$$

\noindent
{\rm (iii)} In case $m=-i$ or $i\leq j\leq -i$:
$$
b^l_{i,j}(\ta,\si) = d^l_{-i,-j}(-\ta,-\si) b^{-i}_{i,j}(\ta,\si)
p^{(j-i,-i-j)}_{l+i}(\rts;q^{-\ta}, q^{-\si} \mid q^2).
$$

\noindent
{\rm (iv)} In case $m=-j$ or $j\leq i\leq -j$:
$$
b^l_{i,j}(\ta,\si) = d^l_{-j,-i}(-\si,-\ta) b^{-j}_{i,j}(\ta,\si)
p^{(i-j,-i-j)}_{l+j}(\rts;q^{-\si}, q^{-\ta} \mid q^2).
$$

\noindent
Here the constant is given by
$$
d^l_{i,j}(\ta,\si) = {{C^{l,j}(\si)C^{l,i}(\ta)}
\over{C^{i,j}(\si)C^{i,i}(\ta)}} 
{{q^{i-l}}\over{(q^{4l};q^{-2})_{l-i}}} .
$$
\endproclaim

\demo{Proof} The explicit form for the orthogonality measure given
in the proof of Theorem \thmref{\thmabstrorthpolgenmatelt},
the symmetry relations of \thetag{\vglsymmetriesgenmatelt} and
$\pi_{\th/2}(\rts) =\cos\theta$ being independent of $\si$, $\ta$, 
show that it suffices to prove the first statement.

{}From the explicit form of the Askey-Wilson weight measure, cf.
Proposition \thmref{\propAWorthogonality},
we immediately get
$$
(az,a/z;q)_r \, dm(x;a,b,c,d| q) = {{ (ab,ac,ad;q)_r}
\over{(abcd;q)_r}}\, dm(aq^r,b,c,d| q)
$$
for $r\in\Zp$, $x=(z+z^{-1})/2$. A double application shows that 
for $r,s\in\Zp$, $x=(z+z^{-1})/2$,
$$
\multline
(az,a/z;q)_r (dz,d/z;q)_s\,  dm(x;a,b,c,d| q) = \\
(ab,ac;q)_r(bd,cd;q)_s {{(ad;q)_{r+s}}\over{(abcd;q)_{r+s}}} \,
dm(aq^r,b,c,dq^s| q).
\endmultline
$$
Hence, Proposition \thmref{\propexplvalueonedimonmingenmatelt}(i)
and Theorem \thmref{\thmHaarongeneralspherelt}
imply that the looked-for polynomials
are multiples of the Askey-Wilson polynomials
$$
p_{l-i}(\rts;-q^{\si+\ta+1+2i+2j},
-q^{-\si-\ta+1}, q^{\si-\ta+1}, q^{\ta-\si+1+2i-2j}| q^2),
$$
which we rewrite in the shorthand notation.

It remains to calculate the constant. We apply the one-dimensional
$\ast$-representation $\pi_{\th/2}$ to both sides of (i), and next 
we compare the coefficient of $e^{-il\theta}$ on both sides. The
coefficient of $e^{-il\theta}$ on the left hand side is
$v^{l,i}_l(\si)\overline{v^{l,j}_l(\ta)}q^{-l} = C^{l,j}(\si)
C^{l,i}(\ta)q^{-l}$. The coefficient of $e^{-i(l-i)\theta}$ of 
$p_{l-i}$ is $(q^{2l+2i+2};q^2)_{l-i}= (q^{4l};q^{-2})_{l-i}$, so
that the coefficient of $e^{-il\theta}$ on the right hand side
equals $C^{i,j}(\si)C^{i,i}(\ta)q^{-i}(q^{4l};q^{-2})_{l-i}$, from 
which we obtain the value for $d^l_{i,j}(\ta,\si)$.
\qed\enddemo

Since the generalised matrix elements $a^l_{ij}(\ta,\si)$,
see Lemma \thmref{\lemintrogenmatrixelts}, are obtained from
$b^l_{ij}(\ta,\si)$ by applying the simple
algebra isomorphism $D.$, we also have explicit expressions
for the generalised matrix elements in terms of
Askey-Wilson polynomials.

%%%%%%%%%%%%%%%%%%%%%%%%%%%%%%%%%%%%%%%%%%%%%%%%%%%%%%%%%%%%%%%%%%%%
%%N E W   S U B S E C T I O N%%%%%%%%%%%%%%%%%%%%%%%%%%%%%%%%%%%%%%%
%%%%%%%%%%%%%%%%%%%%%%%%%%%%%%%%%%%%%%%%%%%%%%%%%%%%%%%%%%%%%%%%%%%%
\subhead\newsubsection 
Limit cases
\endsubhead
Theorem \thmref{\thmdefgenmateltAskeyWpol} remains valid
for the limiting cases $\si\to\infty$ or $\ta\to\infty$ or even
$\si=\ta\to\infty$, cf. Remark 
\thmref{\rempropsphericaleltsalgebra}.
If we let either $\si$ or $\ta$ tend to
infinity, we see a similar expression with the Askey-Wilson
polynomials replaced by the big $q$-Jacobi polynomials, which
are defined by
$$
P_n^{(\al,\be)}(x;c,d;q) = {}_3\vp_2\left(
{{q^{-n},q^{n+\al+\be+1},xq^{\al+1}/c}\atop{q^{\al+1},\, 
-q^{\al+1}d/c}};q,q\right).
\tag\eqname{\vgldefbigqJacobi}
$$
This is due to the limit transition
$$
\multline
\lim_{a\to 0} (a\sqrt{cq/d})^n (-q/a^2;q)_n\, p_n^{(\al,\be)}
\bigl({{xq^{1/2}}\over{2a\sqrt{cd}}};a, \sqrt{ {c\over d}}|q) =\\
(q^{\al+1}, -q^{\be+1}c/d;q)_n\, P_n^{(\al,\be)}(x;c,d;q),
\endmultline
$$
which follows by taking term-wise limits in the ${}_4\vp_3$-series
expression for the Askey-Wilson polynomials.
This limit transition is
motivated from the definition of $\rti$ and $\ris$ in
Remark \thmref{\rempropsphericaleltsalgebra}.

For $\si=\ta\to\infty$ we need the limit transition
of the Askey-Wilson polynomials to little $q$-Jacobi
polynomials defined by
$$
\aligned
p_n^{(\al,\be)}(x;q) &= {}_2\vp_1(q^{-n}, q^{n+\al+\be+1};
q^{\al+1};q,qx) \\
&= {{(q^{-n-\be};q)_n}\over{(q^{\al+1};q)_n}}\,
{}_3\vp_2\left( {{q^{-n}, q^{n+\al+\be+1}, xq^{\be+1}}
\atop{q^{\be +1},\ 0}}; q,q\right).
\endaligned
\tag\eqname{\vgldeflittleqJacobi}
$$
Using the ${}_3\vp_2$-series representation for the little
$q$-Jacobi polynomials we can prove the limit transition
$$
\lim_{a\to 0} (-q^{1/2+\be}a^2)^n\, p_n^{(\al,\be)}\bigl(
{{xq^{1/2}}\over{2a^2}};a,a|q) = (-1)^n q^{n\be}q^{n(n-1)/2}
(q^{\al+1};q)_n\, p_n^{(\al,\be)}(x;q).
$$
This limit transition is motivated from the limit transition of
$\rts$ to $\rii$ in
Remark \thmref{\rempropsphericaleltsalgebra}.
This limit case is stated separately. Note that it suffices
to do the calculations in the first case and then use the symmetry
relations of Remark \thmref{\rempropexplvalueonedimonmingenmatelt}
and the commutation relations in $\Asu$.

\proclaim{Corollary \theoremname{\cortlnmaslittleqJacobipols}}
$$
\aligned
t^l_{n,m} & = c^l_{n,m} \de^{n+m} \ga^{n-m}\,
p_{l-n}^{(n-m,n+m)}(-q^{-1}\be\ga;q^2), \qquad
(n\geq m\geq -n), \\
t^l_{n,m} & = c^l_{m,n} \de^{n+m} \be^{m-n} \,
p_{l-m}^{(m-n,m+n)}(-q^{-1}\be\ga;q^2), \qquad
(m\geq n\geq -m), \\
t^l_{n,m} & = c^l_{-n,-m} \be^{m-n} \al^{-m-n} \,
p_{l+n}^{(m-n,-n-m)}(-q^{2m+2n-1}\be\ga;q^2),
\qquad (-n\geq m\geq n), \\
t^l_{n,m} & = c^l_{-m,-n} \ga^{n-m} \al^{-m-n}\,
p_{l+m}^{(n-m,-n-m)}(-q^{2m+2n-1}\be\ga;q^2),
\qquad (-m\geq n\geq m),
\endaligned
$$
with $p_k^{(\al,\be)}(x;q^2)$ a little
$q$-Jacobi polynomial and
$$
c^l_{n,m} = \left[ {{l-m}\atop{n-m}}\right]_{q^2}^{1/2}
\left[ {{l+n}\atop{n-m}}\right]_{q^2}^{1/2} q^{-(n-m)(l-n)}.
$$
\endproclaim

%%%%%%%%%%%%%%%%%%%%%%%%%%%%%%%%%%%%%%%%%%%%%%%%%%%%%%%%%%%%%%%%%%%%
%%N E W   S U B S E C T I O N%%%%%%%%%%%%%%%%%%%%%%%%%%%%%%%%%%%%%%%
%%%%%%%%%%%%%%%%%%%%%%%%%%%%%%%%%%%%%%%%%%%%%%%%%%%%%%%%%%%%%%%%%%%%
%NOTES AND REFERENCES%%%%%%%%%%%%%%%%%%%%%%%%%%%%%%%%%%%%%%%%%%%%%%%
\subhead Notes and references
\endsubhead
The main Theorem \thmref{\thmdefgenmateltAskeyWpol} is due to
Koornwinder \cite{\KoorZSE} and was already announced in 1990,
see \cite{\KoorOPTA}, in the case $l\in\Zp$ and $i=j=0$, i.e.
in the case of $(\ta,\si)$-spherical elements. This has led
to Koelink \cite{\KoelSIAM} in which the case $j=0$ is calculated
in order to be able to give a quantum group theoretic proof
of the addition formula for continuous $q$-Legendre polynomials.
Noumi and Mimachi \cite{\NoumMPJA}, \cite{\NoumMunpub} then
announced Theorem \thmref{\thmdefgenmateltAskeyWpol} in general.
There are some other ways of proving
Theorem \thmref{\thmdefgenmateltAskeyWpol}, notably by
using the Casimir operator leading to the $q$-difference
equation for Askey-Wilson polynomials \cite{\KoorZSE} or
by the Clebsch-Gordan decomposition leading to the
three-term recurrence relation for the Askey-Wilson polynomials.
This proof uses the techniques of \cite{\KoelSIAM}, see
also \cite{\KoelAAM}.
The symmetry relations in
Remark \thmref{\rempropexplvalueonedimonmingenmatelt}
can already be found in \cite{\KoorIM}.

The limit transitions of the Askey-Wilson polynomials to
big and little $q$-Jacobi polynomials is taken from
Koornwinder \cite{\KoorZSE}. See \cite{\GaspR},
\cite{\KoekS}, \cite{\KoorTrento}
for more information and further references
on the big and little $q$-Jacobi
polynomials, which were originally introduced by Andrews and
Askey.
Corollary \thmref{\cortlnmaslittleqJacobipols} is one of the
first known interactions between $q$-special functions and
quantum groups due to Koornwinder \cite{\KoorIM},
Masuda et al. \cite{\MasuMNNU} and Vaksman and
Soibelman \cite{\VaksS}. A number of other special cases
of Theorem \thmref{\thmdefgenmateltAskeyWpol} in the
case $j=0$ have been considered
before, see Noumi and Mimachi \cite{\NoumMCMP}, \cite{\NoumMDMJ},
\cite{\NoumMCompM}.

%%%%%%%%%%%%%%%%%%%%%%%%%%%%%%%%%%%%%%%%%%%%%%%%%%%%%%%%%%%%%%%%%%%%
%%N E W   S U B S E C T I O N%%%%%%%%%%%%%%%%%%%%%%%%%%%%%%%%%%%%%%%
%%%%%%%%%%%%%%%%%%%%%%%%%%%%%%%%%%%%%%%%%%%%%%%%%%%%%%%%%%%%%%%%%%%%
%EXERCISES%%%%%%%%%%%%%%%%%%%%%%%%%%%%%%%%%%%%%%%%%%%%%%%%%%%%%%%%%%
\subhead Exercises
\endsubhead

\item{\the\sectionno.1} Prove that $h_l$ of 
\thetag{\vgldefinitionhsubl} can be written as
$$
{{(1-q^2)q^{2l}}\over{1-q^{4l+2}}}
\ta_{\sqrt{q}}\bigr( b^l_{j,j}(\si,\si)\bigl)
\ta_{\sqrt{q}}\bigr( b^l_{i,i}(\ta,\ta)\bigl).
$$

\item{\the\sectionno.2} Use
Proposition \thmref{\propexplvalueonedimonmingenmatelt}
to derive the following generating function for the dual
$q$-Krawtchouk polynomials of Exercise~3.3;
$$
\align
\sum_{n=0}^N &t^n q^{n(N+\si)/2} {{(q^{-N};q)_n}\over{(q;q)_n}}
R_n(q^{-x}-q^{x-N-\si};q^\si,N;q) \\
&= (-tq^{-(N+\si)/2};q)_x (tq^{(\si-N)/2};q)_{N-x}
\endalign
$$
for $N\in\N$, $x\in\{ 0,\ldots, N\}$.

\item{\the\sectionno.3} Prove the limit transitions of the
Askey-Wilson polynomials to the big and little $q$-Jacobi
polynomials.

\item{\the\sectionno.4} Prove the second equality in
\thetag{\vgldeflittleqJacobi}.

\item{\the\sectionno.5} Prove the orthogonality relations
for the little $q$-Jacobi polynomials;
$$
\int_0^1 p_n^{(\al,\be)}(x;q)p_m^{(\al,\be)}(x;q)
\, t^\al {{(qt;q)_\infty}\over{(q^{\be+1};q)_\infty}}\, d_qt =
\de_{nm}h_n
$$
and calculate the squared norm $h_n$.

%%%%%%%%%%%%%%%%%%%%%%%%%%%%%%%%%%%%%%%%%%%%%%%%%%%%%%%%%%%%%%%%%%%%
%%N E W   S E C T I O N%%%%%%%%%%%%%%%%%%%%%%%%%%%%%%%%%%%%%%%%%%%%%
%%%%%%%%%%%%%%%%%%%%%%%%%%%%%%%%%%%%%%%%%%%%%%%%%%%%%%%%%%%%%%%%%%%%
\newpage

\head\newsection Addition formulas for Askey-Wilson polynomials
\endhead

As an application of the interpretation of Askey-Wilson polynomials
on the quantum $SU(2)$ group established in the previous section
we derive two addition formulas for the $q$-Legendre
polynomial $p_n^{(0,0)}(\cdot;q^\ta,q^\si|q^2)$ in this section.

%%%%%%%%%%%%%%%%%%%%%%%%%%%%%%%%%%%%%%%%%%%%%%%%%%%%%%%%%%%%%%%%%%%%
%%N E W   S U B S E C T I O N%%%%%%%%%%%%%%%%%%%%%%%%%%%%%%%%%%%%%%%
%%%%%%%%%%%%%%%%%%%%%%%%%%%%%%%%%%%%%%%%%%%%%%%%%%%%%%%%%%%%%%%%%%%%
\subhead\newsubsection 
Abstract addition formulas
\endsubhead
Since $t^l$ defines a unitary representation of the $\ast$-algebra
$\Usu$, we obtain the properties
$$
\gather
\De(t^l_{nm}) = \sum_{k=-l}^l t^l_{nk}\otimes t^l_{km},
\quad \ep(t^l_{nm})=\de_{nm},
\quad S(t^l_{nm})=(t^l_{mn})^\ast, \\
\sum_{p=-l}^l t^l_{ip}(t^l_{jp})^\ast = \de_{ij}
= \sum_{p=-l}^l (t^l_{pi})^\ast t^l_{pj},
\endgather
$$
which are easily verified by testing against appropriate elements
of $\Usu$. The analogous statements for the generalised
matrix elements is the following.

\proclaim{Proposition \theoremname{\propHopfpropgenmatelt}}
The elements $a^l_{i,j}(\ta,\si)$, $\si,\ta\in\R\cup\{\infty\}$,
$i,j=-l,-l+1,\ldots,l$,
$l\in\hZp$, defined in Lemma
\thmref{\lemintrogenmatrixelts}, satisfy
$$
\gather
\De\bigl(a^l_{i,j}(\ta,\si)\bigr) = \sum_{p=-l}^l a^l_{i,p}(\ta,\mu)
\otimes a^l_{p,j}(\mu,\si), \quad \forall\,\mu\in\R\cup\{\infty\},
 \\
\bigl(a^l_{i,j}(\ta,\si)\bigr)^\ast = 
S\bigl(a^l_{j,i}(\si,\ta)\bigr),\qquad
\ep\bigl(a^l_{i,j}(\ta,\si)\bigr) = \langle v^{l,j}(\si),
v^{l,i}(\ta)\rangle ,\\
\sum_{p=-l}^l a^l_{i,p}(\ta,\mu)
\bigl( a^l_{j,p}(\si,\mu)\bigr)^\ast =
\langle v^{l,j}(\si),v^{l,i}(\ta)\rangle =
\sum_{p=-l}^l \bigl(a^l_{p,i}(\mu,\ta)\bigr)^\ast 
a^l_{p,j}(\mu,\si).
\endgather
$$
\endproclaim

\demo{Proof} These statements are proved by testing against 
appropriate elements. Firstly,
$$
\align
\langle X\otimes Y, &\, \De\bigl(a^l_{i,j}(\ta,\si)\bigr) \rangle =
\langle XY, a^l_{i,j}(\ta,\si)\rangle =
\langle t^l(X)t^l(Y) v^{l,j}(\si), v^{l,i}(\ta)\rangle \\
&= \sum_{p=-l}^l \langle t^l(X)v^{l,p}(\mu), v^{l,i}(\ta)\rangle
\langle t^l(Y)v^{l,j}(\si),v^{l,p}(\mu)\rangle \\
&= \sum_{p=-l}^l \langle X, a^l_{i,p}(\ta,\mu)\rangle
\langle Y, a^l_{p,j}(\mu,\si) \rangle , \qquad\forall\, X,Y\in\Usu,
\endalign
$$
by developing $t^l(Y)v^{l,j}(\si)$ in the orthonormal basis $\{
v^{l,p}(\mu)\}_{p=-l,\ldots,l}$.
The next statement follows from
$$
\align
&\langle X,\bigl(a^l_{i,j}(\ta,\si)\bigr)^\ast\rangle =
\overline{\langle S(X)^\ast, (a^l_{i,j}(\ta,\si)\rangle} =
\overline{\langle t^l(S(X)^\ast)v^{l,j}(\si),v^{l,i}(\ta)\rangle}\\
&=\langle t^l(S(X))v^{l,i}(\ta),v^{l,j}(\si)\rangle
=\langle S(X), a^l_{j,i}(\si,\ta)\rangle
=\langle X, S\bigl(a^l_{j,i}(\si,\ta)\bigr)\rangle
\endalign
$$
for arbitrary $X\in\Usu$ and
$$
\ep\bigl( a^l_{i,j}(\ta,\si)\bigr) =
\langle 1, a^l_{i,j}(\ta,\si)\rangle
= \langle v^{l,j}(\si), v^{l,i}(\ta)\rangle.
$$
To prove the last statement we apply $(id\otimes S)$ to
the first statement and we use the Hopf algebra axiom
$m\circ(id\otimes S)=\et\circ\ep$, cf. Definition
\thmref{\defHopfalgebra}.
The first equality then follows from
the second property. The second equality is proved similarly
using the map $m\circ (S\otimes id)$.
\qed\enddemo

\proclaim{Corollary \theoremname{\corpropHopfpropgenmatelt}}
For $l\in\Zp$, $\si,\ta\in\R\cup\{ \infty\}$ the action of
$\De$ on
$(\ta,\si)$-spherical elements is given by
$$
\Delta\bigl(b^l_{00}(\ta,\si)\bigr) = \sum_{n=-l}^l
\bigl( D.b^l_{0n}(\ta,\mu)\bigr) \otimes b^l_{n0}(\mu,\si),
\qquad\forall\;\mu\in\R\cup\{\infty\}.
$$
\endproclaim

\demo{Proof} Use Proposition \thmref{\propHopfpropgenmatelt},
$a^l_{ij}(\ta,\si)=D.b^l_{ij}(\ta,\si)$ and
$\De\circ Z. = (id\otimes Z.)\circ\De$. This follows from
$$
\langle \De (Z.\xi),X\otimes Y\rangle = \langle \xi, XYZ\rangle
= \langle \De X\otimes YZ\rangle = \langle (id\otimes Z.)\De(\xi),
X\otimes Y\rangle,
$$
where we use $\langle Z.\xi,X\rangle = \langle \xi,XZ\rangle$.
\qed\enddemo

Corollary \thmref{\corpropHopfpropgenmatelt} is the starting
point for the derivation of an addition formula for
a two-parameter family of Askey-Wilson polynomials, and in this
sense we may call it an abstract addition formula.
Most of
the ingredients of the proof have already been established.
But let us first note that applying the one-dimensional
representation $\ta_{\sqrt{\la}}\otimes \ta_{\sqrt{\nu}}$ to
the identity in $\Asu\otimes\Asu$ of Corollary
\thmref{\corpropHopfpropgenmatelt} leads to the
addition formula
$$
\aligned
&(q^2;q^2)_l q^{-l} p_l^{(0,0)}(\xi(\la\nu);q^\ta,q^\si| q^2) =
{{p_l^{(0,0)}(\xi(\la);q^\mu,q^\ta| q^2)
p_l^{(0,0)}(\xi(\nu);q^\mu,q^\si|q^2)}\over
{(-q^{2-2\mu},-q^{2+2\mu};q^2)_l}} \\
& + \sum_{p=1}^l {{(1+q^{4p+2\mu})(q^2;q^2)_{l+p} (\la\nu)^{-p}
(\la q^{\mu-\ta},-\la q^{\ta+\mu},\nu q^{\mu-\si}, -\nu
q^{\mu+\si};q^2)_p}\over{(1+q^{2\mu}) (q^2;q^2)_{l-p}
(-q^{2-2\mu};q^2)_{l-p} (-q^{2+2\mu};q^2)_{l+p}}} \\
&\qquad\qquad \times
p_{l-p}^{(p,p)}(\xi(\la);q^\mu,q^\ta| q^2)
p_{l-p}^{(p,p)}(\xi(\nu);q^\mu,q^\si|q^2) \\
& + \sum_{p=1}^l {{(1+q^{4p-2\mu})(q^2;q^2)_{l+p} (\la\nu)^{-p}
(\la q^{\ta-\mu},-\la q^{-\ta-\mu},\nu q^{\si-\mu}, -\nu
q^{-\mu-\si};q^2)_p}\over{(1+q^{-2\mu}) (q^2;q^2)_{l-p}
(-q^{2+2\mu};q^2)_{l-p} (-q^{2-2\mu};q^2)_{l+p}}}\\
&\qquad\qquad \times
p_{l-p}^{(p,p)}(\xi(\la);q^{-\mu},q^{-\ta}|q^2)
p_{l-p}^{(p,p)}(\xi(\nu);q^{-\mu},q^{-\si}|q^2),
\endaligned
\tag\eqname{\vgladdformqLegtwoparam}
$$
with $\xi(\la) = {1\over 2}(q^{-1}\la + q \la^{-1})$.
The one-dimensional
representation $\ta_{\sqrt{\la}}\otimes \ta_{\sqrt{\nu}}$
has a very large kernel, so we may expect that a more general
addition formula than \thetag{\vgladdformqLegtwoparam} holds.
This is indeed the case, and we proceed to derive this
for the special case $\mu\to\infty$.

%%%%%%%%%%%%%%%%%%%%%%%%%%%%%%%%%%%%%%%%%%%%%%%%%%%%%%%%%%%%%%%%%%%%
%%N E W   S U B S E C T I O N%%%%%%%%%%%%%%%%%%%%%%%%%%%%%%%%%%%%%%%
%%%%%%%%%%%%%%%%%%%%%%%%%%%%%%%%%%%%%%%%%%%%%%%%%%%%%%%%%%%%%%%%%%%%
\subhead\newsubsection 
Suitable basis for the right hand side
\endsubhead
Our plan is to derive an addition formula
from Corollary \thmref{\corpropHopfpropgenmatelt}
for the case $\mu\to\infty$. We first determine the explicit form
of the terms in the right hand side of
Corollary \thmref{\corpropHopfpropgenmatelt}.

\proclaim{Lemma \theoremname{\lemgenmateltsinftyasbigqJacobi}}
For $l\in\Zp$ fixed and $n\in\{ 0,\ldots,l\}$ we have
$$
\align
b^l_{n,0}(\infty,\si) &= C_n(\si)
\bigl( \de_{\infty,\si-1}\ga_{\infty,\si} \bigr)^n
P^{(n,n)}_{l-n}(\ris;q^{2\si},1 ; q^2), \\
b^l_{0,n}(\ta,\infty) &= C_n(\ta)
\bigl( \de_{\ta-1,\infty}\be_{\ta,\infty} \bigr)^n
P^{(n,n)}_{l-n}(\rti;q^{2\ta}, 1; q^2),\\
b^l_{-n,0}(\infty,\si) &= C_n(\si)
\bigl( \be_{\infty,\si-1}\al_{\infty,\si} \bigr)^n
P^{(n,n)}_{l-n}(\ris;q^{2\si+2n}, q^{2n}; q^2), \\
b^l_{0,-n}(\ta,\infty) &= C_n(\ta)
\bigl(\ga_{\ta-1,\infty}\al_{\ta,\infty}  \bigr)^n
P^{(n,n)}_{l-n}(\rti;q^{2\ta+2n}, q^{2n}; q^2)
\endalign
$$
with the constant given by
$$
C_n(\si) = q^{(l-n)(l-n-1)/2} q^{-\si l}
{{(q^{2n+2},-q^{2\si-2l};q^2)_{l-n}}\over{
\sqrt{ (q^2,q^{2l+2n+2};q^2)_{l-n}} }}.
$$
\endproclaim

\demo{Proof} Use Theorem \thmref{\thmdefgenmateltAskeyWpol}
and the limit transition of the Askey-Wilson polynomials
to big $q$-Jacobi polynomials in \S 5.3 to obtain
this form. The form of the factors in front follows from
Proposition \thmref{\propshifteigenval}.
\qed\enddemo

Theorem \thmref{\thmdefgenmateltAskeyWpol} and Lemma
\thmref{\lemgenmateltsinftyasbigqJacobi} give explicit
expressions in term of orthogonal polynomials for all
the elements in the identity in the non-commutative
algebra $\Asu\otimes \Asu$ of Corollary
\thmref{\corpropHopfpropgenmatelt}
for $\mu\to\infty$. Although identities in non-commuting
variables have their charm and are useful, we want
to obtain an identity in commuting variables.
In order to do so we use the $\ast$-representations of
$\Asu$ described in Theorem \thmref{\thmirredstarrepsofAsu}.
Now any one-dimensional $\ast$-representation leads
to a trivial identity, so we have to consider
infinite dimensional $\ast$-representations. Since the
result is independent of the $\th$ of $\pi^\infty_\th$,
we may restrict ourselves to $\th=0$, and we
denote $\pi=\pi^\infty_0$. The idea is to let the right hand side
act on suitable eigenvectors of the self-adjoint
operator $\pi(\rti)\otimes \pi(\ris)$,
and to determine the (generalised) eigenvectors of
the self-adjoint operator
$\pi\otimes\pi (\De(\rts))$ in terms of the eigenvectors
of $\pi(\rti)\otimes \pi(\ris)$.

\proclaim{Proposition \theoremname{\propeigvectprti}}
$\Hi$ has an orthogonal basis of
eigenvectors $v_\la(q^\ta)$, where $\la=-q^{2n}$, $n\in\Zp$, 
$\la=q^{2\ta+2n}$, $n\in\Zp$, for the eigenvalue $\la$ of the 
self-adjoint operator $\pi(\rti)$.
For $\la=-q^{2n}$, $\la=q^{2\ta+2n}$, $n\in\Zp$, we have
$$
\gather
\pi(\al_{\ta,\infty})v_\la(q^\ta)
= iq^{1/2-\ta}(1+\la) v_{\la/q^2}(q^{\ta-1}), \qquad
\pi(\be_{\ta,\infty})v_\la(q^\ta)
= iq^{1/2} v_\la (q^{\ta-1}),\\
\pi(\ga_{\ta,\infty})v_\la(q^\ta)
= iq^{1/2} (q^{2\ta}-\la)v_{\la}(q^{\ta+1}),\qquad
\pi(\de_{\ta,\infty})v_\la(q^\ta)
= -iq^{1/2+\ta} v_{\la q^2}(q^{\ta+1}).
\endgather
$$
\endproclaim

\demo{Proof} The first statement is Proposition 
\thmref{\propeigvectpirti} for $\th=0$, and the second statement 
is Exercise~4.4 for $\th=0$. This can be proved
as follows. First observe that the result for 
$\pi(\be_{\ta,\infty})$ and $\pi(\de_{\ta,\infty})$ imply the 
result for $\pi(\ga_{\ta,\infty})$ and $\pi(\al_{\ta,\infty})$ 
by Proposition \thmref{\propeigvectpirti} and the case
$\si\to\infty$ of Proposition \thmref{\propsphericaleltsalgebra}.

Let us prove the statement for $\pi(\de_{\ta,\infty})$. Take
$\si\to\infty$ in Corollary \thmref{\corcommrst}
using Remark \thmref{\rempropsphericaleltsalgebra}, to
see $\de_{\ta,\infty}\rti = q^{-2}\rho_{\ta+1,\infty}
\de_{\ta,\infty}$. By the first part of the proposition we must 
have $\pi(\de_{\ta,\infty})v_\la(q^\ta)=Cv_{\la q^2}(q^{\ta+1})$ 
for some constant $C$. Using the normalisation
$\langle v_\la(q^\ta),e_0\rangle = 1$ we get
$$
\multline
C = \langle \pi(\de_{\ta,\infty})v_\la(q^\ta), e_0\rangle =
\langle v_\la(q^\ta),\pi(\de_{\ta,\infty^\ast}) e_0\rangle = \\
\langle v_\la(q^\ta),\pi(iq^{\ta+1/2}\ga + q^{-1/2}\al) e_0\rangle
 = - iq^{\ta+1/2}.
\endmultline
$$
The statement for $\pi(\be_{\ta,\infty})$ is proved analogously.
\qed\enddemo

\proclaim{Proposition \theoremname{\propeigvectpris}}
$\Hi$ has an orthogonal basis of
eigenvectors $v_\la(q^\si)$, where $\la=-q^{2n}$, $n\in\Zp$, and
$\la = q^{2\si + 2n}$,
$n\in\Zp$, for the eigenvalue $\la$ of the self-adjoint operator 
$\pi(\ris)$. Moreover,
$$
\align
\langle v_\la(q^\si),v_\la(q^\si)\rangle &= q^{-2n} (q^2;q^2)_n
(-q^{2-2\si};q^2)_n (-q^{2\si};q^2)_\infty,\qquad \la=-q^{2n},\\
\langle v_\la(q^\si),v_\la(q^\si)\rangle &= q^{-2n} (q^2;q^2)_n
(-q^{2+2\si};q^2)_n (-q^{-2\si};q^2)_\infty,\qquad \la=q^{2\si+2n}.
\endalign
$$
 For $\la=-q^{2n}$, $\la=q^{2\si+2n}$, $n\in\Zp$, we have
$$
\gather
\pi(\al_{\infty,\si})v_\la(q^\si) =
iq^{1/2-\si}(1+\la) v_{\la/q^2}(q^{\si-1}), \quad
\pi(\be_{\infty,\si})v_\la(q^\si) = iq^{1/2} (q^{2\si}-\la) 
v_\la(q^{\si+1}),\\
\pi(\ga_{\infty,\si})v_\la(q^\si) = iq^{1/2} v_\la(q^{\si-1}),
\qquad \pi(\de_{\infty,\si})v_\la(q^\si) = -iq^{1/2+\si} 
v_{\la q^2}(q^{\si+1}).
\endgather
$$
\endproclaim

\demo{Proof}
Use $-q\pi(\ga)=\pi(\be)$ to find  $\pi(\rts)=\pi(\rho_{\si,\ta})$,
$\pi(\al_{\ta,\si})=\pi(\al_{\si,\ta})$,
$\pi(\be_{\ta,\si})=\pi(\ga_{\si,\ta})$,
$\pi(\ga_{\ta,\si})=\pi(\be_{\si,\ta})$ and
$\pi(\de_{\ta,\si})=\pi(\de_{\si,\ta})$,
to reduce to Propositions \thmref{\propeigvectpirti} and
\thmref{\propeigvectprti}.
\qed\enddemo

\proclaim{Corollary \theoremname{\corexplicitactionbnnuletc}}
For $n\in\{ 0,\ldots,l\}$ we have
$$
\align
\pi(b^l_{n,0}(\infty,\si)) v_\la(q^\si) &=
C_n(\si) q^{n\si} P^{(n,n)}_{l-n}(\la;q^{2\si},1 ; q^2)\,
v_{\la q^{2n}}(q^\si),\\
\pi(D.b^l_{0,n}(\ta,\infty))  v_\mu(q^\ta) &=
C_n(\ta) q^{n(\ta+1)} P^{(n,n)}_{l-n}(\mu;q^{2\ta},1 ; q^2)
\, v_{\mu q^{2n}}(q^\ta)\\
\intertext{and}
\pi(b^l_{-n,0}(\infty,\si)) v_\la(q^\si) &=
C_n(\si)(-1)^n q^{n(\si-1)} (-\la,\la q^{-2\si};q^{-2})_n \\
&\qquad\qquad\qquad\times
P^{(n,n)}_{l-n}(\l q^{-2n};q^{2\si},1 ; q^2)\, 
v_{\la q^{-2n}}(q^\si),\\
\pi(D.b^l_{0,-n}(\ta,\infty)) v_\mu(q^\ta) &=
C_n(\ta)(-1)^n q^{n(\ta-2)} (-\mu,\mu q^{-2\ta};q^{-2})_n\\
&\qquad\qquad\qquad\times
P^{(n,n)}_{l-n}(\mu q^{-2n};q^{2\ta},1 ; q^2)\, 
v_{\mu q^{-2n}}(q^\ta).
\endalign
$$
\endproclaim

\demo{Proof} Use $D.b^{1/2}_{ij}(\ta,\infty)=q^j 
b^{1/2}_{ij}(\ta,\infty)$, so $D.\rti=\rti$, and
$P^{(\al,\be)}_n(Ax;Ac,Ad;q)=P^{(\al,\be)}_n(x;c,d;q)$ for $A>0$ 
and Propositions \thmref{\propeigvectprti}
and \thmref{\propeigvectpris} and Lemma 
\thmref{\lemgenmateltsinftyasbigqJacobi}.
\qed\enddemo

%%%%%%%%%%%%%%%%%%%%%%%%%%%%%%%%%%%%%%%%%%%%%%%%%%%%%%%%%%%%%%%%%%%%
%%N E W   S U B S E C T I O N%%%%%%%%%%%%%%%%%%%%%%%%%%%%%%%%%%%%%%%
%%%%%%%%%%%%%%%%%%%%%%%%%%%%%%%%%%%%%%%%%%%%%%%%%%%%%%%%%%%%%%%%%%%%
\subhead\newsubsection 
Suitable basis for the left hand side
\endsubhead
The basis $v_\mu(q^\ta)\otimes v_\la(q^\si)$ of $\Hi\otimes\Hi$ is
very well suited for the action of the right hand side under
$\pi\otimes\pi$, cf. Corollary \thmref{\corexplicitactionbnnuletc}.
Next we study $\pi\otimes\pi\bigl(\De(\rts)\bigr)$ with
respect to this basis.

\proclaim{Lemma \theoremname{\lemmaactionrtsinbasiseigvecrtiris}}
$$
\multline
2(\pi\otimes\pi)\De(\rts)v_\mu(q^\ta)\otimes v_\la(q^\si) =
 q^2v_{\mu q^2}(q^\ta)\otimes v_{\la q^2}(q^\si)
+ \\ q^{-2}(1+\la)(1+\mu)(1-\la q^{-2\si})(1-\mu q^{-2\ta})
v_{\mu q^{-2}}(q^\ta)\otimes v_{\la q^{-2}}(q^\si)\\
+ \Bigl( \la q^{1-\si}(q^{-\ta}-q^\ta)+\mu q^{1-\ta}
(q^{-\si}-q^\si) +\la\mu
q^{1-\ta-\si}(1+q^2)\Bigr) v_\mu(q^\ta)\otimes v_\la(q^\si).
\endmultline
$$
\endproclaim

\demo{Proof} Use
$$
\Delta\bigl(b^{1/2}_{ij}(\ta,\si)\bigr) = \sum_{n=-1/2}^{1/2}
q^n b^{1/2}_{in}(\ta,\infty) \otimes b^{1/2}_{nj}(\infty,\si),
$$
cf. Corollary \thmref{\corpropHopfpropgenmatelt}, to get
$$
\multline
\De(\be_{\ta+1,\si-1}\ga_{\ta,\si}) =
q^{-1}\al_{\ta+1,\infty}\ga_{\ta,\infty}\otimes
\be_{\infty,\si-1}\al_{\infty,\si} +
\al_{\ta+1,\infty}\de_{\ta,\infty}\otimes
\be_{\infty,\si-1}\ga_{\infty,\si} \\
+\be_{\ta+1,\infty}\ga_{\ta,\infty}\otimes
\de_{\infty,\si-1}\al_{\infty,\si}
+q\be_{\ta+1,\infty}\de_{\ta,\infty})\otimes
\de_{\infty,\si-1}\ga_{\infty,\si},
\endmultline
$$
Now the lemma follows from Propositions 
\thmref{\propsphericaleltsalgebra},
\thmref{\propeigvectprti} and \thmref{\propeigvectpris}.
\qed\enddemo

Lemma \thmref{\lemmaactionrtsinbasiseigvecrtiris} implies that
each subspace of the form
$$
\text{span}\{ v_{\mu q^{2m}}(q^\ta)\otimes v_{\la q^{2m}}
(q^\si)\mid m\in\Zp\}
$$
with either $\mu\in\{ -1,q^{2\ta}\}$ or $\la\in \{-1,q^{2\si}\}$
is invariant under $(\pi\otimes\pi)\De(\rts)$. We pick one,
say $\mu=-1$, $\la=-q^{2p}$ for $p\in\Zp$ fixed, and we call
this subspace $W\cong\Hi$. Using  $\mu=-1$, $\la=-q^{2p}$ and
in Lemma \thmref{\lemmaactionrtsinbasiseigvecrtiris} we
obtain a three-term recurrence relation, which can be solved
using a sub-class of the Askey-Wilson polynomials. If viewed
as $q$-Jacobi polynomials as in \S 5.2, the polynomials we need
are the Laguerre case of the $q$-Jacobi polynomials, i.e. we let
$\be\to\infty$.  So we define
$$
l_n^{(\al)}(x;s,t|q) = p_n(x;q^{1/2+\al}s/t,q^{1/2}t/s,
-q^{1/2}/(st), 0| q),.
$$
{}From the three-term recurrence relation for the Askey-Wilson 
polynomials we get
$$
\gather
2xl_n(x)=l_{n+1}(x) +
(1-q^n)(1-q^{\al+n})(1+q^ns^{-2})(1+q^{n+\al}t^{-2}) l_{n-1}(x)\\
+ q^n\Bigl( (t-t^{-1})q^{1/2} s^{-1} +
(s-s^{-1})q^{1/2+\al}t^{-1} +(1+q)q^{1/2+n+\al}s^{-1}t^{-1}\Bigr)
l_n(x),
\endgather
$$
where $l_n(x)=l_n^{(\al)}(x;s,t| q)$. The orthogonality measure 
for the $q$-Laguerre polynomials follows from
\thetag{\vglnormalisedAWmeasure} and we denote the
normalised orthogonality measure by
$$
dm^{(\al)}(\cdot;s,t| q) = dm(\cdot;q^{1/2+\al}s/t,q^{1/2}t/s,
-q^{1/2}/(st), 0| q).
$$.

\proclaim{Proposition \theoremname{\propintertwingpipirts}}
Define $\Lambda \colon W\to L^2(dm^{(p)}(\cdot;q^\ta,q^\si|q^2))$ by
$$
 v_{-q^{2m}}(q^\ta)\otimes v_{-q^{2m+2p}}(q^\si) \mapsto
q^{-2m-p} \sqrt{(q^2,-q^{2-2\si};q^2)_p 
(-q^{2\ta},-q^{2\si};q^2)_\infty} \ 
l_m^{(p)}(\cdot;q^\ta,q^\si|q^2),
$$
then $\Lambda$ is a unitary operator intertwining 
$(\pi\otimes\pi)\De(\rts)$ on $W$ with the multiplication operator 
on $L^2(dm^{(p)}(\cdot;q^\ta,q^\si|q^2))$.
\endproclaim

\demo{Proof} $(\pi\otimes\pi)\De(\rts)$ is a Jacobi matrix on $W$ 
and the recurrence relation can be matched with the three-term 
recurrence relation for the $q$-Laguerre polynomials 
$l_m^{(p)}(\cdot;q^\ta,q^\si|q^2)$, so we can use the technique 
of the proof of Theorem \thmref{\thmFavard}. The constants involved
follow from Propositions \thmref{\propAWorthogonality},
\thmref{\propeigvectpirti} and \thmref{\propeigvectpris}.
\qed\enddemo

%%%%%%%%%%%%%%%%%%%%%%%%%%%%%%%%%%%%%%%%%%%%%%%%%%%%%%%%%%%%%%%%%%%%
%%N E W   S U B S E C T I O N%%%%%%%%%%%%%%%%%%%%%%%%%%%%%%%%%%%%%%%
%%%%%%%%%%%%%%%%%%%%%%%%%%%%%%%%%%%%%%%%%%%%%%%%%%%%%%%%%%%%%%%%%%%%
\subhead\newsubsection 
Addition formula for Askey-Wilson polynomials
\endsubhead
We now have sufficient information on the action of both sides
of Corollary \thmref{\corpropHopfpropgenmatelt} for $\mu\to\infty$.

\proclaim{Theorem \theoremname{\thmfirstaddform}}
The following addition
formula holds for $l,m,p\in\Zp$, $\si,\ta\in\R$, $x\in\C$;
$$
\align
&p^{(0,0)}_l(x;q^\ta,q^\si|q^2) l_m^{(p)}(x;q^\ta,q^\si| q^2) = \\
& \sum_{n=0}^l D^{n,l}(\ta,\si)
P^{(n,n)}_{l-n}(-q^{2m};q^{2\ta},1;q^2)
P^{(n,n)}_{l-n}(-q^{2m+2p};q^{2\si},1;q^2)
l_{m+n}^{(p)}(x;q^\ta,q^\si|q^2)  \\
&+ \sum_{n=1}^l D^{n,l}(\ta,\si)
(q^{2m},q^{2m+2p},-q^{2m+2p-2\si},-q^{2m-2\ta};q^{-2})_n \\
&\qquad\times P^{(n,n)}_{l-n}(-q^{2m-2n};q^{2\ta},1;q^2)
P^{(n,n)}_{l-n}(-q^{2m+2p-2n};q^{2\si},1;q^2)
l_{m-n}^{(p)}(x;q^\ta,q^\si|q^2),
\endalign
$$
with the constant given by
$$
D^{n,l}(\ta,\si) = (-q^{2\si-2l},-q^{2\ta-2l};q^2)_{l-n}
{{(q^{2(l-n+1)};q^2)_n}\over{(q^2;q^2)_n}} (q^{2(n+1)};q^2)_l
q^{(l-n)(l-n-2\si-2\ta)}.
$$
\endproclaim

\demo{Proof} Let $\pi\otimes\pi$ act on
Corollary \thmref{\corpropHopfpropgenmatelt}, $\mu\to\infty$,
and restrict the action to $W$, which is also invariant under the
action of the right hand side by Corollary
\thmref{\corexplicitactionbnnuletc}, as it should be. Let the
resulting operator identity act on
$ v_{-q^{2m}}(q^\ta)\otimes v_{-q^{2m+2p}}(q^\si)$ and apply
$\Lambda$ of Proposition \thmref{\propintertwingpipirts}
to get as an identity in $ L^2(dm^{(p)}(\cdot;q^\ta,q^\si|q^2))$;
$$
\align
&q^{-2m-p} \sqrt{(q^2,-q^{2-2\si};q^2)_p (-q^{2\ta},
-q^{2\si};q^2)_\infty}\,
p^{(0,0)}_l(x;q^\ta,q^\si|q^2) l_m^{(p)}(x;q^\ta,q^\si| q^2) = \\
&\Lambda\Bigl( p^{(0,0)}_l(\pi\otimes\pi\bigl(\De(\rts)\bigr);
q^\ta,q^\si|q^2)
v_{-q^{2m}}(q^\ta)\otimes v_{-q^{2m+2p}}(q^\si) \Bigr) = \\
& \sum_{n=0}^l C_n(\si)C_n(\ta) q^{n(\ta+\si+1)}
P^{(n,n)}_{l-n}(-q^{2m};q^{2\ta},1;q^2)
P^{(n,n)}_{l-n}(-q^{2m+2p};q^{2\si},1;q^2) \\
&\qquad \times
\Lambda\bigl(v_{-q^{2m+2n}}(q^\ta)\otimes v_{-q^{2m+2p+2n}}
(q^\si)\bigr) + \\
& \sum_{n=1}^l C_n(\si)C_n(\ta) q^{n(\ta+\si-3)}
(q^{2m},q^{2m+2p},-q^{2m+2p-2\si},-q^{2m-2\ta};q^{-2})_n \\
& \qquad \times
 P^{(n,n)}_{l-n}(-q^{2m-2n};q^{2\ta},1;q^2)
P^{(n,n)}_{l-n}(-q^{2m+2p-2n};q^{2\si},1;q^2)\\
& \qquad \times
\Lambda\bigl(v_{-q^{2m-2n}}(q^\ta)\otimes v_{-q^{2m+2p-2n}}
(q^\si)\bigr).
\endalign
$$
Use Proposition \thmref{\propintertwingpipirts} and simplify
to get the results, which holds for all $x\in\C$ by continuity.
\qed\enddemo

\demo{Remark \theoremname{\remsymmaddform}}
It seems that in deriving Theorem \thmref{\thmfirstaddform} there
is some arbitrariness in the choice of the subspace $W$ of
$\Hi\otimes\Hi$, to which we have restricted our attention.
Apart from the choice of $p\in\Zp$ there are $8$ of such possible
choices for $W$, which lead to $8$ of such addition formulas.
These are obviously four by four equivalent by interchanging $\ta$
and $\si$ and noting that the $q$-Legendre polynomial involved is 
also invariant under such a change. The remaining four types of 
addition formulas differ, since the big $q$-ultraspherical 
polynomials are evaluated at other points of the spectrum. 
However, from
$$
P^{(\al,\be)}_n (-x;c,d;q) = (-q^{\al-\be}d/c)^n
{{(q^{\be+1},-q^{\be+1}c/d;q)_n}\over{(q^{\al+1},
-q^{\al+1}d/c;q)_n}} P^{(\be,\al)}_n (x;d,c;q),
\tag\eqname{\vglbigqJacflip}
$$
(which is a direct consequence of the orthogonality relations)
and the trivial relation
$P^{(\al,\be)}_n (Ax;Ac,Ad;q) = P^{(\al,\be)}_n (x;c,d;q)$, for 
$A>0$, we see that
$$
P_{l-n}^{(n,n)}(-q^m;q^{\ta},1;q) = (-q^{-\ta})^{l-n}
{{(-q^{n+1+\ta};q)_{l-n}}\over{(-q^{n+1-\ta};q)_{l-n}}}
P_{l-n}^{(n,n)}(q^{m-\ta};q^{-\ta},1;q).
\tag\eqname{\vglswitchpoints}
$$
If we use this transformation for the big $q$-Jacobi polynomials
in Theorem \thmref{\thmfirstaddform} and we next change $x$
into $-x$ using
$p_n^{(\al,\be)}(-x;s,t|q)=(-1)^n p_n^{(\al,\be)}(x;-s,t| q)$
for the $q$-Legendre polynomial, $\al=\be=0$, and for the 
$q$-Laguerre polynomial, $\be\to\infty$, and we change
$q^\ta$ to $-q^{-\ta}$, then we obtain the
same addition formula as if we had started off with the space
$W$ spanned by the vectors type
$v_{q^{2\ta+2m}}(q^\ta)\otimes v_{-q^{2m+2p}}(q^\si)$, $m\in\Zp$.
A similar approach can be used for the other big $q$-ultraspherical
polynomial, and a combination of both shows that Theorem
\thmref{\thmfirstaddform} contains all the other possibilities by 
symmetry considerations.
\enddemo

Using the orthogonality relations for the $q$-Laguerre polynomials,
cf. Proposition \thmref{\propAWorthogonality}, we can pick
out the term for $n=0$ from the right hand side of
Theorem \thmref{\thmfirstaddform} and
we obtain the following product formula
for big $q$-Legendre polynomials.

\proclaim{Corollary \theoremname{\corprodformeen}} For $l,n\in\Zp$,
$0\leq n\leq l$ we have the product formula
$$
\multline
P^{(0,0)}_l(-q^{2m};q^{2\ta},1;q^2)
P^{(0,0)}_l(-q^{2(m+p)};q^{2\si},1;q) = \\{1\over C}
\int_\R p^{(0,0)}_l(x;q^\ta,q^\si| q^2)
\Bigl( l_m^{(p)}(x;q^\ta,q^\si |q^2)
\Bigr)^2 \, dm^{(p)}(x;q^\ta,q^\si| q^2) ,
\endmultline
$$
with
$$
C = D^{0,l}(\ta,\si) (q^2,q^{2+2p},-q^{2-2\ta},-q^{2+2p-2\si};q)_m.
$$
\endproclaim

Product formulae for big $q$-Legendre polynomials at other points
of the spectrum are obtained from \thetag{\vglswitchpoints}.

%%%%%%%%%%%%%%%%%%%%%%%%%%%%%%%%%%%%%%%%%%%%%%%%%%%%%%%%%%%%%%%%%%%%
%%N E W   S U B S E C T I O N%%%%%%%%%%%%%%%%%%%%%%%%%%%%%%%%%%%%%%%
%%%%%%%%%%%%%%%%%%%%%%%%%%%%%%%%%%%%%%%%%%%%%%%%%%%%%%%%%%%%%%%%%%%%
%NOTES AND REFERENCES%%%%%%%%%%%%%%%%%%%%%%%%%%%%%%%%%%%%%%%%%%%%%%%
\subhead Notes and references
\endsubhead
There are much more applications of the interpretation of 
Askey-Wilson polynomials on the quantum $SU(2)$ group than just 
addition formulas, see e.g. \cite{\KoelAAM}.
The addition formula for the Askey-Wilson polynomials in
\thetag{\vgladdformqLegtwoparam} is due to Noumi
and Mimachi \cite{\NoumMPJA}, see also \cite{\KoelSIAM},
\cite{\KoelAAM}. The addition formula derived in this section
is taken from \cite{\Koelunpub}, see also \cite{\KoelFIC}.
The motivation for choosing $\mu\to\infty$ is that we still
can obtain a polynomial identity.
In \cite{\KoelFIC} a very general identity, i.e. coming
from Corollary \thmref{\corpropHopfpropgenmatelt} for
general $\si$, $\ta$ and $\mu$, is derived, but for this
a non-symmetric Poisson kernel for Al-Salam and Chihara
polynomials is needed. The result from \cite{\KoelFIC}
contains Theorem \thmref{\thmfirstaddform} as a limit case.

Some applications of identities for $q$-special
functions in non-commuting variables, of which Lemma
\thmref{\lemqbinomiallemma} is a simple example,
as well as further references can be
found in Koornwinder \cite{\KoorMLI}.

The case $\si=\ta=\mu=0$ of Corollary 
\thmref{\corpropHopfpropgenmatelt}
has been used in \cite{\KoelSIAM} to derive the addition formula
for the continuous $q$-Legendre polynomial, which is a special
case of the addition formula for the continuous $q$-ultraspherical
polynomials analytically proved by Rahman and Verma \cite{\RahmV}.
{}From the case $\si=\ta=\mu\to\infty$ of
Corollary \thmref{\corpropHopfpropgenmatelt}
Koornwinder \cite{\KoorAF} has proved an addition formula for
the little $q$-Legendre polynomial, which can be obtained
from Theorem \thmref{\thmfirstaddform} as the limit case
$\si=\ta\to\infty$. Also the addition formula for big $q$-Legendre
polynomials of \cite{\KoelCJM} is a special case of
Theorem \thmref{\thmfirstaddform}.
The limit case $q\uparrow 1$ of Theorem \thmref{\thmfirstaddform}
to the addition formula for Legendre polynomials, cf.
Exerxise~6.1, is highly non-trivial. It uses the asymptotic
behaviour of a class of orthogonal polynomials including
the $q$-Laguerre polynomials, see Van Assche and Koornwinder
\cite{\VAsscK} for the general theorem and \cite{\Koelunpub}
for the explicit application to this case.
In \cite{\VAsscK} it is also shown how the weak
convergence of orthogonality measures can be used to show
that the product formulas as the one in Corollary 
\thmref{\corprodformeen}
tend to the product formula for the Legendre polynomials, cf.
\cite{\Koelunpub} for the details for this case.

%%%%%%%%%%%%%%%%%%%%%%%%%%%%%%%%%%%%%%%%%%%%%%%%%%%%%%%%%%%%%%%%%%%%
%%N E W   S U B S E C T I O N%%%%%%%%%%%%%%%%%%%%%%%%%%%%%%%%%%%%%%%
%%%%%%%%%%%%%%%%%%%%%%%%%%%%%%%%%%%%%%%%%%%%%%%%%%%%%%%%%%%%%%%%%%%%
%EXERCISES%%%%%%%%%%%%%%%%%%%%%%%%%%%%%%%%%%%%%%%%%%%%%%%%%%%%%%%%%%
\subhead Exercises
\endsubhead

\item{\the\sectionno.1} Let $R^{(\al,\be)}$ be the Jacobi polynomial
normalised by $R^{(\al,\be)}(1)=1$, i.e. $R^{(\al,\be)}(x) =$
$$
\sum_{k=0}^n {{(-n)_k(n+\al+\be+1)_k}\over{k!\, (\al+1)_k}}
\left({{1-x}\over 2}\right)^k =
{}_2F_1\left( {{-n,n+\al+\be+1}\atop{\al+1}}; {{1-x}\over 2}\right).
$$
The Legendre polynomial $R_n^{(0,0)}$ satisfies the addition formula
$$
\multline
 R_l^{(0,0)}\bigl( xy+t\sqrt{(1-x^2)(1-y^2)}\bigr)
= R^{(0,0)}_l(x)R^{(0,0)}_l(y) \\
+2\sum_{m=1}^l {{(l+m)!}\over{(l-m)!(m!)^2}}2^{-2m}
\biggl(\sqrt{(1-x^2)(1-y^2)}\biggr)^m R^{(m,m)}_{l-m}(x)
R^{(m,m)}_{l-m}(y)T_m(t),
\endmultline
$$
where $T_m(\cos\th)=\cos m\th$ is the Chebyshev polynomial of the
first kind. Prove that for $q\uparrow 1$ in
\thetag{\vgladdformqLegtwoparam} leads to the addition formula
for Legendre polynomials.

\item{\the\sectionno.2} Derive an addition formula for the little
$q$-Legendre polynomial. Do this by redoing the proof of Theorem
\thmref{\thmfirstaddform} using Corollary
\thmref{\cortlnmaslittleqJacobipols}, or by taking suitable
limits in Theorem \thmref{\thmfirstaddform}.

\item{\the\sectionno.3} Given the $q$-integral
$$
\int_{-d}^c {{(qx/c,-qx/d;q)_\infty}\over
{(q^{\al+1}x/c, -q^{\be+1}x/d;q)_\infty}}\, d_qx
= (1-q)c {{(q,-d/c,-qc/d, q^{\al+\be+2};q)_\infty}\over
{(q^{\al+1},q^{\be+1},-q^{\be+1}c/d,-q^{\al+1}d/c;q)_\infty}},
$$
prove the orthogonality relations for the big $q$-Jacobi polynomials
$$
\int_{-d}^c P^{(\al,\be)}_n (x;c,d;q) P^{(\al,\be)}_m(x;c,d;q)
{{(qx/c,-qx/d;q)_\infty}\over
{(q^{\al+1}x/c, -q^{\be+1}x/d;q)_\infty}}\, d_qx = \de_{nm} h_n
$$
Calculate $h_n$ and prove \thetag{\vglbigqJacflip}.

%%%%%%%%%%%%%%%%%%%%%%%%%%%%%%%%%%%%%%%%%%%%%%%%%%%%%%%%%%%%%%%%%%%%
%%N E W   S E C T I O N%%%%%%%%%%%%%%%%%%%%%%%%%%%%%%%%%%%%%%%%%%%%%
%%%%%%%%%%%%%%%%%%%%%%%%%%%%%%%%%%%%%%%%%%%%%%%%%%%%%%%%%%%%%%%%%%%%
\newpage

\head\newsection Convolution theorem for Al-Salam and Chihara
polynomials\endhead

In this section we consider the $\ast$-operator on $\U$ leading to
the real form $\Unc$, see Theorem \thmref{\thmstarstructureonU}.
Then we show how we can interpret Al-Salam and Chihara
polynomials as overlap coefficients similarly as $q$-Krawtchouk
polynomials for the real form $\Usu$. This interpretation can be
used to find a very general convolution theorem for the
Al-Salam and Chihara polynomials.

%%%%%%%%%%%%%%%%%%%%%%%%%%%%%%%%%%%%%%%%%%%%%%%%%%%%%%%%%%%%%%%%%%%%
%%N E W   S U B S E C T I O N%%%%%%%%%%%%%%%%%%%%%%%%%%%%%%%%%%%%%%%
%%%%%%%%%%%%%%%%%%%%%%%%%%%%%%%%%%%%%%%%%%%%%%%%%%%%%%%%%%%%%%%%%%%%
\subhead\newsubsection 
Positive discrete series representations of $\Unc$
\endsubhead
In Exercise~2.6 the positive discrete series representations
had to be calculated. We recall the result. For every
$k>0$ the representation $\pi_k$ of $\Unc$ acting in $\Hi$
is given by
$$
\pi_k(A)\, e_n =q^{k+n}\, e_n, \qquad
\pi_k(C)\, e_n = q^{1/2-k-n}
{{\sqrt{(1-q^{2n})(1-q^{4k+2n-2})}}\over{q-q^{-1}}}\, e_{n-1},
$$
where $\{ e_n\}_{n=0}^\infty$ is the standard orthonormal basis
of $\Hi$. The action of $B$ follows from $B=-C^\ast$,
$$
\pi_k(B) \, e_n =
q^{-1/2-k-n}
{{\sqrt{(1-q^{2n+2})(1-q^{4k+2n})}}\over{q^{-1}-q}}\, e_{n+1}.
$$
Note that $\pi_k(D)$ is an unbounded operator, but that
$\pi_k(A),\pi_k(B),\pi_k(C)\in {\Cal B}(\Hi)$. The operators that
we consider will be bounded. Obviously, $\pi_k$ is an
irreducible unitary (i.e. $\ast$-)representation of $\Unc$.

\proclaim{Lemma \theoremname{\lemmadecomptensorsunc}}
$$
\pi_{k_1}\otimes \pi_{k_2} \cong \bigoplus_{j=0}^\infty 
\pi_{k_1+k_2+j}
$$
\endproclaim

\demo{Sketch of Proof} The left hand side is a unitary 
representation of $\Unc$, hence completely reducible. The 
decomposition follows by counting lowest weight vectors, i.e. 
considering the kernel of $C$. See Exercise~7.1-2 for more details.
\qed\enddemo

Although the basis of the representation space $\Hi$ of $\pi_k$
is independent of $k$ we use the notation $e_n^k$ to stress
the $k$-dependence. Then Lemma \thmref{\lemmadecomptensorsunc}
implies that there exists a unitary matrix mapping
the orthogonal basis $e_{n_1}^{k_1}\otimes e_{n_2}^{k_2}$ onto
$e_n^{k_1+k_2+j}$ intertwining the action of $\Unc$ on both sides.
The matrix elements of this unitary mapping are
the Clebsch-Gordan coefficients.

\proclaim{Lemma \theoremname{\lemCGCforUnc}} The Clebsch-Gordan
coefficients are defined by
$$
e_n^k = \sum_{n_1,n_2=0}^\infty C^{k_1,k_2,k}_{n_1,n_2, n}\,
e_{n_1}^{k_1}\otimes e_{n_2}^{k_2},
$$
where $k=k_1+k_2+j$ for $j\in\Zp$. The sum is finite;
$n_1+n_2=n+j$. The Clebsch-Gordan coefficients are normalised
by $\langle e_0^k,e_0^{k_1}\otimes e_j^{k_2}\rangle >0$.
\endproclaim

\demo{Proof} Let $A$ act on both sides to get
$$
\sum_{n_1,n_2=0}^\infty C^{k_1,k_2,k}_{n_1,n_2, n}\,
q^{n+k}\,
e_{n_1}^{k_1}\otimes e_{n_2}^{k_2} =
q^{n+k}\, e_n^k =
\sum_{n_1,n_2=0}^\infty C^{k_1,k_2,k}_{n_1,n_2, n}\,
q^{k_1+n_1+k+_2+n_2}\,
e_{n_1}^{k_1}\otimes e_{n_2}^{k_2},
$$
so $C^{k_1,k_2,k}_{n_1,n_2, n}$ is zero unless
$n+k=k_1+n_1+k_2+n_2$ or $n_1+n_2=n+j$.

Using the action of $B$ and $C$ we can derive recurrence relations 
for the Clebsch-Gordan coefficients. Using 
$\De(C)=A\otimes C+C\otimes D$ we get
$$
\multline
q^{1/2-k-n}\sqrt{(1-q^{2n})(1-q^{4k+2n-2})} \, 
C^{k_1,k_2,k}_{n_1,n_2,n-1} = \\
q^{k_1+n_1-1/2-k_2-n_2}\sqrt{(1-q^{2n_2+2})(1-q^{4k_2+2n_2})} \,
C^{k_1,k_2,k}_{n_1,n_2+1,n} \\
+ q^{-1/2-k_1-n_1-k_2-n_2}\sqrt{(1-q^{2n_1+2})(1-q^{4k_1+2n_1})} \,
C^{k_1,k_2,k}_{n_1+1,n_2,n}.
\endmultline
\tag\eqname{\vglthreetermCGC}
$$
Using the action of $B$ we get a similar recursion for
$C^{k_1,k_2,k}_{n_1,n_2,n+1}$ showing that the Clebsch-Gordan
coefficients are completely determined by 
$C^{k_1,k_2,k}_{n_1,n_2,0}$. Take $n=0$ in 
\thetag{\vglthreetermCGC} to find a two-term recurrence
which can be solved by iteration;
$$
C^{k_1,k_2,k}_{n_1,n_2,0} = (-1)^{n_1} q^{2n_1k_1} q^{n_1(n_1-1)}
\sqrt{ {{(q^{2j}, q^{4k_2+2j-2};q^{-2})_{n_1}}\over
{(q^2,q^{4k_1};q^2)_{n_1}}}} \, C^{k_1,k_2,k}_{0,j,0}
$$
where $n_1+n_2=j$, $k=k_1+k_2+j$.
Since the transition matrix is unitary we have
$$
\multline
1 = \sum_{n_1+n_2=j} |C^{k_1,k_2,k}_{n_1,n_2,0}|^2 =
|C^{k_1,k_2,k}_{0,j,0}|^2 {}_2\vp_1\left( {{q^{-2j}, q^{2-2j-4k_2}}
\atop{q^{4k_1}}};q^2, q^{4k_1+4k_2+4j-2}\right) \\
= |C^{k_1,k_2,k}_{0,j,0}|^2
{{(q^{4k_1+4k_2+2j-2};q^2)_j}\over{(q^{4k_1};q^2)_j}}
\endmultline
$$
by the Gau\ss\ summation formula of Exercise~3.6.
This determines $C^{k_1,k_2,k}_{0,j,0}$ up to a phase factor, and
given the normalisation
$\langle e_0^k, e_0^{k_1}\otimes e_j^{k_2}\rangle >0$ we have
$C^{k_1,k_2,k}_{0,j,0}>0$ and its value as well
as the value of $C^{k_1,k_2,k}_{n_1,n_2,n}$ are completely 
determined.
\qed\enddemo

\demo{Remark \theoremname{\remlemCGCforUnc}} The value obtained
for $C^{k_1,k_2,k}_{0,j,0}$ and hence for 
$C^{k_1,k_2,k}_{n_1,n_2,0}$,
together with the three-term recurrence relation obtained from the
action of $B$ completely determine the Clebsch-Gordan
coefficients. This can be used to find an
explicit expression in terms of $q$-hypergeometric series.
We come back to the explicit expression in \S 7.3.
\enddemo

%%%%%%%%%%%%%%%%%%%%%%%%%%%%%%%%%%%%%%%%%%%%%%%%%%%%%%%%%%%%%%%%%%%%
%%N E W   S U B S E C T I O N%%%%%%%%%%%%%%%%%%%%%%%%%%%%%%%%%%%%%%%
%%%%%%%%%%%%%%%%%%%%%%%%%%%%%%%%%%%%%%%%%%%%%%%%%%%%%%%%%%%%%%%%%%%%
\subhead\newsubsection 
Twisted primitive elements and eigenvectors
\endsubhead
{}From Proposition \thmref{\propclassgroupprimelts}(ii) we have
a description of the twisted primitive elements in $\U$. We
now choose
$$
Y_s = q^{1/2}B-q^{-1/2}C + {{s^{-1}+s}\over{q^{-1}-q}}(A-D)
$$
in the space of twisted primitive elements. Then $Y_sA$
is a self-adjoint element in $\Unc$ for $s\in\R\backslash\{ 0\}$, 
or $s\in\T$, the unit circle. We study the bounded self-adjoint 
operator $\pi_k(Y_sA)$. In order to formulate the following result,
we recall the Al-Salam and Chihara polynomials defined in
Exercise~3.4. By $S_n$ we denote the orthonormal Al-Salam and
Chihara polynomials;
$$
S_n(x;a,b|q) = {1\over{\sqrt{(q,ab;q)_n}}}\, s_n(x;a,b|q)
$$
and corresponding normalised orthogonality measure
$dm(x;a,b|q) = dm(x;a,b,0,0|q)$, cf. Exercise~3.4 and
\thetag{\vglnormalisedAWmeasure}. We also use the
notation $\mu(x)=(x+x^{-1})/2=\mu(x^{-1})$ for $x\not= 0$
in this section.

\proclaim{Proposition \theoremname{\propspectralYsA}}
Define $\Lambda\colon \Hi \to
L^2(\R, dm(\cdot;q^{2k}s,q^{2k}/s|q^2))$
by
$$
\Lambda \colon e^k_n \mapsto S_n (\cdot;q^{2k}s,q^{2k}/s|q^2),
$$
then $\Lambda$ is a unitary mapping intertwining $\pi_k(Y_sA)$
with $2(M-\mu(s))/(q^{-1}-q)$, where $M$ is the multiplication
map on $L^2(\R, dm(\cdot;q^{2k}s,q^{2k}/s|q^2))$.
\endproclaim

\demo{Remark \theoremname{\rempropspectralYsA}}
Proposition \thmref{\propspectralYsA} says that formally
$$
v^k(x) = \sum_{n=0}^\infty S_n( \mu(x); q^{2k}s, q^{2k}/s|q^2)
\, e_n^k
$$
is an eigenvector of the bounded self-adjoint operator
$\pi_k(Y_sA)$ for the eigenvalue
$$
\la_x = {{ x+x^{-1} - s - s^{-1}}\over{q^{-1}-q}} =
2{{\mu(x)-\mu(s)}\over{q^{-1}-q}}.
$$
The spectrum of $\pi_k(Y_sA)$ is
$\{\la_x\mid \mu(x)\in\text{supp}
(dm(\cdot;q^{2k}s,q^{2k}/s|q^2))\}$.
\enddemo

\demo{Proof} We have
$$
\multline
\bigl( (q^{-1}-q)\pi_k(Y_sA) + s +s^{-1}\bigr)\, e_n =
\sqrt{ (1-q^{2n+2})(1-q^{4k+2n})}\, e_{n+1} \\ + q^{2k+2n}
(s+s^{-1})\, e_n + \sqrt{ (1-q^{2n})(1-q^{4k+2n-2})}\, e_{n-1},
\endmultline
$$
and comparing with the three-term recurrence relation for the 
orthonormal Al-Salam and Chihara polynomials $S_n(x)=S_n(x;a,b|q)$,
cf. Exercise~3.4,
$$
\aligned
2x\, S_n(x) &= a_{n+1}\, S_{n+1}(x) + q^n(a+b)\, S_n(x) + a_n\, 
S_{n-1}(x),\\ a_n &= \sqrt{ (1-abq^{n-1})(1-q^n)}.
\endaligned
\tag\eqname{\vgldrietermAlSC}
$$
leads to the result, cf. proof of Theorem \thmref{\thmFavard}.
\qed\enddemo

Next we consider the action of $Y_sA$ in the tensor product 
representation $\pi_{k_1}\otimes \pi_{k_2}$ acting in 
$\Hi\otimes\Hi$. This can be done using orthogonal polynomials in 
two variables.

\proclaim{Proposition \theoremname{\propspectralYsAintensor}}
Define $\Upsilon\colon \Hi\otimes\Hi
\to L^2(\R^2, dm(x,y))$,
where
$$
dm(x,y) = dm(x;q^{2k_1}w,q^{2k_1}/w|q^2)\, dm(y;q^{2k_2}s,
q^{2k_2}/s|q^2), \qquad y=\mu(w),
$$
by
$$
\Upsilon \colon
e^{k_1}_{n_1} \otimes e^{k_2}_{n_2}
\mapsto S_{n_1} (x;q^{2k_1}w,q^{2k_1}/w|q^2)\,
S_{n_2}(y;q^{2k_2}s,q^{2k_2}/s|q^2)
$$
then $\Upsilon$ is a unitary mapping intertwining
$\pi_{k_1}\otimes\pi_{k_2}(\De(Y_sA))$
with $2(M_x-\mu(s))/(q^{-1}-q)$, where $M_x$ is multiplication by 
$x$ in $L^2(dm(x,y))$.
\endproclaim

\demo{Remark \theoremname{\rempropspectralYsAintensor}}(i) Put, 
$y=\mu(w)$,
$$
R_{l,m}(x,y) = S_l (x;q^{2k_1}w,q^{2k_1}/w|q^2)\,
S_m(y;q^{2k_2}s,q^{2k_2}/s|q^2)
$$
then $R_{l,m}$ are orthonormal polynomials in two variables
of degree $l$ in $x$ and $l+m$ in $y$;
$$
\iint R_{l,m}(x,y)\, R_{r,s}(x,y)\, dm(x,y) =
\de_{lr}\de_{ms},
$$
as a straightforward consequence of the orthogonality relations of 
the Al-Salam and Chihara polynomials.

(ii) Proposition \thmref{\propspectralYsAintensor} states that 
formally the vector $w(x;y) =$
$$
\multline
\sum_{n_1,n_2=0}^\infty
S_{n_1} (\mu(x);q^{2k_1}y,q^{2k_1}/y|q^2)\,
S_{n_2}(\mu(y);q^{2k_2}s,q^{2k_2}/s|q^2)\,
e^{k_1}_{n_1} \otimes e^{k_2}_{n_2} \\ =
\sum_{n_1=0}^\infty S_{n_1} (\mu(x);q^{2k_1}y,q^{2k_1}/y|q^2)\,
e^{k_1}_{n_1} \otimes v^{k_2}(y)
\endmultline
$$
is an eigenvector of $\pi_{k_1}\otimes\pi_{k_2}(\De(Y_sA))$
for the eigenvalue $\la_x$. This last observation is essentially
the way to obtain Proposition \thmref{\propspectralYsAintensor}, 
since $\De(Y_sA)=A^2\otimes Y_sA+Y_sA\otimes 1$ acts as a 
three-term recurrence operator in 
$e^{k_1}_{n_1} \otimes v^{k_2}(y)$.
\enddemo

\demo{Proof} We use $\De(Y_sA)=A^2\otimes Y_sA+Y_sA\otimes 1$
and Proposition \thmref{\propspectralYsA} to
define for fixed $y$ the map $\Lambda\colon\Hi\otimes\Hi\to \Hi$
by
$$
\Lambda \colon e^{k_1}_{n_1} \otimes e^{k_2}_{n_2}
\mapsto
S_{n_2}(y;q^{2k_2}s,q^{2k_2}/s|q^2)\, e^{k_1}_{n_1}
$$
to obtain the recurrence in $n_1$
$$
\multline
\Lambda \Bigl( (q^{-1}-q)\bigl(\pi_{k_1}\otimes\pi_{k_2}
(\Delta(Y_sA)\bigr)
 + s +s^{-1}\Bigr)\,e^{k_1}_{n_1} \otimes e^{k_2}_{n_2}  = \\
S_{n_2}(y;q^{2k_2}s,q^{2k_2}/s|q^2)
\Bigl( q^{2n_1} \bigl( (s+s^{-1})q^{2k_1} + \la_y q^{2k_1}
(q^{-1}-q)\bigr)\, e^{k_1}_{n_1} \\
 +\sqrt{ (1-q^{2n_1+2})(1-q^{4k_1+2n})}\, e_{n_1+1}^{k_1}
+ \sqrt{ (1-q^{2n_1})(1-q^{4k_1+2n_1-1})}\, e_{n_1-1}^{k_1}\Bigr)
\endmultline
$$
Use the explicit expression for $\la_y$ as in Remark
\thmref{\rempropspectralYsA} and the three-term recurrence
relation \thetag{\vgldrietermAlSC} to obtain the result.
\qed\enddemo

By Lemma \thmref{\lemCGCforUnc} the representation space 
$\Hi\otimes\Hi$ has another orthonormal
basis. We now calculate the action of $\Upsilon\, e^k_n$.

\proclaim{Proposition \theoremname{\propactioUpsionlonekn}}
Let $k=k_1+k_2+j$ for $j\in\Zp$, and $x=\mu(z)$ then
$$
\align
\bigl( \Upsilon\, e^k_n\bigr)(x,y) &= S_n(x;q^{2k}s,q^{2k}/s|q^2)
\, \bigl( \Upsilon\, e^k_0\bigr)(x,y), \\
\bigl( \Upsilon\, e^k_0\bigr)(x,y) &= C\,
p_j(y;q^{2k_1}z, q^{2k_1}/z,q^{2k_2}s,q^{2k_2}/s|q^2),\\
C^{-1} &= \sqrt{ (q^2,q^{4k_1},q^{4k_2},q^{4k_1+4k_2+2j-2};q^2)_j}.
\endalign
$$
\endproclaim

\demo{Proof} By Proposition \thmref{\propspectralYsAintensor}
$\Upsilon$ intertwines $\pi_{k_1}\otimes\pi_{k_2}(\De(Y_sA))$
with $2(M_x-\mu(s))/(q^{-1}-q)$, and by Lemmas
\thmref{\lemmadecomptensorsunc} and \thmref{\lemCGCforUnc}
the action of $\pi_{k_1}\otimes\pi_{k_2}(\De(Y_sA))$ on
$e_n^k$ is $\pi_k(Y_sA)$, so
$$
2{{x-\mu(s)}\over{q^{-1}-q}} \, \Upsilon e_n^k(x,y) =
\Upsilon (\pi_k(Y_sA)\, e_n^k) (x,y)
$$
and we obtain the three-term recurrence relation as
in Proposition \thmref{\propspectralYsA}, but with
initial conditions $\Upsilon e^k_{-1} = 0$ and
$\Upsilon e^k_0(x,y)$ some polynomial in two variables
$x$ and $y$. Hence, the first statement follows.

Now $\Upsilon e^k_n(x,y)$ is a polynomial in two variables,
and since $\Upsilon$ is unitary we have the orthogonality relations
$\de_{nm}\de_{kl} = \langle \Upsilon e^k_n,\Upsilon e^l_m\rangle
=$
$$
\int S_n(x;q^{2k}s,q^{2k}/s|q^2)S_m(x;q^{2l}s,q^{2l}/s|q^2)
\int \Upsilon e^k_0(x,y) \Upsilon e^l_0(x,y)\, dm(x,y),
$$
by our first observation. Since the orthogonality measure
for the Al-Salam and Chihara polynomials is unique (i.e.
the corresponding moment problem is determined), we find
that for $k=l$ the inner
integral must give the orthogonality measure
$dm(x;q^{2k}s,q^{2k}/s|q^2)$
as measures on $\R$ with respect to $x$. Take $k\not= l$, then
the inner integral as function of $x$ gives zero when integrated
against an arbitrary polynomial. Since the support of the measure
is compact, the polynomials are dense in the corresponding
$L^2$-space and the inner integral is zero for $k\not= l$.

Since
$\Upsilon e^{k_1}_{n_1}\otimes e^{k_2}_{n_2}(x,y)$ is a polynomial
of degree $n_1+n_2$ in $y$, see Remark 
\thmref{\rempropspectralYsAintensor},
the Clebsch-Gordan decomposition of Lemma \thmref{\lemCGCforUnc}
implies that $\Upsilon e^k_0(x,y)$ is a polynomial of degree $j$ 
in $y$, where $k=k_1+k_2+j$ for $j\in\Zp$, say 
$\Upsilon e^k_0(x,y)=p_j(y)$.
Let $l=k_1+k_2+i$ for $i\in\Zp$, then we obtain the orthogonality
relations
$$
\int_y p_j(y)p_i(y)\, dm(x,y) = \de_{ij}
dm(x;q^{2k}s,q^{2k}/s|q^2), \quad\text{for almost all $x$.}
$$

We now assume for ease of presentation
that $dm(x,y)$ is absolutely continuous, the general case being
proved similarly. This is the case
if we take $q^{2k_2}<|s|<q^{-2k_2}$ since $k_1,k_2>0$.
Put $x=\cos\th$, $y=\cos\psi$, and use Theorem 
\thmref{\thmAWpolsareorthogonal}
and \thetag{\vglnormalisedAWmeasure}, to find
$$
\multline
{1\over{2\pi}}\int_0^\pi p_i(\cos\psi)p_j(\cos\psi)
{{ (e^{\pm 2i\psi}, e^{\pm 2i\th};q^2)_\infty}\over
{(q^{2k_2}se^{\pm i\psi}, q^{2k_2}e^{\pm i\psi}/s,
q^{2k_1}e^{\pm i\psi \pm i\th};q^2)_\infty}} d\psi
= \\ \de_{ij}
{{ (q^{4k_1+4k_2+4j};q^2)_\infty}\over{
(q^2,q^{4k_1},q^{4k_2};q^2)_\infty}}
{{ (e^{\pm 2i\th};q^2)_\infty}\over{
(q^{2k_1+2k_2+2j}se^{\pm i\th}, 
q^{2k_1+2k_2+2j}e^{\pm i\th}/s;q^2)_\infty}}
\endmultline
$$
for almost all $\th$. The $\pm$-signs mean that we take all
possible combinations in the infinite $q$-shifted factorials.
Cancelling the $(e^{\pm 2i\th};q^2)_\infty$ on both sides and 
comparing the result with Theorem \thmref{\thmAWpolsareorthogonal} 
we see that $p_j$ is a multiple of
$p_j(\cdot;q^{2k_1}e^{i\th}, q^{2k_1}e^{-i\th},q^{2k_2}s,
q^{2k_2}/s|q^2)$.
The (real) constant in front follows up to a sign by comparing the
squared norms. To determine the sign we recall the
normalisation of Lemma \thmref{\lemCGCforUnc},
$$
\multline
0<\langle e^k_0,e^{k_1}_0\otimes e^{k_2}_j\rangle =
\langle \Upsilon e^k_0, \Upsilon e^{k_1}_0\otimes e^{k_2}_j\rangle 
= \\ C \iint  p_j(y;q^{2k_1}z, q^{2k_1}/z,q^{2k_2}s,q^{2k_2}/s|q^2)
S_j(y;q^{2k_2}s,q^{2k_2}/s|q^2)\, dm(x,y).
\endmultline
$$
The Askey-Wilson polynomials are orthogonal with respect to the 
integration over $y$, and the Al-Salam and Chihara polynomial has 
positive leading coefficient, so the inner integral is positive and
the remaining measure over $x$ is positive as well. Hence the 
double integral is positive and we conclude that $C>0$.
\qed\enddemo

\demo{Remark \theoremname{\rempropactioUpsionloneknone}}(i) The 
proof shows that the polynomials, $x=\mu(z)$,
$$
P_{l,m}(x,y) = S_l(x;q^{2k_1+2k_2+2m}s,q^{2k_1+2k_2+2m}/s|q^2)\,
p_m(y;q^{2k_1}z, q^{2k_1}/z,q^{2k_2}s,q^{2k_2}/s|q^2),
$$
of degree $m$ in $y$ and $l+m$ in $x$ are orthogonal with respect
to the measure $dm(x,y)$;
$$
\iint P_{l,m}(x,y)P_{r,s}(x,y)\, dm(x,y) =
\de_{lr}\de_{ms}
(q^2,q^{4k_1},q^{4k_2},q^{4k_1+4k_2+2m-2};q^2)_m.
$$

(ii) $\Upsilon e^k_0$ can also be calculated explicitly using the
Clebsch-Gordan coefficients, so, $y=\mu(w)$,
$$
\Upsilon e^k_0(x,y) = \sum_{n_1+n_2=j} C^{k_1,k_2,k}_{n_1,n_2,0}\,
S_{n_1} (x;q^{2k_1}w,q^{2k_1}/w|q^2)\,
S_{n_2}(y;q^{2k_2}s,q^{2k_2}/s|q^2)
$$
by Lemma \thmref{\lemCGCforUnc} and
Proposition \thmref{\propspectralYsAintensor}. This can be 
evaluated directly by using the ${}_3\vp_2$-series representation 
for the Al-Salam and Chihara polynomials and the explicit 
expression for $C^{k_1,k_2,k}_{n_1,n_2,0}$ derived in the proof of 
Lemma \thmref{\lemCGCforUnc}, interchanging summations and using
summation formulas of \S 3. The sum can also be evaluated by
viewing it as a convolution, and using suitable generating
functions for the Al-Salam and Chihara and Askey-Wilson polynomials.
\enddemo

%%%%%%%%%%%%%%%%%%%%%%%%%%%%%%%%%%%%%%%%%%%%%%%%%%%%%%%%%%%%%%%%%%%%
%%N E W   S U B S E C T I O N%%%%%%%%%%%%%%%%%%%%%%%%%%%%%%%%%%%%%%%
%%%%%%%%%%%%%%%%%%%%%%%%%%%%%%%%%%%%%%%%%%%%%%%%%%%%%%%%%%%%%%%%%%%%
\subhead\newsubsection 
Convolution theorem for Al-Salam and Chihara polynomials
\endsubhead
The convolution formula for the Al-Salam and Chihara polynomials
is obtained by applying $\Upsilon$ to
Lemma \thmref{\lemCGCforUnc} using the results of Propositions
\thmref{\propspectralYsAintensor} and 
\thmref{\propactioUpsionlonekn}.
The results holds as an identity in a weighted $L^2$-space, but 
since it is a polynomial identity it holds for all $x$, $y$.

\proclaim{Lemma \theoremname{\lemconvtheoremstepone}}
For $x=\mu(z)$, $y=\mu(w)$, and $k=k_1+k_2+j$ we have
$$
\multline
\sum_{n_1+n_2=n+j} C^{k_1,k_2,k}_{n_1,n_2,n}
\, S_{n_1} (x;q^{2k_1}w,q^{2k_1}/w|q^2)\,
S_{n_2}(y;q^{2k_2}s,q^{2k_2}/s|q^2) = \\
{{S_n(x;q^{2k}s,q^{2k}/s|q^2)
p_j(y;q^{2k_1}z, q^{2k_1}/z,q^{2k_2}s,q^{2k_2}/s|q^2)}\over{
 \sqrt{ (q^2,q^{4k_1},q^{4k_2},q^{4k_1+4k_2+2j-2};q^2)_j}}}.
\endmultline
$$
\endproclaim

We have not yet calculated the Clebsch-Gordan
coefficients explicitly, but we can now use
Lemma \thmref{\lemconvtheoremstepone} to
determine $C^{k_1,k_2,k}_{n_1,n_2,n}$ by specialising to
a generating function for the Clebsch-Gordan coefficients.
The result is phrased in terms of $q$-Hahn polynomials, which are
defined as follows;
$$
Q_n(q^{-x},a,b,N;q) = {}_3\vp_2 \left( {{q^{-n}, 
q^{-x},abq^{n+1}}\atop{ aq,\ q^{-N}}};q,q\right).
$$

\proclaim{Lemma \theoremname{\lemexplicitCGCgenfunct}}
With $n_1+n_2=n+j$ we get
$$
C^{k_1,k_2,k_1+k_2+j}_{n_1,n_2,n} =
C\ Q_j(q^{-2n_1};q^{4k_1-2},q^{4k_2-2},n+j;q^2),
$$
with the constant $C$ given by
$$
{{ q^{2n_1k_1-2n(k+j)}
(q^2;q^2)_{n+j}
\sqrt{(q^{4k_1};q^2)_{n_1} (q^{4k_2};q^2)_{n_2} 
(q^{4k_1};q^2)_j}}\over
{\sqrt{ (q^2;q^2)_n(q^2;q^2)_{n_1}(q^2;q^2)_{n_2}(q^2;q^2)_j
(q^{4k_1+4k_2+4j};q^2)_n (q^{4k_2};q^2)_j 
(q^{4k_1+4k_2+2j-2};q^2)_j}}} .
$$
\endproclaim

\demo{Proof} Observe that $C^{k_1,k_2,k}_{n_1,n_2,n}$ is
independent of $s$, $x=\mu(z)$ and $y=\mu(w)$.
Specialise $w=q^{2k_2}s$ and $z=q^{2k_1}/w=q^{2k_1-2k_2}/s$, then
the Al-Salam and Chihara polynomials in the summand on the
left hand side of Lemma \thmref{\lemconvtheoremstepone}
can be evaluated explicitly, since the ${}_3\vp_2$-series
reduces to $1$. For this choice the Askey-Wilson polynomial
on the right hand side can also be evaluated explicitly,
and we obtain the generating function for the
Clebsch-Gordan coefficients
$$
\multline
\sum_{n_1+n_2=n+j} C^{k_1,k_2,k}_{n_1,n_2,n}
\, q^{2n_1(k_2-k_1)-2n_2k_2} s^{n_1-n_2}
{{ \sqrt{ (q^{4k_1};q^2)_{n_1}(q^{4k_2};q^2)_{n_2}}}\over
 { \sqrt{ (q^2;q^2)_{n_1}(q^2;q^2)_{n_2}}}} = \\
{{q^{-2jk_2 -2n(k_1+k_2+j)} s^{n-j}
(q^{4k_1},q^{4k_2}, q^{4k_2}s^2;q^2)_j }\over{
 \sqrt{ (q^2,q^{4k_1},q^{4k_2},q^{4k_1+4k_2+2j-2};q^2)_j}}}
\sqrt{ {{ (q^{4k_1+4k_2+4j};q^2)_n}\over{(q^2;q^2)_n}}}\\ \times
\ {}_3\vp_2 \left( {{q^{-2n},q^{4k_2+2j}, q^{4k_1+2j}/s^2}\atop{
q^{4k_1+4k_2+4j},\ 0}};q^2,q^2\right).
\endmultline
$$
This determines $C^{k_1,k_2,k}_{n_1,n_2,n}$, but it takes some
work to find the expression in terms of $q$-Hahn polynomials.
First, take $n_1$ as the summation parameter in the sum and
multiply both sides by $s^{n+j}$ to find that both sides are
polynomials of degree $n+j$ in $s^2$. Apply
$$
{}_2\vp_1(q^{-n},b;c;q,z) = {{(c/b;q)_n}\over{(c;q)_n}} (bz/q)^n
\ {}_3\vp_2 \left( {{q^{-n}, c^{-1}q^{1-n}, q/z}\atop
{bq^{1-n}/c,\ 0}};q,q\right),
$$
which is just Exercise~3.6 with $a=q^{-n}$ and the series in the
terminating ${}_2\vp_2$-series inverted, and the $q$-binomial
theorem to $(q^{4k_2}s^2;q^2)_j$ to find
$$
\align
&s^{n+j} (q^{4k_2}s^2;q^2)_j
\ {}_3\vp_2 \left( {{q^{-2n},q^{4k_2+2j}, q^{4k_1+2j}/s^2}\atop{
q^{4k_1+4k_2+4j},\ 0}};q^2,q^2\right) = \\
&{{(q^{2-2n-4k_2-2j};q^2)_n}\over{(q^{2-2n-4k_1-4k_2-4j};q^2)_n}}
\ {}_1\vp_0\left( {{q^{-2j}}\atop{-}};q^2,s^2 q^{2j+4k_2}\right)
\ {}_2\vp_1\left( {{q^{-2n},q^{4k_1+2j}}\atop{q^{2-2n-4k_2-2j}}};
q^2, s^2 q^{2-4k_1-2j}\right).
\endalign
$$
The product of the $q$-hypergeometric series is written as a
polynomial in $s^2$ by
$$
\align
&\sum_{n_1=0}^{n+j} s^{2n_1}
\sum_r {{(q^{-2j};q^2)_{n_1-r} (q^{-2n},q^{4k_1+2j};q^2)_r}\over
{(q^2;q^2)_{n_1-r} (q^2,q^{2-2n-4k_2-2j};q^2)_r}}
q^{2(n_1-r)(j+2k_2) + 2r(1-2k_1-j)} = \\
&\sum_{n_1=0}^{n+j} s^{2n_1}
{{(q^{-2j};q^2)_{n_1}}\over{(q^2;q^2)_{n_1}}} q^{n_1(2j+4k_2)}
\ {}_3\vp_2 \left( {{q^{-2n},q^{4k_1+2j},q^{-2n_1}}\atop
{q^{2-2n-4k_2-2j}, q^{2+2j-2n_1}}};q^2, q^{2-4k_1-4k_2-2j}\right).
\endalign
$$

This gives an explicit expression for the Clebsch-Gordan 
coefficients in terms of a terminating ${}_3\vp_2$-series. To put 
it into the required form in terms of $q$-Hahn polynomials, we 
need to apply some transformations for ${}_3\vp_2$-series, namely
\cite{\GaspR, (III.13), (III.11)}. Then the ${}_3\vp_2$-series
can be rewritten as
$$
{{(q^{2-2n_1-4k_1};q^2)_{n_1}}\over{(q^{2+2j-2n_1};q^2)_{n_1}}}
\, {}_3\vp_2\left( {{q^{-2n_1}, q^{4k_1+2j}, q^{2-2j-4k_2}}\atop{
q^{2-2j-4k_2-2n}, q^{4k_1}}};q^2,q^2\right)
$$
by \cite{\GaspR, (III.13)}. And by \cite{\GaspR, (III.11)}
this ${}_3\vp_2$-series can be written as
$$
{{(q^{-2j-2n};q^2)_{n_1}}\over{(q^{2-2j-4k_2-2n};q^2)_{n_1}}} 
q^{2(1-2k_2)n_1}
\, {}_3\vp_2 \left( {{q^{-2n_1},q^{-2j}, q^{4k_1+4k_2+2j-2}}\atop
{q^{4k_1},\ q^{-2n-2j}}};q^2,q^2\right),
$$
which is of the desired form. The constant follows by
a straightforward calculation.
\qed\enddemo

Now all ingredients for the general convolution theorem for
the Al-Salam and Chihara polynomials are known. Applying
Lemma \thmref{\lemexplicitCGCgenfunct} in Lemma
\thmref{\lemconvtheoremstepone}, and rewriting the result
proves the following theorem.

\proclaim{Theorem \theoremname{\thmgenconvoforAlSCpols}}
With $x=\mu(z)$ and $y=\mu(w)$ and $n,j\in\Zp$, $k_1,k_2>0$ we have
$$
\multline
\sum_{l=0}^{n+j}
q^{2lk_1}
\left[ {{n+j}\atop l}\right]_{q^2}\
Q_j(q^{-2l};q^{4k_1-2}, q^{4k_2-2},n+j;q^2)\\
\times
s_l(x;q^{2k_1}w,q^{2k_1}/w|q^2)\,
s_{n+j-l}(y;q^{2k_2}s,q^{2k_2}/s|q^2) = \\
{{q^{2n(k_1+k_2+2j)}}\over{(q^{4k_1};q^2)_j}}
\, s_n(x;q^{2k_1+2k_2+2j}s,q^{2k_1+2k_2+2j}/s|q^2)\,
p_j(y;q^{2k_1}z,q^{2k_1}/z,q^{2k_2}s,q^{2k_2}/s|q^2).
\endmultline
$$
\endproclaim

\demo{Remark \theoremname{\remthmgenconvoforAlSCpols}}(i)
Theorem \thmref{\thmgenconvoforAlSCpols} is a connection
coefficient formula for orthogonal polynomials in two
variables, orthogonal for the same measure, cf.
Remarks \thmref{\rempropspectralYsAintensor}(i)
and \thmref{\rempropactioUpsionloneknone}(i). The
connection coefficients being given by the
$q$-Hahn polynomials. Since the Clebsch-Gordan coefficients
form a unitary matrix, we also have
$e_{n_1}^{k_1}\otimes e_{n_2}^{k_2} = \sum_{n,k}
 C^{k_1,k_2,k}_{n_1,n_2, n}\, e_n^k$, and from this we can obtain
the inverse connection coefficient problem. This also follows
from the orthogonality relations for the dual $q$-Hahn polynomials,
cf. Exercise~7.4.

(ii) The case $j=0$ gives a simple convolution
property for the Al-Salam and Chihara polynomials, since
the $q$-Hahn and the Askey-Wilson polynomial reduce to $1$.
The case $n=0$ is also of interest, since then the
$q$-Hahn polynomial can be evaluated and the Al-Salam and
Chihara polynomial on the right hand side reduces to $1$.
In both cases we have a free parameter in the sum.

(iii) Formally, in the representation space $\Hi\otimes\Hi$ we have
two bases of (generalised) eigenvectors for the action of $Y_sA$,
namely $v^k(x)$ as in Remark \thmref{\rempropspectralYsA} and
$w(x;y)$ as in Remark \thmref{\rempropspectralYsAintensor}(ii).
They are connected by Clebsch-Gordan coefficients, which are
now expressible as Askey-Wilson polynomials;
$$
w(x;y) = \sum_{j=0}^\infty {{
p_j(\mu(y);q^{2k_1}x, q^{2k_1}/x,q^{2k_2}s,q^{2k_2}/s|q^2)}\over{
 \sqrt{ (q^2,q^{4k_1},q^{4k_2},q^{4k_1+4k_2+2j-2};q^2)_j}}}
\, v^{k_1+k_2+j}(x).
$$
This follows immediately from Lemma \thmref{\lemconvtheoremstepone}
and the orthogonality relations for the Clebsch-Gordan coefficients,
cf. Remark \thmref{\remthmgenconvoforAlSCpols}(i), using the 
explicit expressions for $v^k(x)$ and $w(x;y)$ as in Remarks
\thmref{\rempropspectralYsA} and 
\thmref{\rempropspectralYsAintensor}(ii).

The orthogonality
relations for the Askey-Wilson polynomials then imply
$$
\align
v^{k_1+k_2+j}(x) &= {1\over{h_j}} \int w(x;y)
{{p_j(\mu(y);q^{2k_1}x, q^{2k_1}/x,q^{2k_2}s,q^{2k_2}/s|q^2)}\over{
 \sqrt{ (q^2,q^{4k_1},q^{4k_2},q^{4k_1+4k_2+2j-2};q^2)_j}}}
\\ & \qquad\qquad \quad \times
dm(\mu(y);q^{2k_1}x, q^{2k_1}/x,q^{2k_2}s,q^{2k_2}/s|q^2), \\
h_j &= {{1-q^{2j-2+4k_1+4k_2}}\over{1-q^{4j-2+4k_1+4k_2}}}
{{(q^2,q^{4k_1}, q^{4k_2}, q^{2k_1+2k_2}x^{\pm 1}s^{\pm 1};q^2)_j}
\over{(q^{4k_1+4k_2};q^2)_j}},
\endalign
$$
where we have to take all possible choices of signs at the $\pm$.
This can also be proved using Lemma \thmref{\lemconvtheoremstepone}.
\enddemo

%%%%%%%%%%%%%%%%%%%%%%%%%%%%%%%%%%%%%%%%%%%%%%%%%%%%%%%%%%%%%%%%%%%%
%%N E W   S U B S E C T I O N%%%%%%%%%%%%%%%%%%%%%%%%%%%%%%%%%%%%%%%
%%%%%%%%%%%%%%%%%%%%%%%%%%%%%%%%%%%%%%%%%%%%%%%%%%%%%%%%%%%%%%%%%%%%
%NOTES AND REFERENCES%%%%%%%%%%%%%%%%%%%%%%%%%%%%%%%%%%%%%%%%%%%%%%%
\subhead Notes and references
\endsubhead
The results are obtained in
joint work with Van der Jeugt \cite{\KoelVdJ} elaborating
an idea of Granovskii and Zhedanov \cite{\GranZ}.
For the case of the
Lie algebra ${\frak{su}}(1,1)$ see \cite{\KoelVdJ} and
Van der Jeugt \cite{\VdJeug}. A similar approach can be used in the
case $\Usu$ to lead to a formula for $q$-Krawtchouk polynomials,
see \cite{\KoelVdJ} for details. The method can also be extended
to three-fold tensor products, and then we use the $q$-Racah 
coefficients, see \cite{\KoelVdJ}. Some related results on 
overlap coefficients for quantum algebras can be found in Klimyk
 and Kachurik \cite{\KlimK}.

Lemma \thmref{\lemmadecomptensorsunc} is well-known, see e.g.
Kalnins, Manocha and Miller \cite{\KalnMMJMP}, where also the
explicit expression for the Clebsch-Gordan coefficients is
derived. Exercises~7.1-2 are also taken from \cite{\KalnMMJMP}.
The generating function for the Askey-Wilson polynomials,
see Remark \thmref{\rempropactioUpsionloneknone}(ii), is given
by Ismail and Wilson \cite{\IsmaW}. Taking $c=d=0$ also gives a
generating function for the Al-Salam and Chihara polynomials,
the other generating function needed is
$$
\sum_{n=0}^\infty {{q^{n(n-1)/2} t^n}\over {(q,ab;q)_n}}\,
s_n(\mu(x);a,b\mid q) = (-t/a;q)_\infty {}_2\vp_1\left(
{{ax,a/x}\atop{ab}};q, -t/a\right).
$$

The case $j=0$ of Theorem \thmref{\thmgenconvoforAlSCpols}
was the reason for Al-Salam and Chihara \cite{\AlSaC} to
introduce the Al-Salam and Chihara polynomials as the
most general set of orthogonal polynomials still
satisfying a convolution property, see also
Al-Salam \cite{\AlSa, \S 8}.
Al-Salam and Chihara
obtained the three-term recurrence relation, and later
Askey and Ismail \cite{\AskeIMAMS} determined the
orthogonality measure.

The complete representation theory of the quantised
universal enveloping algebra $\Unc$ can be found in
Burban and Klimyk \cite{\BurbK}. The representation theory
is similar to that of ${\frak{su}}(1,1)$, but there is an
extra series of representations, the so-called
strange series. The dual Hopf $\ast$-algebra $A_q(SU(1,1))$
is defined in Theorem \thmref{\thmstarstructureonA}, and
Masuda et al. \cite{\MasuMNNSU} give explicit expressions
for the matrix elements in terms of ${}_2\vp_1$-series
of argument $-q^{-1}\be\ga$, which are $q$-analogues of
the Jacobi functions. There are no corresponding
orthogonality relations for the spherical elements, but the
result has to be stated in terms of transform pairs. The transform
pair is a $q$-analogue of
the Mehler-Fock transform, see Vaksman and
Korogodski\u\i\ \cite{\VaksKFAA} and Kakehi, Masuda
and Ueno \cite{\KakeMU}, Kakehi \cite{\Kake} for
analytic proofs.

%%%%%%%%%%%%%%%%%%%%%%%%%%%%%%%%%%%%%%%%%%%%%%%%%%%%%%%%%%%%%%%%%%%%
%%N E W   S U B S E C T I O N%%%%%%%%%%%%%%%%%%%%%%%%%%%%%%%%%%%%%%%
%%%%%%%%%%%%%%%%%%%%%%%%%%%%%%%%%%%%%%%%%%%%%%%%%%%%%%%%%%%%%%%%%%%%
%EXERCISES%%%%%%%%%%%%%%%%%%%%%%%%%%%%%%%%%%%%%%%%%%%%%%%%%%%%%%%%%%
\subhead Exercises
\endsubhead

\item{\the\sectionno.1} An explicit model for the positive discrete
series $\pi_k$ can be obtained as follows. First prove that the
operators $A=q^kT_{q^2}$,
$$
B= {M\over{q-q^{-1}}}\bigl( q^{-2k}T_q^{-1} - q^{2k}T_q\bigr),
\qquad
C = T_q^{-1} D_{q^2},
$$
acting on the formal power series give a representation of
$\U$. Here $M$ is the multiplication operator, $Mf(z)=zf(z)$,
$T_a$, $a\not= 0$, is the shift operator, $T_af(z) = f(az)$,
and $D_q$ is the $q$-derivative,
$D_qf(z) = (f(z)-f(qz))/((1-q)z)$, see
Exercise~3.5. Then show that taking $f_n(z)=z^n$ as an orthogonal
basis gives a unitary representation of $\Unc$. Prove that
$$
\parallel f_n \parallel^2 = q^{n(2k-1)}
{{(q^2;q^2)_n}\over{(q^{4k};q^2)_n}}
\parallel f_0 \parallel^2 .
$$
So the representation space is the space of formal power
series $\sum_{n=0}^\infty c_n z^n$ such that
$\sum_{n=0}^\infty |c_n|^2 \parallel f_n\parallel^2 <\infty$,
which are the functions analytic on the disc with radius
$q^{1/2-k}$.

\item{\the\sectionno.2} By Exercise~7.1
$f_{n_1,n_2}(z,w) = z^{n_1}w^{n_2}$ forms an orthogonal basis
for the representation space of the tensor product of
the model for $\pi_{k_1}\otimes\pi_{k_2}$. Show that all
eigenvectors of $\De(A)$ in the kernel of $\De(C)$ are
given by
$$
p_{j,0}(z,w) =  z^j (q^{2-k_1-k_2-3s}w/z;q^2)_j, \qquad s\in\Zp
$$
And $\De(A) p_{j,0} = q^{k_1+k_2+j} p_{j,0}$. Now define
inductively
$$
p_{j,n} = {{q-q^{-1}}\over{(q^{-2k_1-2k_2-2j-n}- 
q^{2k_1+2k_2+2j+n})}} \De(B) p_{j,n-1}
$$
and prove that the span of $p_{j,n}$, $n\in\Zp$, realises an 
irreducible unitary representation of $\Unc$ equivalent to
the discrete series representation $\pi_{k_1+k_2+j}$.
Finish the proof of Lemma \thmref{\lemmadecomptensorsunc}
by proving that $p_{n,j}$, $n,j\in\Zp$ form an
orthogonal basis.

\item{\the\sectionno.3} Derive the recurrence relation
mentioned in Remark \thmref{\remlemCGCforUnc} and use
Lemma \thmref{\lemexplicitCGCgenfunct} to obtain a
contiguous relation for ${}_3\vp_2$-series.

\item{\the\sectionno.4} Use the unitarity of the
Clebsch-Gordan coefficients, i.e.
$$
\sum_{n_1=0}^{n+j} C^{k_1,k_2,k_1+k_2+j}_{n_1,n+j-n_1,n}
C^{k_1,k_2,k_1+k_2+i}_{n_1,n+j-n_1,n+j-i} = \de_{ij}
$$
with Lemma \thmref{\lemexplicitCGCgenfunct} to
derive the orthogonality relations for
the $q$-Hahn polynomials $Q_n(q^{-x}) =
Q_n(q^{-x};a,b,N;q)$; for $n,m\in\{0,1,\ldots,N\}$,
$$
\sum_{x=0}^N Q_m(q^{-x})Q_n(q^{-x})
{{(aq,q)_x(bq;q)_{N-x}}\over{(q;q)_x(q;q)_{N-x}}} (aq)^{-x}
=\de_{nm} h_n.
$$
What are the dual orthogonality relations?

%%%%%%%%%%%%%%%%%%%%%%%%%%%%%%%%%%%%%%%%%%%%%%%%%%%%%%%%%%%%%%%%%%%%
%%N E W   S E C T I O N%%%%%%%%%%%%%%%%%%%%%%%%%%%%%%%%%%%%%%%%%%%%%
%%%%%%%%%%%%%%%%%%%%%%%%%%%%%%%%%%%%%%%%%%%%%%%%%%%%%%%%%%%%%%%%%%%%
\newpage

\head\newsection More examples\endhead

In the last section we shortly discuss some more known examples
between quantum groups and $q$-special functions, which are
similarly using the duality of Hopf $\ast$-algebras. We also
consider the quantum algebra approach, by considering it
for a special case. Finally, we give some open problems.
The references are given separately for each subsection, and
there more details can be found.

%%%%%%%%%%%%%%%%%%%%%%%%%%%%%%%%%%%%%%%%%%%%%%%%%%%%%%%%%%%%%%%%%%%%
%%N E W   S U B S E C T I O N%%%%%%%%%%%%%%%%%%%%%%%%%%%%%%%%%%%%%%%
%%%%%%%%%%%%%%%%%%%%%%%%%%%%%%%%%%%%%%%%%%%%%%%%%%%%%%%%%%%%%%%%%%%%
\subhead\newsubsection 
The quantum group of plane motions and $q$-Bessel functions
\endsubhead
In Exercise~2.4 the Hopf $\ast$-algebra $U_q({\frak m}(2))$
is obtained by a contraction procedure from $\Usu$. We use the
same letters $A$, $B$, $C$ and $D$ for the generators, then
the comultiplication, counit, antipode and $\ast$-operator
are unchanged, i.e. given by \thetag{\vglHopfstructureUq}
and Theorem \thmref{\thmstarstructureonU}(ii). The relations
among the generators change;
$$
AB=qBA,\quad AC=q^{-1}CA,\quad AD=1=DA,\quad BC=CB.
$$
We transpose the contraction procedure to the dual Hopf 
$\ast$-algebra
by demanding that the duality on the level of generators,
cf. \thetag{\vglabgdonbasis}, is not changed. Denoting the
resulting generators of the Hopf $\ast$-algebra $A_q(M(2))$ by
the same letters $\al$, $\be$, $\ga$ and $\de$, we see that
the action of the counit, antipode and $\ast$-operator
remains, cf. \thetag{\vglSenepsilonopAq} and Theorem
\thmref{\thmstarstructureonA}(ii). The commutation relations
and the action of the comultiplication change;
$$
\gather
\al\be =q\be\al ,\quad \al\ga = q\ga\al ,\quad \be\de = q\de\be ,
\quad \ga\de = q\de\ga ,\quad
\be\ga =\ga\be ,\quad \al\de = \de\al  =1 .
\endgather
$$
Then $A_q(M(2))$ is Hopf algebra.
The comultiplication $\Delta$, the counit $\varepsilon$ and the 
antipode $S$ given on the generators by
$$
\gather
\De(\al)=\al\otimes\al,\quad
\De(\be )=\al\otimes\be + \be\otimes\de ,\\
\De(\ga )=\ga\otimes\al + \de\otimes\ga ,\quad
\De(\de )=\de\otimes\de,\\
\ep(\al)=\ep(\de)=1, \qquad \ep(\be)=\ep(\ga)=0, \\
S(\al)=\de,\quad S(\be)=-q^{-1}\be,\quad S(\ga)=-q\ga,
\quad S(\de)=\al.
\endgather
$$

\proclaim{Lemma \theoremname{\lempairingforqM}} For
$p,l\in\Z$, $n,m,r,s\in\Zp$ we have
$\langle A^pB^rC^s, \al^l\be^m\ga^n\rangle = $
$$\de_{rm}\de_{sn}
q^{(pl+l(m+n)+p(m-n))/2}
q^{-m(m-1)/2}q^{-n(n-1)/2} {{(q^2;q^2)_m(q^2;q^2)_n}
\over{(1-q^2)^{m+n}}}.
$$
\endproclaim

As a consequence, we see that the extension of $A_q(SU(2))
\subset \bigl( U_q({\frak m}(2))\bigr)^\ast$, the linear
dual of $U_q({\frak m}(2))$, given by
$$
\sum_{l\in\Z}\sum_{n,m=0}^\infty c_{lmn} \al^l\be^m\ga^n
$$
with $c_{lmn}$ non-zero for only finitely many $l$, is still
well-defined.

Let us now consider a special unitary representation $t^R$ of
$U_q({\frak m}(2))$ acting in $\ell^2(\Z)$ equipped with
an orthonormal basis $\{ e_n\}_{n=-\infty}^\infty$ by
$$
t^R(A)\, e_n =q^n e_n, \quad t^R(B)\, e_n =Re_{n+1},
\quad t^R(C)\, e_n =Re_{n-1},
\quad t^R(D)\, e_n =q^{-n} e_n,
$$
where $R>0$.
Note that the operators for $A$ and $D$ are unbounded. Denote
the corresponding matrix elements by $t^R_{nm}$.

\proclaim{Theorem \theoremname{\thmHEqBesselasmatelts}}
The matrix elements $t^R_{ij}$ are contained in the extension
of $A_q(M(2))$. For
$i\geq j$
$$
t^R_{ij} =\bigl( {{R(1-q^2)}\over{q^{j+1/2}}}\bigr)^{i-j}
{{\al^{i+j}\be^{i-j}}\over{(q^2;q^2)_{i-j}}}
\, {}_1\vp_1\left( {{0}\atop{q^{2(i-j+1)}}};q^2,
-(1-q^2)R^2q^{-2j}\be\ga\right)
$$
and for $i\leq j$
$$
t^R_{ij} =\bigl( {{R(1-q^2)}\over{q^{i+1/2}}}\bigr)^{j-i}
{{\al^{i+j}\ga^{j-i}}\over{(q^2;q^2)_{j-i}}}
\, {}_1\vp_1\left( {{0}\atop{q^{2(j-i+1)}}};q^2,
-(1-q^2)R^2q^{-2i}\be\ga\right).
$$
\endproclaim

\demo{Proof} Use Lemma \thmref{\lempairingforqM}
and the definition of the representation $t^R$ to see
that both sides agree on $A^pB^rC^s$.
\qed\enddemo

The ${}_1\vp_1$-series in Theorem \thmref{\thmHEqBesselasmatelts}
are known as Hahn-Exton $q$-Bessel functions, which are
defined by
$$
J_\nu (z;q) = z^\nu {{(q^{\nu+1};q)_\infty}\over{(q;q)_\infty}}
\, {}_1\vp_1\left( {0\atop{q^{\nu+1}}};q, qz^2\right).
$$
We can now use this interpretation of the Hahn-Exton $q$-Bessel
function as matrix elements of irreducible unitary representations
of $U_q({\frak m}(2))$ to obtain identities for these
$q$-analogues of the Bessel function.

\proclaim{Theorem \theoremname{\thmrelationsforHEqBessel}}
The Hahn-Exton $q$-Bessel function satisfies
the Hansen-Lommel type orthogonality relations;
$$
\sum_{i=-\infty}^\infty q^i J_{k+i}(z;q)\, J_{j+i}(z;q) = 
\de_{kj}q^{-k}, \qquad |z|<q^{-1/2},\ k,j\in\Z,
$$
the Hankel type orthogonality relations;
$$
\sum_{p=-\infty}^\infty q^{2p} J_n(q^{k+p};q^2)\,
J_n(q^{l+p};q^2) = C q^{-2k} \de_{kl},\qquad k,l\in\Z,
$$
for some non-zero constant $C$,
and the Graf type addition formula; for $n,y,x,z\in\Z$, $R>0$
$$
(-q)^n J_{x-n}(q^{z-y};q^2) J_n(Rq^{n+z};q^2) =
\sum_{k=-\infty}^\infty q^{2k} J_k(Rq^{y+x};q^2)
 J_{k-n}(Rq^y;q^2) J_x(q^{z+k-y};q^2),
$$
for $|R^2q^{2x+2y+2}|<1$.
\endproclaim

\demo{Remark \theoremname{\remthmrelationsforHEqBessel}}
The Hahn-Exton $q$-Bessel function of negative integer order
is defined by using
$$
(q^{1-n};q)_\infty \sum_{m=0}^\infty {{c_m }
\over{(q^{1-n};q)_\infty}}
= \sum_{m=n}^\infty c_m (q^{1+m-n};q)_\infty
$$
for $n\in\Zp$. This leads to $J_{-n}(z;q)=(-1)^n q^{n/2}
J_n(zq^{n/2};q)$, which is valid for $n\in\Z$.
\enddemo

\demo{Sketch of Proof} The first relation is a consequence
of the unitarity of the representation $t^R$. The second relation
uses the fact that there exists a Haar functional on
$A_q(M(2))$ as well as on its extension, which is only defined
on a suitable subspace, and for which Schur type orthogonality
relations can be derived. The Hankel type orthogonality
relations then follow from the Schur orthogonality relations.
The constant $C$ is involved, since there is no proper
normalisation for the Haar functional. However, $C=1$, see
Exercise~8.3. Finally, the addition formula follows by
representing the identity
$\De(t^R_{0n})=\sum_{k\in\Z} t^R_{0k}\otimes t^R_{kn}$ in
$A_q(M(2))\otimes A_q(M(2))$ using a suitable representation
of $A_q(M(2))$ in $\ell^2(\Z)$ given by $\al\, e_n =e_{n-1}$,
$\ga\, e_n = q^ne_n$.
\qed\enddemo

\demo{Remark \theoremname{\remthmrelationsforHEqBessel}}
Since the contraction procedure is a limit, we can show that
some of the results of Theorem \thmref{\thmrelationsforHEqBessel}
can be obtained from a suitable limit in the corresponding
result for the little $q$-Jacobi polynomials, which occur
as matrix elements on the quantum $SU(2)$ group, cf.
Corollary \thmref{\cortlnmaslittleqJacobipols}.
\enddemo

\demo{Notes and references} The interpretation of the Hahn-Exton
$q$-Bessel functions on the quantum group of plane motions is
due to Vaksman and Korogodski\u\i\ \cite{\VaksKSMD}. The 
presentation is taken from \cite{\KoelDMJ}. For the 
$C^\ast$-algebra approach to the quantum group of plane motions, 
see Woronowicz \cite{\Worodrie}, and Pal \cite{\Pal} for the 
interpretation of the Hahn-Exton $q$-Bessel functions in this 
context. The addition formula in Theorem 
\thmref{\thmrelationsforHEqBessel} is also derived by Kalnins, 
Miller and Mukherjee \cite{\KalnMMSIAM} using only the quantum 
algebra. See Koornwinder and Swarttouw \cite{\KoorS} for more 
information on the Hahn-Exton $q$-Bessel function, as well as 
for the limit transition from the little $q$-Jacobi polynomials.

For the quantum group of plane motions we can also speak
of the analogue of generalised matrix elements, see
\cite{\KoelIM}, and we can then obtain the analogue of the
addition formula \thetag{\vgladdformqLegtwoparam} for
this situation. Here the Jackson $q$-Bessel functions \cite{\Isma}
play the role of transition coefficients.
\enddemo

%%%%%%%%%%%%%%%%%%%%%%%%%%%%%%%%%%%%%%%%%%%%%%%%%%%%%%%%%%%%%%%%%%%%
%%N E W   S U B S E C T I O N%%%%%%%%%%%%%%%%%%%%%%%%%%%%%%%%%%%%%%%
%%%%%%%%%%%%%%%%%%%%%%%%%%%%%%%%%%%%%%%%%%%%%%%%%%%%%%%%%%%%%%%%%%%%
\subhead\newsubsection 
The quantum algebra approach
\endsubhead
Let us now briefly consider the quantum algebra approach.
The quantum algebra approach is based on the observation that
for $X$ in the Lie algebra $\frak g$ the exponential mapping 
$\exp tX$ gives a function on the corresponding group $G$.
The representation theory
of $U_q{\frak g}$ is usually similar to the representation theory
of ${\frak g}$. Classically we can obtain elements of the 
corresponding group $G$ by exponentiating Lie algebra elements. 
In the quantum algebra approach the action of 
$\exp_q(\al_1X_1)\ldots \exp_q(\al_nX_n)$ is calculated in a
representation of $U_q{\frak g}$. Here $\exp_q$ can be one of the
$q$-analogues of the exponential function, see e.g.
Corollary \thmref{\corstothmqbinomialthm},
$\al_i$ are scalars
and $X_i$ are generators of $U_q{\frak g}$. For a
suitable basis $\{ f_m\}$ of the representation space we get
$$
\exp_q(\al_1X_1)\ldots \exp_q(\al_nX_n)\, f_m =
\sum_k U_{m,k}(\al_1,\ldots,\al_n)\, f_k.
$$
The matrix coefficients $U_{m,k}(\al_1,\ldots,\al_n)$ can
be calculated in terms of special functions in $\al_1,\ldots,\al_n$.

Let us discuss shortly the example of $U_q({\frak m}(2))$. Let $t^R$
be the representation of $U_q({\frak m}(2))$ as in
the previous section. Using the
notation of \S 8.1 and Corollary \thmref{\corstothmqbinomialthm}
we get e.g.
$$
t^R\bigl( e_q(b B) E_q(c C)\bigr)\, e_m =
\sum_{k=-\infty}^\infty U_{m,k}(b,c)\, e_k.
\tag\eqname{\vgldefmateltforHEqBes}
$$
Then it is straightforward to calculate the matrix elements
in terms of $q$-hypergeometric series;
$$
U_{m,k}(b,c) = (Rb)^{k-m} {{(q^{k-m};q)_\infty}\over{(q;q)_\infty}}
\, {}_1\vp_1(0;q^{k-m+1};q,-bcR^2)
$$
which are the Hahn-Exton $q$-Bessel functions. For
$k<m$ we apply Remark \thmref{\remthmrelationsforHEqBessel}.

Having this interpretation, a number of properties
like orthogonality relations and addition formulas
can be derived. Explicit models for the representation space
in terms of function spaces are usually needed. As
a very simple example we derive the Hansen-Lommel
orthogonality relations of Theorem 
\thmref{\thmrelationsforHEqBessel}.
Observe that $e_q(z)E_q(-z)=1$, and hence
$$
e_l = t^R\bigl( e_q(bB) E_q(cC) e_q(-cC) E_q(-bB)\bigr)\, e_l
= \sum_{k=-\infty}^\infty
\Bigl( \sum_{p=-\infty}^\infty T_{l,p}(-b,-c)\, U_{p,k}(b,c)\Bigr)
e_k,
$$
where the matrix coefficient $T_{l,p}(b,c)$ is defined by
$$
t^R\bigl(  E_q(bB) e_q(cC)\bigr) \, e_m =
\sum_{k=-\infty}^\infty T_{m,k}(b,c)\, e_k.
$$
Since interchanging $B$ and $C$ and $A$ and $D$ is a symmetry,
say $J$, of $U_q({\frak m}(2))$ for which
$t^R_{nm}(J(X))= t^R_{mn}(X)$ we find that
$T_{m,k}(b,c)=U_{k,m}(c,b)$. Hence,
$$
\sum_{p=-\infty}^\infty U_{p,l}(-c,-b)\, U_{p,k}(b,c) = \de_{kl}
$$
and this is equivalent to the Hansen-Lommel orthogonality relations
of Theorem \thmref{\thmrelationsforHEqBessel}.

\demo{Notes and references}
The method sketched here
is motivated by the classical relation between Lie algebras and 
special functions as described in Miller's book \cite{\Mill}.
The example is taken from Kalnins, Miller
and Mukherjee \cite{\KalnMMSIAM}.
There exists a huge amount of papers
on this approach, in particular Kalnins, Miller et al.
\cite{\KalnMMJMP}, \cite{\KalnMJMP}, \cite{\KalnMMuJMP}, 
\cite{\KalnMMSIAM} and Floreanini and Vinet \cite{\FlorVPLA}, 
\cite{\FlorVLMP} and references given there, see also the 
references in \cite{\KoelFIC}.

Using the model of $t^R$ in the space of Laurent series with the
actions of $A$, $B$ and $C$ given by $T_q$, the shift operator
defined in Exercise~7.1, $RM_z$ and $RM_{1/z}$, where $M_f$ is
multiplication by $f$. Then the basis $e_m$ corresponds to
$z^m$, $m\in\Z$. In this model \thetag{\vgldefmateltforHEqBes}
is the generating function for the Hahn-Exton $q$-Bessel
function, and a large number of results can already be obtained
by working only with the generating function, see
Koornwinder and Swarttouw \cite{\KoorS} and Koelink \cite{\KoelITSF}
for more general generating functions leading to $q$-Bessel 
functions. The quantum algebra becomes important when dealing with
different models, and in particular when using the tensor
product of representations, see Kalnins, Miller
and Mukherjee \cite{\KalnMMSIAM} for a derivation of the
addition formula of Theorem \thmref{\thmrelationsforHEqBessel}
by using the Clebsch-Gordan decomposition of
$t^R\otimes t^S$.
\enddemo

%%%%%%%%%%%%%%%%%%%%%%%%%%%%%%%%%%%%%%%%%%%%%%%%%%%%%%%%%%%%%%%%%%%%
%%N E W   S U B S E C T I O N%%%%%%%%%%%%%%%%%%%%%%%%%%%%%%%%%%%%%%%
%%%%%%%%%%%%%%%%%%%%%%%%%%%%%%%%%%%%%%%%%%%%%%%%%%%%%%%%%%%%%%%%%%%%
\subhead\newsubsection 
Quantum spheres and spherical functions
\endsubhead
In sections 4 and 5 we discussed the analogue of functions on
$K\backslash SU(2)/H$ for $K$, $H$ one-parameter subgroups,
which are precisely the
elements $b^l_{00}(\ta,\si)$. We can also consider the
elements $b^l_{n0}(\ta,\si)$ as the analogues of the
functions on the sphere $S^2= SU(2)/H$. These
$q$-analogues of the sphere are parametrised by $\si$
and were first introduced by Podle\'s \cite{\Podl},
see Dijkhuizen and Koornwinder \cite{\DijkKGD},
Noumi and Mimachi \cite{\NoumMLNM}.

This can be generalised to give $q$-analogues of $U(n)/U(n-1)$,
the sphere in $\C^n$. There is a quantised universal
enveloping algebra $U_q({\frak{gl}}(n,\C))$, which is in duality
as Hopf algebras with $A_q(GL(n,\C))$. These Hopf algebras can be
made into Hopf $\ast$-algebras, and we obtain a Hopf $\ast$-algebra
$A_q(U(n))$, cf. \cite{\Flor}, \cite{\KoelComp}, \cite{\NoumYM}.
The algebra structure of $A_q(U(n))$ is as follows.
We have $n^2$ generators $t_{ij}$ satisfying the
relations
$$
\gather
t_{ij}t_{il}=qt_{il}t_{ij}, \quad j<l;\qquad
t_{ij}t_{kj}= qt_{kj}t_{ij}, \quad i<k; \\
t_{ij}t_{kl}=t_{kl}t_{ij}, \quad i>k,\, j<l;\qquad
t_{ij}t_{kl}-t_{kl}t_{ij} = (q-q^{-1})t_{il}t_{kj}, 
\quad i<k,\, j<l.
\endgather
$$
Then the $t_{ij}$ generate a bi-algebra, which is an analogue of 
the semigroup of $n\times n$-matrices. In order to obtain a Hopf 
algebra we have to localise along the quantum determinant, which 
is a central element, see \cite{\Flor}, \cite{\KoelComp}, 
\cite{\NoumYM} for detailed relations. The important thing to 
notice is that we have a surjective (Hopf $\ast$-)algebra 
homomorphism $\pi\colon A_q(U(n))\to A_q(U(n-1))$, which is the 
identity on $t_{ij}$, $1\leq i,j<n$ and 
$\pi(t_{ni})=\de_{ni}=\pi(t_{in})$.
On the dual level we have a natural inbedding
$$
U_q({\frak{gl}}(n-1,\C))\oplus U_q({\frak{gl}}(1,\C))
\hookrightarrow U_q({\frak{gl}}(n,\C)),
$$
and we can talk of
$U_q({\frak{gl}}(n-1,\C))\oplus U_q({\frak{gl}}(1,\C))$-invariant
vectors in irreducible representations of $U_q({\frak{gl}}(n,\C))$.
The space of such invariant vectors is at most one-dimensional,
and the corresponding representations can be labeled by two
integers $l$ and $m$. Let $\psi_{l,m}$ be the matrix element
with respect to the invariant vector
of the dual algebra $A_q(U(n))$, then this means precisely
$$
(id\otimes \pi)\De(\psi_{l,m}) = \psi_{l,m}\otimes 1, \qquad
(\pi\otimes id)\De(\psi_{l,m}) = 1\otimes \psi_{l,m}.
$$
So we can view $\psi_{l,m}$ as an element of the deformed
algebra of $U(n-1)$-biinvariant functions
on $U(n)$. Using the Schur orthogonality
for the Haar functional it is possible to derive a very
explicit expression for the zonal spherical elements $\psi_{l,m}$;
$$
\psi_{l,m} = \cases t_{nn}^{l-m} p_m^{(n-2,l-m)}
(1-t_{nn}t_{nn}^\ast;q^2), &\text{if $l\geq m$,} \\
 p_l^{(n-2,m-l)}(1-t_{nn}t_{nn}^\ast;q^2)(t_{nn}^\ast)^{m-l},
&\text{if $m\geq l$,} \endcases
$$
where $p_l^{(n-2,m-l)}$ are little $q$-Jacobi polynomials as in 
\S 5, see Noumi, Yamada and Mimachi \cite{\NoumYM}.
For $n=2$ this corresponds to
Corollary \thmref{\cortlnmaslittleqJacobipols}.

Floris \cite{\Flor} obtains an abstract addition formula, i.e.
in non-commuting variables, for the little $q$-Jacobi polynomials
$p_l^{(n-2,m-l)}(\cdot;q^2)$ by calculating
$\De(\psi_{l,m})$ explicitly
modulo $U_q({\frak{gl}}(n-2,\C))$-invariance in each
factor of the tensor product. Using the representation theory
of the Hopf $\ast$-algebra $A_q(U(n))$, cf. \cite{\KoelComp},
Floris and Koelink \cite{\FlorK} derive an explicit addition
and product formula for the little $q$-Jacobi polynomials
$p_l^{(\al,m-l)}(\cdot;q^2)$, which contains as a special case
the addition formula for the little $q$-Legendre polynomial
$p_l^{(0,0)}(\cdot;q^2)$ discussed in \S 6.

{}From the $n=2$ case discussed in \S\S 4-5, we may suspect that
we can also have $(\ta,\si)$-spherical elements in this case.
This is indeed the case, as shown by Dijkhuizen and Noumi
\cite{\DijkN}, and they obtain an interpretation of
$p_m^{(n-2,0)}(\cdot;q^\ta,q^\si|q^2)$ as spherical functions
on the quantum analogue of $U(n)/U(n-1)$. See also
Dijkhuizen and Koornwinder \cite{\DijkKGD} for a general
discussion of quantum homogeneous spaces.

Instead of considering an analogue of the sphere in
$\C^n$, Sugitani \cite{\Sugi} considers the analogue of
the sphere $SO(n)/SO(n-1)$
in $\R^n$. The zonal spherical functions can be expressed
using big $q$-Jacobi polynomials and continuous
$q$-ultraspherical polynomials, i.e. $p_n^{(\al,\al)}(\cdot;1,1|q)$
with the notation for the Askey-Wilson polynomials in \S 5.

%%%%%%%%%%%%%%%%%%%%%%%%%%%%%%%%%%%%%%%%%%%%%%%%%%%%%%%%%%%%%%%%%%%%
%%N E W   S U B S E C T I O N%%%%%%%%%%%%%%%%%%%%%%%%%%%%%%%%%%%%%%%
%%%%%%%%%%%%%%%%%%%%%%%%%%%%%%%%%%%%%%%%%%%%%%%%%%%%%%%%%%%%%%%%%%%%
\subhead\newsubsection 
Multi-variable orthogonal polynomials as spherical functions
\endsubhead
The multi-variable orthogonal polynomials of importance are the
Macdonald polynomials \cite{\Macdpp}, \cite{\MacdSB}, which are
associated with root systems. The polynomials
depend on $n$ variables and are invariant under the
Weyl group, see Macdonald \cite{\Macdbook, Ch.~VI} for the case
of symmetric functions, i.e. for the root system of type $A$.
Koornwinder \cite{\KoorBC} has obtained multivariable
analogues of the Askey-Wilson polynomials by considering the
non-reduced root system $BC_n$. The orthogonality measures for
these polynomials are absolutely continuous.
For the special case
$n=1$ we obtain the Askey-Wilson polynomials, whereas the case 
$n=1$, i.e. for root system $A_1$, of the Macdonald polynomials 
gives the continuous $q$-ultraspherical polynomials.

Noumi \cite{\NoumAM} shows that for appropriate analogues
of $GL(n)/SO(n)$ and $GL(2n)/Sp(n)$, where $Sp(n)$
is the symplectic group, the Macdonald polynomials for root
system $A_{n-1}$ and a specific choice of the free
parameter $t$ arise as the zonal spherical functions.
Noumi, Dijkhuizen and Sugitani \cite{\NoumDS} have announced
that they can parametrise the quantum homogeneous spaces
continuously, similarly as for the quantum $U(n)/U(n-1)$ space, and
the corresponding zonal spherical functions are expressible
in terms of Koornwinder's \cite{\KoorBC} multivariable
Askey-Wilson polynomials. It is expected that multivariable
orthogonal polynomials with (partly) discrete orthogonality
measure introduced by Stokman \cite{\StokSIAM},
\cite{\Stokrep} can be obtained as limit cases from the general
setting, see Stokman and Koornwinder \cite{\StokK} for the
analytic proof of this limit transition.

%%%%%%%%%%%%%%%%%%%%%%%%%%%%%%%%%%%%%%%%%%%%%%%%%%%%%%%%%%%%%%%%%%%%
%%N E W   S U B S E C T I O N%%%%%%%%%%%%%%%%%%%%%%%%%%%%%%%%%%%%%%%
%%%%%%%%%%%%%%%%%%%%%%%%%%%%%%%%%%%%%%%%%%%%%%%%%%%%%%%%%%%%%%%%%%%%
\subhead\newsubsection 
Some open problems
\endsubhead

\demo{Problem \theoremname{\probsltwoR}}
The quantum $SU(2)$ group and the quantum $SU(1,1)$ group
have been treated in some detail. The other real form,
the quantum $SL(2,\R)$ group, cf. Theorems
\thmref{\thmstarstructureonU} and  \thmref{\thmstarstructureonA},
presents us with some difficulties, although its
representation theory is known, cf. Schm\"udgen \cite{\SchmCMP},
\cite{\Schmrep}. Is it possible to associate special functions
of $q$-hypergeometric type for $|q|=1$ to this quantum group?
In particular, is there a relation with sieved orthogonal
polynomials?
\enddemo

\demo{Problem \theoremname{\probmomentHaar}}
Is it possible to give a proof of Theorem 
\thmref{\thmHaarongeneralspherelt}
along the lines of the proof of Lemma \thmref{\lemHaaronbasis}?
Or, can we derive a simple recurrence relation for
$h(p_n(\rts))$ for some suitable choosen polynomial $p_n$, and
next identify it with the corresponding Askey-Wilson integral?
\enddemo

\demo{Problem \theoremname{\probaddform}}
Derive an addition formula for Askey-Wilson polynomials
by first deriving an abstract addition formula using the
interpretation as zonal spherical function on
the quantum analogue of $U(n)/U(n-1)$, cf. \S 8.3.
Then use the representation theory of $A_q(U(n))$ to derive
an addition formula in commuting variables.
Classically, i.e. for $q=1$ and working in the group case, an
addition formula for the Jacobi
polynomials $R_n^{(\al,0)}$ (notation of Exercise~6.1)
is derived in this way and from this an addition formula
for general Jacobi polynomials $R_n^{(\al,\be)}$
can be obtained by differentiation, cf. \cite{\Koorold}.
For which cases of $\si$ and $\ta$ can we obtain a general
addition formula, i.e. an addition formula for
$p_n^{(\al,\be)}(\cdot,s,t|q)$ for all $\al$, $\be$?
The Rahman-Verma \cite{\RahmV} addition formula suggests that
it might be possible for continuous $q$-Jacobi polynomials, i.e.
$s=t=1$.
\enddemo

%%%%%%%%%%%%%%%%%%%%%%%%%%%%%%%%%%%%%%%%%%%%%%%%%%%%%%%%%%%%%%%%%%%%
%%N E W   S U B S E C T I O N%%%%%%%%%%%%%%%%%%%%%%%%%%%%%%%%%%%%%%%
%%%%%%%%%%%%%%%%%%%%%%%%%%%%%%%%%%%%%%%%%%%%%%%%%%%%%%%%%%%%%%%%%%%%
%EXERCISES%%%%%%%%%%%%%%%%%%%%%%%%%%%%%%%%%%%%%%%%%%%%%%%%%%%%%%%%%%
\subhead Exercises
\endsubhead

\item{\the\sectionno.1} Prove Lemma \thmref{\lempairingforqM}.
Use Exercise~2.5, or apply the
contraction procedure
to Theorem \thmref{\thmexplicitdualityonbasis}.

\item{\the\sectionno.2} Use the $q$-gamma function as in
Exercise~3.9 to see that $J_\nu((1-q)z;q)$ tends
to $J_\nu(2z)$ as $q\uparrow 1$. The Bessel function is
defined by $J_\nu(z) = \sum_{k=0}^\infty (-1)^k z^{\nu+2k}/
(k! \, \Gamma(\nu+k+1))$.

\item{\the\sectionno.3} Prove
$(w;q)_\infty {}_1\vp_1(0;w;q,z) =(z;q)_\infty {}_1\vp_1(0;z;q,w)$
and use this symmetry to derive the Hankel type
orthogonality relations of Theorem 
\thmref{\thmrelationsforHEqBessel} from the Hansen-Lommel type 
orthogonality relations. Show that $C=1$.

\item{\the\sectionno.4} Define matrix coefficients by
$$
t^R\bigl( E_q(b B) E_q(c C)\bigr)\, e_m =
\sum_{k=-\infty}^\infty U_{m,k}(b,c)\, e_k
$$
and
$$
t^R\bigl( e_q(b B) e_q(c C)\bigr)\, e_m =
\sum_{k=-\infty}^\infty T_{m,k}(b,c)\, e_k.
$$
Calculate the matrix coefficients explicitly, and derive the
corresponding Hansen-Lommel (bi-)orthogonality relations.
The $q$-Bessel functions are known as Jackson's $q$-Bessel
functions, see Ismail \cite{\Isma}.

%%%%%%%%%%%%%%%%%%%%%%%%%%%%%%%%%%%%%%%%%%%%%%%%%%%%%%%%%%%%%%%%%%%%
%%R E F E R E N C E S%%%%%%%%%%%%%%%%%%%%%%%%%%%%%%%%%%%%%%%%%%%%%%%
%%%%%%%%%%%%%%%%%%%%%%%%%%%%%%%%%%%%%%%%%%%%%%%%%%%%%%%%%%%%%%%%%%%%

\enddocument